\documentclass[12pt, reqno]{amsart}
\usepackage{amsmath, amsthm, amscd, amsfonts, amssymb, graphicx, color, mathrsfs}
\usepackage[bookmarksnumbered, colorlinks, plainpages]{hyperref}
\usepackage[all]{xy} 
\usepackage{slashed}
\usepackage{cite}
\usepackage{xcolor}
\usepackage{listings}
\usepackage{xcolor}
\usepackage{comment}
\usepackage{pgfplots}
\usepgfplotslibrary{patchplots}
\pgfplotsset{compat=1.18}
\usepackage{subcaption}

\renewcommand{\Im}{\mathrm{Im}\,}
\renewcommand{\Re}{\mathrm{Re}\,}

\newtheorem{theorem}{Theorem}[section]
\newtheorem{lemma}[theorem]{Lemma}
\newtheorem{proposition}[theorem]{Proposition}
\newtheorem{corollary}[theorem]{Corollary}
\theoremstyle{definition}
\newtheorem{definition}[theorem]{Definition}

\theoremstyle{remark}
\newtheorem{remark}[theorem]{Remark}
\numberwithin{equation}{section}

\newcommand*\diff{\mathop{}\!\mathrm{d}}

\begin{document}
\setcounter{page}{1}

\color{darkgray}{
\noindent \centering
{\small   }\hfill    {\small }\\
{\small }\hfill  {\small }}

\centerline{}

\centerline{}


\title[Bilinear Fourier integral operators with complex phases]{Boundedness properties of bilinear Fourier integral operators with complex phases}

\author[D. Cardona]{Duv\'an Cardona S\'anchez$^{1,**}$}
\address{
  Duv\'an Cardona:
  \endgraf
   King Fahd University of Petroleum and Minerals 
   \endgraf
    College of Computing and Mathematics
  \endgraf $\&$
   Department of  Mathematics
  \endgraf 
   IR Center for Industrial Nuclear Energy, 
  \endgraf
   Dhahran,  Saudi Arabia.
  \endgraf
  Pontificia Universidad Javeriana, 
  \endgraf
 Department of Mathematics
 \endgraf
  Bogot\'a-Colombia.
  \endgraf
  {\it E-mail address:} {\rm duvanc306@gmail.com, duvan.sanchez@kfupm.edu.sa}   }

\author[M. MART\'INEZ-Flores]{ MANUEL  MART\'INEZ-Flores$^{2,*}$}
\address{
  Manuel  Mart\'inez-Flores
  \endgraf
  Department of Mathematics
  \endgraf
  Bonn University, 
  \endgraf
  Bonn, Germany.
  \endgraf
  {\it E-mail address:} {\rm  manuelalejandromartinezf@gmail.com } {\rm manuel.martinezflores@uni-bonn.de}
  }

\author[F. Opoku]{Frederick Opoku$^{3,*}$}
\address{
  Frederick Opoku
  \endgraf
  Department of Mathematics
  \endgraf
  Kwame Nkrumah University of Science and Technology
  \endgraf
  (KNUST)- Ghana.
  \endgraf
  {\it E-mail address:} {\rm frederick@aims.edu.gh}
  }

\thanks{Duv\'an Cardona S\'anchez  has been  supported  by the Department of Mathematics at the King Fahd University of Petroleum and Minerals and by the IR Center of Industrial Nuclear Energy. He also was supported in part by the FWO  Odysseus  1  grant  G.0H94.18N:  Analysis  and  Partial Differential Equations and by the Methusalem programme of the Ghent University Special Research Fund (BOF)
(Grant number 01M01021), by the FWO Fellowship
Grant No 1204824N and by the FWO Grant K183725N of the Belgian Research Foundation FWO. He also has been supported by the Oberwolfach Leibniz Fellow of the Mathematical Institute of Research of Oberwolfach, MFO-Germany, Project F2511 (2026) and by the Department of Mathematics of the Pontificia Universidad Javeriana, Bogot\'a-Colombia. 
This work was partially supported by the Friends of IMU (FIMU), under the IMU Breakout Fellowship Program. Frederick Opoku has been supported by the Breakout Graduate Fellowship of the International Mathematical Union.
He has also been supported by the Kwame Nkrumah University of Science and Technology (KNUST), Ghana. The authors have been benefited by the support of the Research Group: Analysis and PDE in Developing Countries of ICMAM Latin America (International Digital Center of Mathematical Research for Latin America).}

\begin{abstract} In this paper, we investigate mapping properties on Lebesgue spaces of bilinear Fourier integral operators (FIOs) with complex phases. Indeed, we provide for bilinear Fourier integral operators with complex phases boundedness theorems for the corresponding  Seeger-Sogge-Stein orders, including the endpoint estimates. Our results extend the same boundedness property in the case of real-valued phases proved by Rodriguez-L\'opez,  Rule, and Staubach, and provide the bilinear version for the linear one proved earlier by Ruzhansky for phases of positive type. Also, our Seeger-Sogge-Stein theorem for bilinear operators with complex phases cannot be directly deduced from its real counterpart.
A new approach used in this article consists of providing directly the boundedness theorems on Lebesgue spaces without the use of the multilinear Fefferman-Stein interpolation theorem. For this, a maximal function theorem due to Bourgain is applied. 
As a byproduct, we also consider transference boundedness theorems of FIOs with complex phases from $\mathbb{R}^n$ to the torus $\mathbb{T}^n\equiv \mathbb{R}^n/\mathbb{Z}^n$ and then, we provide new criteria for the boundedness of  Fourier integral operators on the torus with complex phases (Fourier series operators in the sense of Ruzhansky and Turunen).
\end{abstract} \maketitle
\allowdisplaybreaks
\tableofcontents

\section{Introduction}

\subsection{Overview}
This paper continues the analysis in \cite{Cardona:2026} that was primarily dedicated to obtaining the local smoothing property for bilinear Fourier integral operators and to proposing the corresponding {\it bilinear local smoothing conjecture} for Fourier integrals. One of the novelties of this manuscript is that we use a boundedness theorem for maximal functions due to Bourgain \cite{Bourgain1985}, to provide the boundedness of bilinear Fourier integral operators without using, as usual, the bilinear Fefferman-Stein interpolation method.

To complete the perspective on the subject, this paper addresses the following problems:
\begin{itemize}
    \item providing the analogue of the Seeger-Sogge-Stein theorem \cite{seeger1991regularity} of Rodriguez-L\'opez, Rule, and Staubach \cite{rodriguez2011boundedness} for bilinear Fourier integral operators with complex-valued phase functions, and,
    \item establishing transference theorems from $\mathbb{R}^n$ to the torus $\mathbb{T}^n:=(\mathbb{R}/\mathbb{Z})^n$, to provide boundedness properties for linear and general multilinear Fourier integral operators on the torus in the setting of complex-valued phases. 
\end{itemize} It is important to mention that, when considering the case of phase functions with zero imaginary part, we get analogue results in the literature in the case of $\mathbb{R}^n$ and for bilinear operators, the resultant boundedness criteria are new in the torus $\mathbb{T}^n.$ Particularly, the Seeger-Sogge-Stein theorem for bilinear Fourier integral operators with real-valued phase functions has been established by Rodriguez-Lopez,  Rule and Staubach in \cite{rodriguez2011boundedness}.

There are several reasons why the theory of Fourier integral operators with complex phases is more appropriate than that of their real counterparts. Indeed, from the point of view of symplectic geometry, in the complex case, complex canonical relations can be parametrised by a single phase function. Moreover, this is not, in general, possible for real canonical relations due to topological obstructions related to the Maslov index \cite{melin2006fourier}. However, any real canonical relation admits an almost analytic extension contained in a complex canonical relation, which is then parametrised by a single phase, see Laptev, Safarov, and Vassiliev \cite{Laptev:Safarov:Vassiliev}.

As recently observed by Martini and M\"uller \cite{MMNG2023}, Fourier integral operators with complex phases appear when constructing the parametrix for the wave operator associated to the sub-Laplacian defined on M\'etivier groups. These classes of groups contain, for instance, Heisenberg-type groups.
On the other hand, the theory of FIOs with complex canonical relations was developed by Melin and Sj\"ostrand \cite{melin2006fourier}, motivated by the analysis of parametrices for problems with oblique derivatives \cite{melin1976fourier}. 

In complex geometry, H\"ormander's work on asymptotics for weakly pseudoconvex or hyperconvex domains led to the existence of the complete asymptotic expansion at the boundary by Fefferman \cite{Fefferman1974} on the diagonal. Meanwhile, Boutet de Monvel and Sj\"ostrand \cite{de1976singularite} described the singularity of the full Bergman kernel by showing that it is a Fourier integral operator with complex phases. Recently, FIOs with complex phases have also arisen in the investigation of asymptotics for the Bergman projection kernel, since its full kernel is a Fourier integral operator with complex phases; see Hsiao and Marinescu \cite{hsiao2025singularities}.

The pioneering work of Ruzhansky \cite{Ruzhansky2001} has established the Seeger-Sogge-Stein theorem for Fourier integral operators with complex phases. The recent work by the first author and Ruzhansky \cite{cardona2026weak} has obtained for these operators the corresponding weak (1,1) inequality previously obtained in the setting of real-valued phases by Tao \cite{Tao2004}. 

Part of the motivation for this work lies in providing corresponding estimates for bilinear Fourier integral operators with phases that parametrise complex canonical relations. Application of these estimates includes, for instance, estimates for products $u(x,t)v(x,t)$ of solutions to the wave equation and other hyperbolic partial differential equations.

\subsection{Main results}
The first main result of this paper regarding the boundedness of bilinear Fourier integral operators with complex phases is Theorem \ref{Main:1:th} below. We introduce the required notation as follows. In estimates, we write $A\lesssim B$ as a shorthand for $A \le CB$ for $C>0$ and will use the notation $A \approx B$ if $A\lesssim B$ and $B\lesssim A$. We present the notion of bilinear FIOs investigated in this work.

\begin{definition}[Bilinear Fourier integral operator]\label{Def:1}
    A bilinear Fourier integral operator $ T_a^{\phi_1,\phi_2}:C^\infty_0(\mathbb{R}^n)\times C^\infty_0(\mathbb{R}^n)\rightarrow C^\infty(\mathbb{R}^n),$  is an operator which is defined to act on test functions $f$ and $g$ by the formula
    \begin{align*}
T^{\phi_{1}, \phi_{2}}_{a}(f,g)(x)
&= \int_{\mathbb{R}^n}\int_{\mathbb{R}^n}
e^{2\pi i (\phi_1(x,\xi)+ \phi_2(x,\eta))}
\, a(x,\xi,\eta)\, \widehat{f}(\xi)\, \widehat{g}(\eta)
\, d\xi \, d\eta,
    \end{align*}
    where:
    
    $\bullet$ 
    $\widehat{f}$ and $\widehat{g}$ are the Fourier transform of $f\in C^\infty_0(\mathbb{R}^n)$ and $g\in C^\infty_0(\mathbb{R}^n),$ respectively. 
    
     $\bullet$  The bilinear symbol $a:=a(x,\xi,\eta) \in C^\infty({\mathbb{R}^n \times \mathbb{R}^n \times \mathbb{R}^n})$ is assumed to have compact support in $x\in \mathbb{R}^n.$

      $\bullet$  The phase  $\phi_1:=\phi_1(x,\xi)\in C^\infty( \mathbb{R}^n\times \mathbb{R}^n\setminus\{0\},\mathbb{C})$ is positively homogeneous of degree one in $\xi\neq0$ and  $\phi_2:=\phi_2(x,\eta)\in C^\infty( \mathbb{R}^n\times \mathbb{R}^n\setminus\{0\},\mathbb{C}),$ is positively homogeneous
of degree one in $\eta\neq0$.

      $\bullet$ We assume that $\textnormal{Im}(\phi_1), \textnormal{Im}(\phi_2)\geq 0.$
\end{definition}

To formulate our main Theorem~\ref{Main:1:th} of bilinear Fourier integral operators with complex phases, we consider the following assumption for the phases as in the calculus of FIOs with complex phases developed by Melin and Sj\"ostrand \cite{melin2006fourier,melin1976fourier}.

\begin{definition}\label{MS:C:tau} We say that $\phi_1$ and $\phi_2,$ as in Definition \ref{Def:1} are of {\it Melin-Sj\"ostrand type} if they satisfy the following property:

$\bullet$ assume that for some $\tau\in \mathbb{R},$ we have that $$\phi_{1,\tau}(x,\xi) = \Re \phi_1(x,\xi) + \tau \Im \phi_1(x,\xi), $$
and $$\phi_{2,\tau}(x,\eta) = \Re \phi_2(x,\eta) + \tau \Im \phi_2(x,\eta),$$ satisfy the non-degeneracy condition, meaning this that: 
$$\forall (x,\xi)\in \mathbb{R}^n\times \mathbb{R}^n\setminus\{0\}, \det\left(
\partial_{x,\xi}^{2}\phi_{1,\tau}(x,\xi)
\right)_{1\leq i,j\leq n}\neq 0,$$ and that 
$$\forall (x,\eta)\in \mathbb{R}^n\times \mathbb{R}^n\setminus\{0\}, \det\left(
\partial_{x,\eta}^{2}\phi_{2,\tau}(x,\eta)
\right)_{1\leq i,j\leq n}\neq 0.$$ 
We will also say, equivalently, that $\phi_1$ and $\phi_2$ satisfy the {\it Melin-Sj\"ostrand condition} for some $\tau \in \mathbb{R}.$ In addition, we say that these phases are of {\it strictly positive type}, if 
there exists
$c>0$ such that
$$
\operatorname{Im}\phi_{j}(x,\zeta)
\geq c|\zeta|,
\qquad
x\in\mathcal K,
\quad
\zeta\neq0.
$$
\end{definition}
\begin{remark} Note that if   $\phi_1$ and $\phi_2$ satisfy the {\it Melin-Sj\"ostrand condition} for some $\tau \in \mathbb{R},$ then both phases also satisfy the same condition for almost all $\tau\in \mathbb{C},$ see Remark  \ref{Rem:MS:condition}.
    Note also that the homogeneity condition on $\phi_{1,\tau}$ and on $\phi_{2,\tau}$ holds under the fact that the same condition holds for $\phi_{1}$ and $\phi_{2}.$ Note also that standard real-valued non-degenerate phase functions as considered in Seeger, Sogge, and Stein \cite{seeger1991regularity} hold the {\it Melin-Sj\"ostrand condition} with $\tau=0.$
\end{remark}
\begin{definition}\label{Def:bili:m}
The Kohn--Nirenberg bilinear class $$S^{m}_{1,0}:=S^{m}_{1,0}\bigl(\mathbb{R}^{n}_{x}\times\mathbb{R}^{n}_{\xi}\times\mathbb{R}^{n}_{\eta}\bigr)$$  is the set formed by those bilinear symbols $a \in C^\infty({\mathbb{R}^n \times \mathbb{R}^n \times \mathbb{R}^n})$ satisfying the symbol inequalities
\begin{align*}
\left|
\partial_{\xi}^{\alpha}
\partial_{\eta}^{\beta}
\partial_{x}^{\gamma}
a(x,\xi,\eta)
\right|
\lesssim_{\alpha,\beta,\gamma}
( 1+ |\xi|+|\eta|)^{m-|\alpha|-|\beta| }.
\end{align*} Also, a bilinear Fourier integral operator as in Definition \ref{Def:1} is of order $m\in \mathbb{R},$ if it is associated to a bilinear symbol $a\in S^{m}_{1,0}.$    
\end{definition}

The following is our main theorem. 

\begin{theorem}[Bilinear Seeger--Sogge--Stein theorem for FIOs with complex phases]\label{Main:1:th}
Let $
T^{\phi_{1}, \phi_{2}}_{a} : C^\infty_0(\mathbb{R}^n)\times C^\infty_0(\mathbb{R}^n)\rightarrow C^\infty(\mathbb{R}^n),$ 
be a bilinear Fourier integral operator of order $m\in \mathbb{R}$ in the sense of Definition \ref{Def:bili:m}. Assume that simultaneously $\phi_1$ and $\phi_2$ satisfy the Melin-Sj\"ostrand condition for some $\tau\in \mathbb{R},$ and are phases of {\it strictly positive type}.
    Then, when $$1< p_1,p_2< \infty,\quad 1<p < \infty, \quad\frac{1}{p}=\frac{1}{p_1}+\frac{1}{p_2}$$ and
\begin{align}\label{SSS:order:bilinear}
m\le-(n-1)\left(
\left|\frac{1}{p_1}-\frac{1}{2}\right|
+\left|\frac{1}{p_2}-\frac{1}{2}\right|
\right), \quad n\ge 2,
\end{align} the bilinear Fourier integral operator $T^{\phi_{1}, \phi_{2}}_{a}
:L^{p_1}(\mathbb{R}^n)\times L^{p_2}(\mathbb{R}^n)
\longrightarrow L^p(\mathbb{R}^n)$ extends to a bounded bilinear operator. Equivalently, $$ \|T^{\phi_{1},\phi_{2}}_{a}(f,g)\|_{L^p}
\lesssim \|f\|_{L^{p_1}}\|g\|_{L^{p_2}}.$$
\end{theorem}
\begin{remark}
    We observe that Theorem \ref{Main:1:th} cannot be derived from its real-counterpart, see Remark \ref{Rem:fund} and Figure \ref{fig:order-comparison}.
\end{remark}
\pgfplotsset{
    colormap={blues}{
        rgb255(0cm)=(0, 34, 78)      
        rgb255(4cm)=(200, 210, 255) 
    }}
\pgfplotsset{
    colormap={oranges}{
        rgb255(0cm)=(100, 50, 3)      
        rgb255(4cm)=(254, 227, 145)   
    }
}
\definecolor{kiteLightBlue}{RGB}{180, 190, 240}      
\definecolor{kiteLightOrange}{RGB}{240, 210, 130}

\begin{tikzpicture}
    \begin{axis}[
        view={55}{45}, 
        width=0.4\textwidth, 
        tick label style={font=\tiny},
        xlabel={$\frac{1}{p_1}$},
        ylabel={$\frac{1}{p_2}$},
        xtick={0, 0.5, 1},
        ytick={0, 0.5, 1},
        ztick={-2.5, -1.25, 0},
        yticklabels={0, $\frac{1}{2}$, 1},
        xticklabels={0, $\frac{1}{2}$, 1},
        grid=major,
        colormap name=blues, 
        scale only axis=true,
        faceted color=white,
        xlabel shift=-10pt,
        ylabel shift=-10pt,
        enlargelimits=false
    ]

    \addplot3[
        patch,
        patch type=bilinear,
        shader=flat,  
        point meta=explicit,
        point meta min=-2.2,
    ] coordinates {
        (0.5, 0.5,  0.0)[0] (1.0, 0.5, -1.25)[0] (1.0, 1.0, -2.5)[0] (0.5, 1.0, -1.25)[0]
        (0.0, 0.5, -1.25)[-0.5] (0.5, 0.5,  0.0)[-0.5] (0.5, 1.0, -1.25)[-0.5] (0.0, 1.0, -2.5)[-0.5]
        (0.0, 0.0, -2.5)[-0.7] (0.5, 0.0, -1.25)[-0.7] (0.5, 0.5,  0.0)[-0.7] (0.0, 0.5, -1.25)[-0.7]
        (0.5, 0.0, -1.25)[-0.2] (1.0, 0.0, -2.5)[-0.2] (1.0, 0.5, -1.25)[-0.2] (0.5, 0.5,  0.0)[-0.2]
    };
    \addplot3[
        color=lightgray,
    ] coordinates {
        (1, 0.5, -2.5)
        (1, 0.5, -1.25)
    };
    \addplot3[
        color=lightgray,
    ] coordinates {
        (0.5, 0., -2.5)
        (0.5, 0., -1.25)
    };
    \addplot3[
        draw=black,
        thick,
        mark=none
    ] coordinates {
        (0.5, 0.9, -1) 
        (0.9, 0.5, -1) 
        (0.5, 0.1, -1) 
        (0.1, 0.5, -1) 
    };
    \end{axis}
\end{tikzpicture}
\begin{tikzpicture}
    \begin{axis}[
        view={55}{45}, 
        width=0.4\textwidth, 
        tick label style={font=\tiny},
        xlabel={$\frac{1}{p_1}$},
        ylabel={$\frac{1}{p_2}$},
        xtick={0, 0.5, 1},
        ytick={0, 0.5, 1},
        ztick={-2.5, -2, -1, 0},
        yticklabels={0, $\frac{1}{2}$, 1},
        xticklabels={0, $\frac{1}{2}$, 1},
        grid=major,
        colormap name=oranges, 
        scale only axis=true,
        faceted color=white,
        xlabel shift=-10pt,
        ylabel shift=-10pt,
        enlargelimits=false
    ]

    \addplot3[
        patch,
        patch type=bilinear,
        shader=flat,  
        point meta=explicit,
        point meta min=-2.2,
    ] coordinates {
        (0.5, 0.5,  0.0)[0] (1.0, 0.5, -1)[0] (1.0, 1.0, -2)[0] (0.5, 1.0, -1)[0]
        (0.0, 0.5, -1)[-0.5] (0.5, 0.5,  0.0)[-0.5] (0.5, 1.0, -1)[-0.5] (0.0, 1.0, -2)[-0.5]
        (0.0, 0.0, -2)[-0.7] (0.5, 0.0, -1)[-0.7] (0.5, 0.5,  0.0)[-0.7] (0.0, 0.5, -1)[-0.7]
        (0.5, 0.0, -1)[-0.2] (1.0, 0.0, -2)[-0.2] (1.0, 0.5, -1)[-0.2] (0.5, 0.5,  0.0)[-0.2]
    };
    \addplot3[
        color=lightgray,
    ] coordinates {
        (1, 0.5, -2.5)
        (1, 0.5, -1)
    };
    \addplot3[
        color=lightgray,
    ] coordinates {
        (0.5, 0., -2.5)
        (0.5, 0., -1)
    };
    \addplot3[
        color=lightgray,
    ] coordinates {
        (1, 0., -2.5)
        (1, 0., -2)
    };
    \addplot3[
        draw=black, 
        thick,
        mark=none
    ] coordinates {
        (0.5, 1.0, -1) 
        (1.0, 0.5, -1) 
        (0.5, 0.0, -1) 
        (0.0, 0.5, -1) 
    };
    \end{axis}
\end{tikzpicture}
\begin{figure}[htbp]
    \centering
    \include{order-comparison.tex}
    \caption{Comparison of the order  \eqref{nonsharp} provided by the theory of FIOs with real-valued phases for the operators of this manuscript (left in blue) and the Seeger-Sogge-Stein order \eqref{SSS:order:bilinear} (right in orange) reached in our analysis for $n=3$. }
    \label{fig:order-comparison}
\end{figure}

\begin{remark}
The assumption $n\ge 2$ is necessary.
For $n=1$, some of these estimates for Fourier integral operators fail, see the analysis in \cite{rodriguez2011boundedness}.
\end{remark}
\begin{remark}\label{Rem:MS:condition}
    Note also that if in Theorem \ref{Main:1:th}, the phases $\phi_1,\phi_2$ satisfy the Melin-Sj\"ostrand condition for some $\tau\in \mathbb{R},$ that is, each phase $\phi_{j,\tau}(x,\xi)$ in Definition \ref{MS:C:tau} is non-degenerate for some $\tau\in \mathbb{R},$ the same conditions hold for almost every $\tau\in \mathbb{C}.$ In consequence, if the phases $\phi_1,\phi_2$ satisfy the Melin-Sj\"ostrand condition for some $\tau\in \mathbb{R},$ they also satisfy the same condition for almost every $\tau\in \mathbb{C}.$ 
\end{remark}
\begin{remark}\label{Rem:fund}
Our Seeger--Sogge--Stein theorem for bilinear Fourier integral operators with complex phases does not follow directly from the real-phase counterpart established in \cite{rodriguez2011boundedness}. In fact, a bilinear Fourier integral operator with complex phase and symbol in the class $S^{m}_{1,0}:=S^{m}_{1,0}\bigl(\mathbb{R}^{n}_{x}\times\mathbb{R}^{n}_{\xi}\times\mathbb{R}^{n}_{\eta}\bigr)$ can be realized as a bilinear Fourier integral operator with real phase but with symbol in the less regular class $S^{m}_{1/2,1/2}:=S^{m}_{1/2,1/2}\bigl(\mathbb{R}^{n}_{x}\times\mathbb{R}^{n}_{\xi}\times\mathbb{R}^{n}_{\eta}\bigr)$. For such operators, the conjectured order ensuring boundedness from $L^{p_1}(\mathbb{R}^n)\times L^{p_2}(\mathbb{R}^n)$ into $L^p(\mathbb{R}^n)$ is
\begin{equation}\label{nonsharp}
    m \le -(n-\tfrac12)\left( \left|\frac{1}{p_1}-\frac{1}{2}\right| + \left|\frac{1}{p_2}-\frac{1}{2}\right| \right), \quad n\ge 2,
\end{equation}
which, however, is not sharp when compared with the result obtained in the present paper, namely, the one stated in \eqref{SSS:order:bilinear}. See Figures \ref{fig:order-comparison} and \ref{fig:order-comparison2}.
\end{remark}
{\centering
\begin{tikzpicture}
\begin{axis}[
    width=0.6\textwidth,
    xlabel={$\frac{1}{p_1}$},
    ylabel={$\frac{1}{p_2}$},
    xmin=0, xmax=1,
    ymin=0, ymax=1,
    xtick={0, 0.1, 0.5, 0.9, 1}, 
    ytick={0, 0.1, 0.5, 0.9, 1},
    xticklabels={0, $\frac{1}{10}$, $\frac{1}{2}$, $\frac{9}{10}$, 1},
    yticklabels={0, $\frac{1}{10}$, $\frac{1}{2}$, $\frac{9}{10}$, 1},
    grid=major,
    enlargelimits=false,
    axis equal image 
]
\fill[kiteLightOrange, fill opacity=0.8, even odd rule]
    (0.5, 1.0) -- (1.0, 0.5) -- (0.5, 0.0) -- (0.0, 0.5) -- cycle
    (0.5, 0.9) -- (0.1, 0.5) -- (0.5, 0.1) -- (0.9, 0.5) -- cycle;
\addplot[
    fill=kiteLightBlue,
    draw=black,
    opacity=0.8
] coordinates {
    (0.5, 0.9) (0.9, 0.5) (0.5, 0.1) (0.1, 0.5) (0.5, 0.9)
};
\end{axis}
\end{tikzpicture}
\begin{figure}[htbp]
    \centering
    \include{order-comparison2.tex}
    \caption{Contour plot of the non-sharp order \eqref{nonsharp} (in blue) and the Seeger-Sogge-Stein order \eqref{SSS:order:bilinear} (in orange) for the admissible values of $p_1$ and $p_2$ for $m=-1, n=3,$ in each regime. }
    \label{fig:order-comparison2}
\end{figure}

}

\begin{remark}
We observe that the case where $\operatorname{Im}\phi \equiv 0$ in Theorem \ref{Main:1:th},  has been proved by Rodriguez-L\'opez, Rule, and Staubach in \cite{rodriguez2011boundedness}. We also recall that Ruzhansky \cite{Ruzhansky2001} has established the linear version of the Seeger-Sogge-Stein theorem for Fourier integral operators with complex phases.
\end{remark}
\begin{remark} We have assumed that  simultaneously $\phi_1$ and $\phi_2$ satisfy the Melin-Sj\"ostrand condition for some $\tau\in \mathbb{R},$ and that are phases of {\it strictly positive type}. The latter condition is of qualitative nature, since we require that 
$$
\operatorname{Im}\phi_{j}(x,\zeta)
\geq c|\zeta|,
\qquad
x\in\mathcal K,
\quad
\zeta\neq0.
$$ Particularly, these kind of conditions become important to establish other kind of estimates, as the weak (1,1) property for linear FIOs in the setting of complex-valued functions, see 
    \cite{cardona2026weak}. We also refer to Tao \cite{Tao2004} for the weak (1,1)-estimate of Fourier integral operators associated to real-valued phases. 
\end{remark}

In this paper, we prove the following transference criterion, which will allow us to establish the Seeger-Sogge-Stein $L^p$-boundedness theorem for bilinear operators on the torus.
 \begin{theorem}
    Let $1<p<\infty$, and let $\phi$ a complex valued function defined in $\mathbb{T}^n\times \mathbb{R}^{nr}$ such that $\mathrm{Im}(\phi) \geq M$ where $M\in \mathbb{R}.$ Let us assume that $a:\mathbb{T}^n\times \mathbb{R}^{nr}\to \mathbb{C}$ is bounded and continuous, and the corresponding multi-linear Fourier integral operator
    \begin{equation*}
        Tf(x) = \int_{\mathbb{R}^{nr}} e^{i\phi(x, \xi)}a(x, \xi)\widehat{f_1}(\xi_1)\cdots \widehat{f_r}(\xi_r) \diff \xi
    \end{equation*}
    extends to a bounded linear operator from $L^{p_1}(\mathbb{R}^n)\times\cdots\times L^{p_r}(\mathbb{R}^n)$ into $L^{p}(\mathbb{R}^n)$. Then the periodic multi-linear Fourier integral operator 
    \begin{equation*}
        Af(x) = \sum_{\xi \in\mathbb{Z}^{nr}} e^{i\phi(x, \xi)}a(x, \xi)(\mathcal{F}_{\mathbb{T}^n}f_1)(\xi_1)\cdots (\mathcal{F}_{\mathbb{T}^n}f_r)(\xi_r) 
    \end{equation*}
    extends to a bounded linear operator from $L^{p_1}(\mathbb{T}^n)\times\cdots\times L^{p_r}(\mathbb{T}^n)$ into $L^{p}(\mathbb{T}^n)$, given that
    \begin{equation*}
        \frac{1}{p_1} + \cdots +\frac{1}{p_r} = \frac{1}{p}, \quad 1\leq p_i <\infty.
    \end{equation*}
    Moreover, there exists a positive constant $C_p$ such that 
    \begin{equation*}
        \|A\|_{\mathcal{B}(L^{p_1}(\mathbb{T}^n)\times\cdots\times L^{p_r}(\mathbb{T}^n), L^{p}(\mathbb{T}^n))} \leq C_p \|T\|_{\mathcal{B}(L^{p_1}(\mathbb{R}^n)\times\cdots\times L^{p_r}(\mathbb{R}^n), L^{p}(\mathbb{R}^n))}.
    \end{equation*}
\end{theorem}
Combining this result with the one about bilinear Fourier integral operators in $\mathbb{R}^n$ (see Theorem~\ref{Main:1:th}) allows us to obtain the analogue result in the torus, $\mathbb{T}^n.$ 
\begin{corollary}\label{cor:torus:FIOs}
    Let $T^{\phi_1,\phi_2}_a$ be a bilinear Fourier integral operator defined on the torus $\mathbb{T}^n\cong [0,1)^n,$ as 
    \begin{equation*}
        T^{\phi_1,\phi_2}_a (f,g)(x) = \sum_{\xi,\eta\in\mathbb{Z}^n} e^{2\pi i[\phi_1(x, \xi)+\phi_2(x,\eta)]} a(x, \xi,\eta) (\mathcal{F}_{\mathbb{T}^n}f)(\xi)(\mathcal{F}_{\mathbb{T}^n}g)(\eta),
    \end{equation*}
    where $a\in S^m_{1,0}(\mathbb{T}^n_x \times \mathbb{R}^n_{\xi}\times \mathbb{R}^n_{\eta})$ is compactly supported in $x,$
    and there exists $\tau\in \mathbb{R},$ such that $\phi_{1,\tau}(x,\xi) = \Re \phi_1(x,\xi) + \tau \Im \phi_1(x,\xi) $ and $\phi_{2,\tau}(x,\eta) = \Re \phi_2(x,\eta) + \tau \Im \phi_2(x,\eta)$ satisfy:
\begin{itemize}
        \item[1.] $\phi_{1,\tau}(x,\xi)$ is positively homogeneous of degree one in $\xi\neq 0,$ and $\phi_{2,\tau}(x,\eta)$ is also positively homogeneous of degree one in $\eta\neq 0,$
        \item[2.]  $\phi_{1,\tau}(x,\xi) $ and $\phi_{2,\tau}(x,\eta)$ satisfy the non-degeneracy condition.
        \item[3.] $\phi_1, \phi_2\in C^\infty( \mathbb{T}^n\times \mathbb{R}^n\setminus\{0\},\mathbb{C}).$
        \item[4.] We assume that the phases $\phi_j,$ with $\textnormal{Im}(\phi_1), \textnormal{Im}(\phi_2)\geq 0,$ are of {\it strictly positive type}.
    \end{itemize} 
    Let us consider
    \begin{equation*}
        1< p_1,p_2 < \infty, \quad 1<p< \infty, \quad \frac{1}{p}=\frac{1}{p_1} + \frac{1}{p_2},
    \end{equation*}
    and 
    \begin{equation*}
        m \le - (n-1)\left( \left| \frac{1}{p_1}-\frac{1}{2} \right| + \left| \frac{1}{p_2}-\frac{1}{2} \right|
        \right), \quad n\ge 2.
    \end{equation*}
    Then $T^{\phi_1,\phi_2}_a$ extends to a bounded bilinear operator from $L^{p_1}(\mathbb{T}^n) \times L^{p_2}(\mathbb{T}^n)$ to $L^{p}(\mathbb{T}^n)$. Moreover, we have that 
    \begin{equation*}
        \| T^{\phi_1,\phi_2}_a(f,g)\|_{L^p} \lesssim \|f\|_{L^{p_1}}\|g\|_{L^{p_2}}.
    \end{equation*}
\end{corollary}

Now we briefly discuss our results from the perspective of other results in the literature. 
\begin{remark}
    The theory of pseudo-differential operators and of Fourier series operators on the torus \(\mathbb{T}^n\) and, more generally, on compact manifolds without boundary has attracted considerable attention in recent decades,  see e.g. the historical perspective provided in the monograph \cite{Cardona:Ruzhansky:2026}. 
    
    Foundational results on the boundedness of periodic pseudo-differential operators were established in the monograph of Ruzhansky and Turunen \cite{ruzhansky-turunen}. Then, other boundedness results are available in \cite{Cardona:2014, Cardona:2016, Cardona:2018}, where weak type (1,1) estimates, H\"older--Besov boundedness, and general \(L^p\)-criteria were systematically developed for symbol classes on the torus. These toroidal results were further refined and extended in a comprehensive series of works \cite{Cardona:Martinez:2025a, Cardona:Martinez:2025b, Cardona:Martinez:2025c}, which revisited the boundedness theory on the torus from multiple perspectives, providing sharp conditions and unifying previously known criteria. 
    
    The multilinear setting on the torus and on \(\mathbb{Z}^n\) was investigated in \cite{Cardona:Kumar:2019}, where \(L^p\)-boundedness and nuclearity properties of multilinear pseudo-differential operators were established. Beyond the pseudo-differential framework, periodic Fourier integral operators (FIOs) have been studied in \cite{Cardona:Messiouene:Senoussaoui:2021}, where \(L^p\)-boundedness criteria for such operators were obtained on the torus, bridging the gap between the periodic pseudo-differential theory and the more general setting of FIOs. 
    
    In the specific context of bilinear Fourier integral operators with complex phases on the torus, however, the sharp Seeger-Sogge-Stein type orders remain to be fully characterised, motivating the present investigation to get results as in Corollary \ref{cor:torus:FIOs}.
\end{remark}

\subsection{Structure of the paper}
Apart from this introduction, where we have announced our main results, this paper is organised as follows. 

$\bullet$ In Section \ref{section:prelims}, we present the basic definitions and preliminary results that will be used throughout the paper. In particular, we recall the Hörmander symbol classes, the classes of complex-valued phase functions, and several tools from the theory of Fourier integral operators. We also review the toroidal pseudo-differential calculus and the framework of Fourier analysis on the torus developed by \cite{ruzhansky-turunen}.

$\bullet$  In Section \ref{Bili}, we extend the theory of bilinear Fourier integral operators from the real-valued phase function to the case of complex-valued phase functions. We establish boundedness properties on Lebesgue spaces by adapting the Seeger--Sogge--Stein type estimates to bilinear Fourier integral operators with complex phases. The analysis is divided into low-frequency and high-frequency, where suitable decompositions, oscillatory integral estimates, and Coifman--Meyer type arguments \cite{coifman1978dela} are employed to derive the desired bounds. 

$\bullet$  In Section \ref{endpoint:estimates} we provide the endpoint estimates, these include:
\begin{itemize}
    \item the $(L^2-L^\infty,L^2)$ boundedness theorem,
    \item the $(L^\infty-L^\infty,BMO)$-Fefferman-Stein type boundedness theorem extending a classical estimate due to Fefferman \cite{Fefferman1973}.
    \item  the $(L^\infty-H^1,L^1)$ and the $(H^1-L^\infty,L^1)$ Seeger-Sogge-Stein type boundedness theorems. 
\end{itemize}  Although our method here to get the general $(L^p-L^q,L^r)$ boundedness theorem does not use the standard method of complex interpolation, given that these estimates are proved directly using the Bourgain endpoint maximal function theorem, the end-point estimates cannot be derived from the estimates in Lebesgue spaces, and then, they provide an analogue to the end-point estimates proved by Ruzhansky \cite{Ruzhansky2001}.

$\bullet$  In Section \ref{4}, we study the transference of bilinear Fourier integral operators with complex phases from the Euclidean setting $\mathbb{R}^{n}$ to the torus $\mathbb{T}^{n}$. Using the toroidal quantization developed by \cite{ruzhansky-turunen}, we obtain boundedness results for toroidal pseudo-differential operators and Fourier series operators associated with complex-valued phases. As an application of our results on $\mathbb{R}^n$, we establish new mapping properties for bilinear toroidal Fourier integral operators on periodic function spaces.

\section{Preliminaries}
\label{section:prelims}
\subsection{Overview} In this section, we present the preliminaries about the Fourier analysis on $\mathbb{R}^n$. Our approach follows the framework developed in Staubach et al \cite{rodriguez2011boundedness} for bilinear Fourier integral operators, suitably adapted to the setting of complex-valued phase functions. We also refer to the analysis on the torus and the toroidal pseudo-differential calculus developed in \cite{ruzhansky-turunen}.

Here, we recall several definitions and tools from the linear and bilinear theory of Fourier integral operators that will be used throughout the paper. For the theory of pseudo-differential operators and of Fourier integral operators with smooth symbols see H\"ormander \cite{hormander1967pseudo}.
\subsection{Symbol Classes} 
We begin by recalling the H\"ormander symbol classes, see e.g. \cite{Cardona:2026} where these classes have been used for the analysis of local smoothing estimates. Here, we will show that they play a fundamental role in time-dependent estimates.
\begin{definition}
Let $m \in \mathbb{R}$ and let $0\leq \delta,\rho \leq 1$. A function
\begin{align*}
a:=a(x,\xi,\eta)\in C^{\infty}(\mathbb{R}^{n}\times \mathbb{R}^{n}\times \mathbb{R}^{n})
\end{align*}
belongs to the bilinear symbol class
$S^m_{\rho,\delta}(\mathbb{R}^n_x \times \mathbb{R}^n_{\xi}\times \mathbb{R}^n_{\eta}),$
if for every multi-indices $\alpha,\beta,\gamma \in \mathbb{N}_{0}^{n}$ there exists a constant $C_{\alpha,\beta,\gamma}>0$ such that
\[
\left|
\partial_{\xi}^{\alpha}
\partial_{\eta}^{\beta}
\partial_{x}^{\gamma}
a(x,\xi,\eta)
\right|
\leq C_{\alpha,\beta,\gamma}
( 1+ |\xi|+|\eta|)^{m-\rho(|\alpha|+|\beta|)+\delta|\gamma|}.
\]
When $\rho=1$ and $\delta=0$, we also write
$S^m(\mathbb{R}^n_x \times \mathbb{R}^n_{\xi}\times \mathbb{R}^n_{\eta})=S^m_{1,0}(\mathbb{R}^n_x \times \mathbb{R}^n_{\xi}\times \mathbb{R}^n_{\eta}).$
\end{definition}

\subsection{Phase Functions}

We now introduce the class of phase functions that will be considered throughout the paper for the definition of oscillatory integrals.

\begin{definition}
A function
\[\phi:=
\phi(x,\xi)\in C^{\infty}(\mathbb{R}^{n}\times (\mathbb{R}^{n}\setminus\{0\}))
\]
is called a phase function if it is positively homogeneous of degree one in the frequency variable $\xi \neq 0$, that is,
$
\phi(x,t\xi)=t\phi(x,\xi),
\quad t>0.$
The phase $\phi$ is said to satisfy the non-degeneracy condition if
\[
\det \partial_{x,\xi}^{2}\phi(x,\xi)
\neq 0,
\]
for all $(x,\xi)\in \mathbb{R}^{n}\times (\mathbb{R}^{n}\setminus\{0\})$.
\end{definition}
In this paper, we shall allow the phase functions to be complex-valued. Accordingly, we write
\[
\phi(x,\xi)
=
\Re \phi(x,\xi)
+
i\,\Im \phi(x,\xi),
\]
and assume throughout that
$\Im \phi(x,\xi)\geq 0.$
This positivity assumption guarantees the boundedness of the oscillatory exponential factor.
For a fixed parameter $\tau\in \mathbb{R}$, we define the associated real-valued phase
\[
\phi_\tau\equiv \phi_{\tau}(x,\xi):=\Re \phi(x,\xi)
+\tau \Im \phi(x,\xi).
\]
Consequently,
\[
\phi(x,\xi)=
\phi_{\tau}(x,\xi)
-(\tau-i)\Im \phi(x,\xi).
\]

\subsection{Fourier Integral Operators} 

Let $a:=a(x,\xi)$ be a symbol and $\phi:=\phi(x,\xi)$ be a phase function. The associated linear Fourier integral operator is defined by
\[
T^{\phi}_{a}f(x)
=
\int_{\mathbb{R}^{n}}
e^{2\pi i \phi(x,\xi)}
a(x,\xi)\widehat{f}(\xi)\,d\xi, \quad f\in C_0^{\infty}(\mathbb{R}^n),
\]
 where $$\hat{f}(\xi)= \int_{\mathbb{R}^n}e^{-2\pi i x \cdot \xi} f(x)\diff x, \quad \quad f\in C_0^{\infty}(\mathbb{R}^n),$$ is the Fourier transform of $f.$\\
Similarly, given the symbol
$a(x,\xi,\eta)\in S^m_{1,0}(\mathbb{R}^n_x \times \mathbb{R}^n_{\xi}\times \mathbb{R}^n_{\eta}),$
and phase functions $\phi_{1},\phi_{2}$, we define the bilinear Fourier integral operator by
\begin{align*}
T^{\phi_{1},\phi_{2}}_{a}(f,g)(x)
=
\int_{\mathbb{R}^{n}}
\int_{\mathbb{R}^{n}}
e^{2\pi i(\phi_{1}(x,\xi)+\phi_{2}(x,\eta))}
a(x,\xi,\eta)
\widehat{f}(\xi)\widehat{g}(\eta)
\,d\xi d\eta,\,\,f,g \in C^{\infty}_{0}(\mathbb{R}^n).
\end{align*}
When the phase functions are complex-valued, we rewrite
\[
\phi_{j}(x,\cdot)
=
\phi_{j,\tau}(x,\cdot)
-
(\tau-i)\Im \phi_{j}(x,\cdot),
\qquad j=1,2.
\] 
Hence,
\begin{align}\label{op}
T^{\phi_{1,\tau},\phi_{2,\tau}}_{a_{\tau}}(f,g)(x)
=
\int_{\mathbb{R}^{n}}
\int_{\mathbb{R}^{n}}
e^{2\pi i(\phi_{1,\tau}(x,\xi)+\phi_{2,\tau}(x,\eta))}
a_{\tau}(x,\xi,\eta)
\widehat{f}(\xi)\widehat{g}(\eta)
\,d\xi d\eta,
\end{align}
with modified symbol
$$
a_{\tau}(x,\xi,\eta)
=e^{-2\pi i(\tau-i)(\Im \phi_{1}(x,\xi)+\Im \phi_{2}(x,\eta))}
a(x,\xi,\eta).$$ 
The analysis of endpoint estimates of FIOs requires the following notion of a Carleson measure:
\begin{definition}[Carleson measure]
A Borel measure $d\mu$ on $\mathbb{R}^{n}_{+}$ is called a Carleson measure if
\begin{align*}
\|d\mu\|_{\mathcal{C}}
=\sup_{\epsilon>0}\sup_{y\in\mathbb{R}^{n}}\frac{1}{\epsilon^{n}}\int_{0}^{\epsilon}
\int_{|x-y|<\epsilon}d|\mu|(x,t) < \infty.
\end{align*}
The quantity $\|d\mu\|_{\mathcal{C}}$ is called the Carleson norm of $d\mu$. Here, $d|\mu|,$ denotes the total variation measure of the Borel measure $d\mu.$
\end{definition}
For our further analysis of periodic Fourier integral operators we require the following definition.
\begin{definition}[Periodic functions] 
     A function is 1-periodic if $f(x) = f(x + k)$ for every $ x \in \mathbb{R}^n $ and all $ k \in \mathbb{Z}^n $. We can identify this function with one defined in $\mathbb{T}^n := \mathbb{R}^n / \mathbb{Z}^n $, which is a compact manifold without boundary. Moreover, the space of 1-periodic functions $m$ times continuously differentiable is denoted by $C^m(\mathbb{T}^n)$. The test functions are elements of the space $C^\infty(\mathbb{T}^n) := \bigcap_m C^m(\mathbb{T}^n)$. 
\end{definition}
\subsection{Fourier series operators}

In order to define the class of Fourier series operators on $\mathbb{T}^n$, we first need to define the group Fourier transform of 1-periodic smooth functions.
\begin{definition}[Group Fourier Transform on $\mathbb{T}^n$]
    For $f \in C^\infty(\mathbb{T}^n)$ we define the  \textit{toroidal Fourier transform} $\mathcal{F}_{\mathbb{T}^n}$ as follows
    \begin{equation*}
        (\mathcal{F}_{\mathbb{T}^n}f)(\xi) := \int_{\mathbb{T}^n} e^{-i2\pi x \cdot \xi}f(x)\diff  x, 
    \end{equation*}
    as a function of $\xi \in \mathbb{Z}^n$. Here, we have denoted $\diff x,$ the Haar measure on $\mathbb{Z}^n.$
\end{definition}
In view of the group structure of the torus $\mathbb{T}^n$, the corresponding frequency domain of 1-periodic smooth functions is the lattice $\mathbb{Z}^n$. Also, it can be proved that the corresponding functions $\mathcal{F}_{\mathbb{T}^n}f$ satisfy certain regularity conditions. In fact, they belong to the Schwartz space defined below.
Here we use $\langle\xi\rangle:= (1+|\xi|^2)^{1/2}$.
\begin{definition}[Schwartz space $\mathcal{S}(\mathbb{Z}^n) $]
    We say $\varphi \in \mathcal{S}(\mathbb{Z}^n)$, that is a \textit{rapidly decaying} function $\mathbb{Z}^n \rightarrow \mathbb{C}$, if for given $0<M<\infty$, there exists $C_{\varphi M} > 0$ such that 

    \begin{equation*}
         \forall\, \xi \in \mathbb{Z}^n ,\quad|\varphi(\xi)| \leq C_{\varphi M} \langle\xi\rangle^{-M}.
    \end{equation*}
\end{definition}
It can be proved that the toroidal Fourier transform is a bijection between $C^\infty(\mathbb{T}^n)$ and $\mathcal{S}(\mathbb{Z}^n)$ with inverse  $\mathcal{F}_{\mathbb{T}^n}^{-1}:  \mathcal{S}(\mathbb{Z}^n) \rightarrow C^\infty(\mathbb{T}^n)$ defined by
    \begin{equation*}
        (\mathcal{F}_{\mathbb{T}^n}^{-1}\varphi)(x) := \sum_{\xi \in \mathbb{Z}^n} e^{i2\pi x \cdot \xi} \varphi(\xi).
    \end{equation*}
In the Euclidean case, the Hörmander symbol classes are defined requiring some regularity conditions using partial derivatives in both the space and frequency domains. However, in the case of the torus, the frequency domain is the lattice $\mathbb{Z}^n$. Therefore, we define the natural analogue of the derivative in the discrete case: the partial difference operators.
\begin{definition}[Partial difference operators]
    Let $e_j \in \mathbb{Z}^n$ be such that $(e_j)_i$ is equal to 1, if $i=j$, and let it be equal to 0 otherwise. Then, for $\sigma(\xi):\mathbb{Z}^n \rightarrow \mathbb{C} $, we define 
    \begin{equation*}
        \Delta_{\xi_j}\sigma(\xi) := 
        \sigma(\xi + e_j) - \sigma(\xi) \quad \text{and} \quad \Delta^\alpha_\xi := \Delta_{\xi_1}^{\alpha_1} \cdots \Delta_{\xi_n}^{\alpha_n},
    \end{equation*}
    for any multi-index $\alpha \in \mathbb{N}^n_0$.
\end{definition}
Now we proceed to define the Hörmander classes of symbols on $ \mathbb{T}^n \times \mathbb{Z}^n $, as in Ruzhansky and Turunen \cite{ruzhansky-turunen}, that will be important in the quantization of pseudo-differential operators.

\begin{definition}[Multilinear Hörmander classes on $ \mathbb{T}^n \times \mathbb{Z}^n $]
    Let $m \in \mathbb{R} $, and let $0 \leq \delta, \rho\leq 1$. We say that a function $a:=a(x, \xi_1, \ldots, \xi_r)$ that is smooth on $x$ for any $\xi \in \mathbb{Z}^{nr}$ belongs to the \emph{multilinear toroidal symbol class} $S^m_{\rho,\delta} (\mathbb{T}^n\times \mathbb{Z}^n,r) $ if 
    \begin{equation*}
        \left|\partial^\beta_x \Delta^{\alpha_1}_{\xi_1}\cdots\Delta^{\alpha_r}_{\xi_r} a(x, \xi_1,\ldots,\xi_r)\right| \leq C_{\alpha\beta}(1+|\xi_1|+\cdots+|\xi_r|)^{m-\rho|\alpha| + \delta|\beta|},
    \end{equation*}
    for every $x \in \mathbb{T}^n$, $\alpha\in\mathbb{N}^{nr}_0$, $\beta \in \mathbb{N}^n_0$ and $\xi \in \mathbb{Z}^{nr}$. Moreover, we say that $a(x, \xi)$ has order $m$ and we define $S^{-\infty}_{\rho,\delta} (\mathbb{T}^n\times \mathbb{Z}^n, r) = \bigcap_{m\in\mathbb{R}}S^m_{\rho,\delta} (\mathbb{T}^n\times \mathbb{Z}^n, r) $. When $\rho=1$ and $\delta=0$, we simply denote $S^m(\mathbb{T}^n \times \mathbb{Z}^n, r).$
\end{definition}
\begin{definition}[Multilinear FSO on the torus]
    Let $a(x, \xi) \in S^m_{\rho,\delta}(\mathbb{T}^n\times \mathbb{Z}^n, r)$ and $\phi_1(x,\xi_1),\ldots,\phi_r(x,\xi_r)$ be phase functions. We define the \emph{multilinear Fourier series operator} as 
    \begin{equation*}
        T^{\phi_1,\ldots,\phi_r}_a (f_1, \ldots, f_r)(x) := \sum_{\xi\in\mathbb{R}^{nr}}  e^{2\pi i\phi(x,\xi)}a(x,\xi)(\mathcal{F}_{\mathbb{T}^n}f_1)\cdots(\xi_1)(\mathcal{F}_{\mathbb{T}^n}f_r)(\xi_r)
    \end{equation*}
    where $\phi(x,\xi):= \phi_1(x,\xi_1)+\cdots+\phi_r(x,\xi_r).$
\end{definition}

\section{Bilinear Fourier integral operators on $\mathbb{R}^n$} \label{Bili}
In this section, we study the boundedness of bilinear Fourier integral operators on $\mathbb{R}^n$ with complex-valued phase functions. All the section will be dedicated to the proof of Theorem \ref{Main:1:th} presented below in the form of Theorem \ref{Th:FIO:C:2}.

\subsection{Extension to Complex-Valued Phase Functions}

Here we extend the boundedness results for bilinear Fourier integral operators, originally established for real-valued phase functions, to the setting of complex-valued phase functions.
 We have defined the bilinear Fourier integral operator $T^{\phi_1, \phi_2}_a$ as
\begin{align}\label{aa}
&T^{\phi_1, \phi_2}_a(f,g)(x)
\nonumber\\&=\int_{\mathbb{R}^n} \int_{\mathbb{R}^n} 
e^{2\pi i \bigl(\phi_1(x,\xi)+\phi_2(x,\eta)\bigr)}
\, a(x,\xi,\eta)\, \widehat{f}(\xi)\, \widehat{g}(\eta)
\, d\xi \, d\eta, \text{ for } f,g \in C^{\infty}_{0}(\mathbb{R}^n),
\end{align}
where the symbol $a(x,\xi,\eta)$ is scalar-valued and compactly supported in the spatial variable $x$. To extend the boundedness analysis of bilinear Fourier integral operators to the complex-valued setting, we assume that the phase functions $\phi_1(x,\xi)$ and $\phi_2(x,\eta)$ are complex-valued and can be expressed as
\begin{align*}
\phi_1(x,\xi) &= \Re \phi_1(x,\xi) + i\, \Im \phi_1(x,\xi),
\qquad \Im \phi_1(x,\xi) \ge 0; \\
\phi_2(x,\eta) &= \Re \phi_2(x,\eta) + i\, \Im \phi_2(x,\eta), 
\qquad \Im \phi_2(x,\eta) \ge 0.
\end{align*}
Let $\tau \in \mathbb{R}$ be a fixed real parameter such that we express $\phi_1(x,\xi)$ as
\begin{align*}
\phi_1(x,\xi)
&= \Re \phi_1(x,\xi) + i\, \Im \phi_1(x,\xi) \\
&= \Re \phi_1(x,\xi) + \tau \Im \phi_1(x,\xi) - \Im \phi_1(x,\xi)(\tau - i) \\
&= \phi_{1,\tau}(x,\xi) - \Im \phi_1(x,\xi)(\tau - i),
\end{align*}
where
\begin{align*}
\phi_{1,\tau}(x,\xi) = \Re \phi_1(x,\xi) + \tau \Im \phi_1(x,\xi).
\end{align*}
Similarly, we write
\begin{align*}
\phi_2(x,\eta)
&= \Re \phi_2(x,\eta) + \tau \Im \phi_2(x,\eta) - \Im \phi_2(x,\eta)(\tau - i) \\
&= \phi_{2,\tau}(x,\eta) - \Im \phi_2(x,\eta)(\tau - i),
\end{align*}
where
\begin{align*}
\phi_{2,\tau}(x,\eta) = \Re \phi_2(x,\eta) + \tau \Im \phi_2(x,\eta).
\end{align*}
Our hypothesis are that $\phi_{1,\tau}(x,\xi)$ and $\phi_{2,\tau}(x,\eta)$ are non-degenerate phases for some $\tau\in \mathbb{R}.$ That is, 
\begin{align}
\mathrm{det}[\partial_x\partial_\xi\phi_{1,\tau}\label{assump}(x,\xi)]\neq 0 \text{ and } \mathrm{det}[\partial_x\partial_\eta\phi_{2,\tau}(x,\eta)]\neq 0 .
\end{align}
\begin{remark}
    The above assumptions are not restrictive, since if \eqref{assump} holds for some $\tau\in \mathbb{R},$ then it holds for almost all $\tau\in\mathbb{C},$ see Remark \ref{Rem:MS:condition}.
\end{remark}
We now express the bilinear Fourier integral operator in \eqref{aa} in terms of the complex phases as follows:
\begin{align*}
&T^{\phi_1, \phi_2}_a(f,g)(x)\\
&= \int_{\mathbb{R}^n}\int_{\mathbb{R}^n} 
e^{2\pi i \bigl(\phi_1(x,\xi)+\phi_2(x,\eta)\bigr)}
\, a(x,\xi,\eta)\, \widehat{f}(\xi)\, \widehat{g}(\eta)
\, d\xi \, d\eta \nonumber\\
&= \int_{\mathbb{R}^n}\int_{\mathbb{R}^n} 
e^{2\pi i \bigl(\Re\phi_1(x,\xi) + i\,\Im\phi_1(x,\xi) 
+ \Re\phi_2(x,\eta) + i\,\Im\phi_2(x,\eta)\bigr)}
\, a(x,\xi,\eta)\, \widehat{f}(\xi)\, \widehat{g}(\eta)
\, d\xi \, d\eta \nonumber\\
&= \int_{\mathbb{R}^n}\int_{\mathbb{R}^n} 
e^{2\pi i \bigl(\phi_{1,\tau}(x,\xi) - (\tau - i)\Im\phi_1(x,\xi)
+ \phi_{2,\tau}(x,\eta) - (\tau - i)\Im\phi_2(x,\eta)\bigr)}
\, a(x,\xi,\eta)\, \widehat{f}(\xi)\, \widehat{g}(\eta)
\, d\xi \, d\eta \nonumber\\
&= \int_{\mathbb{R}^n}\int_{\mathbb{R}^n} 
e^{2\pi i \bigl(\phi_{1,\tau}(x,\xi)+\phi_{2,\tau}(x,\eta)\bigr)}
e^{-2\pi i (\tau - i)\bigl(\Im\phi_1(x,\xi)+\Im\phi_2(x,\eta)\bigr)}
\, \\& \quad \times a(x,\xi,\eta)\, \widehat{f}(\xi)\, \widehat{g}(\eta)
\, d\xi \, d\eta.
\end{align*}
We may write the bilinear Fourier integral operator as another bilinear Fourier integral operator defined as
\begin{align}
T^{\phi_{1,\tau}, \phi_{2,\tau}}_{a_{12,\tau}}(f,g)(x)
&= \int_{\mathbb{R}^n}\int_{\mathbb{R}^n}
e^{2\pi i \phi_{12,\tau}(x,\xi,\eta)}
\, a_{12,\tau}(x,\xi,\eta)\, \widehat{f}(\xi)\, \widehat{g}(\eta)
\, d\xi \, d\eta,
\end{align}
where the new symbol is defined by
$$
a_{12,\tau}(x,\xi,\eta)
= e^{-2\pi i(\tau - i)\bigl(\Im\phi_1(x,\xi)+\Im\phi_2(x,\eta)\bigr)}\, a(x,\xi,\eta),$$
and the combined phase is given by
$$
\phi_{12,\tau}(x,\xi,\eta)
= \phi_{1,\tau}(x,\xi)+\phi_{2,\tau}(x,\eta).$$
For the boundedness of bilinear Fourier integral operators $T_{a}^{\phi_1, \phi_2}$, we break our analysis into low and high frequencies.
To this end, we introduce a smooth function 
$\mu:\mathbb{R}^{n}\to\mathbb{R}$
such that
$\mu(\xi)=0, \text{for } |\xi|\leq \frac{1}{4\lambda},$
and $\mu(\xi)=1,  \text{for } |\xi|\geq \frac{1}{3\lambda},$ where $\lambda$ is a parameter satisfying $0<\lambda\leq 1.$
Now, we decompose the operator $T_{a}^{\phi_1, \phi_2}$ as follows
\begin{align}\label{modified:1}
    T_{a}^{\phi_1, \phi_2}(f,g)&= T_{a}^{\phi_1, \phi_2}(\mu(D)f, \mu(D)g)\nonumber\\&+T_{a}^{\phi_1, \phi_2}(\mu(D)f, (1-\mu)(D)g)+T_{a}^{\phi_1, \phi_2}((1-\mu(D))f,g)\nonumber\\&
    =T_{a_{12,\tau}}^{\phi_{1,\tau}, \phi_{2,\tau}}(\mu(D)f, \mu(D)g)\nonumber\\&+T_{a_{12,\tau}}^{\phi_{1,\tau}, \phi_{2,\tau}}(\mu(D)f, (1-\mu)(D)g)+T_{a_{12,\tau}}^{\phi_{1,\tau}, \phi_{2,\tau}}((1-\mu(D))f,g),
\end{align}
where $\mu(D)$ and $(1-\mu)(D)$ denote the Fourier multiplier operators with symbols $\mu(\xi)$ and $\bigl(1-\mu(\xi)\bigr),$ respectively.

\subsection{Structure of the high frequency part}

The aim of this subsection is to prove that we have the following representation
$$T^{{\phi_{1,\tau}},{\phi_{2,\tau}}}_{a_{12,\tau;1}}(\mu(D)f,\mu(D)g)(x)$$
\begin{equation*}\label{representation:high:freq}
=\displaystyle{ \int\limits_{0}^1\int_{\mathbb{R}^n}\int_{\mathbb{R}^n}}T^{\phi_{1,\tau}}_{\nu_{1}^{t,u}} f(x)
T^{\phi_{2,\tau}}_{\mu_{1}^{t,v}} g(x)
\frac{m_{\tau}(t,x,u,v)}{(1+4\pi^2|u|^2+4\pi^2|v|^2)^N}
\,du\,dv\,\frac{dt}{t},
\end{equation*} 
holds for the high frequency part. We will prove that for any integer $N,$ there exists a function $m_{\tau}:=m_\tau^{(N)}$ depending on $N$ such that  \eqref{representation:high:freq} holds.

\subsection{Construction of $m_{\tau}(t, x,u,v)$}
In this subsection, we construct the function $m_{\tau}(t, x,u,v)$ associated with the modified symbol of the bilinear Fourier integral operator, $T_{a_{12,\tau}}^{\phi_{1,\tau}, \phi_{2,\tau}}$.
Consider the new symbol 
\begin{align*}
a_{12,\tau}(x,\xi,\eta)&=e^{-2\pi i(\Im\phi_{1,\tau}(x,\xi)+\Im\phi_{2,\tau}(x,\eta))(\tau-i)}a(x,\xi,\eta)\\&
=e^{-2\pi\tau i(\Im\phi_{1,\tau}(x,\xi)+\Im\phi_{2,\tau}(x,\eta))}e^{-2\pi (\Im\phi_{1,\tau}(x,\xi)+\Im\phi_{2,\tau}(x,\eta))}a(x,\xi,\eta).
\end{align*}
To deal with the term
\begin{align*}
T_{a_{12,\tau}}^{\phi_{1,\tau}, \phi_{2,\tau}}\bigl(\mu(D)f,\mu(D)g\bigr),
\end{align*}from \eqref{modified:1} let us introduce two smooth cut-off functions
$\chi,\nu:\mathbb{R}^{2n}\to\mathbb{R},$
such that
$\chi(\xi,\eta)=1,
\quad \text{for } |(\xi,\eta)|\le 1,$
and
$\chi(\xi,\eta)=0,
\quad \text{for } |(\xi,\eta)|\ge 2,$
and
$\nu(\xi,\eta)=0,
\quad \text{for } \lambda^{2}|\xi|\le 16|\eta|,
$
and
$\nu(\xi,\eta)=1,\quad \text{for } 64|\eta|\le \lambda^{2}|\xi|.$
Defining
\begin{align*}
a_{12,\tau;1}(x,\xi,\eta)=\bigl(1-\chi(\xi,\eta)\bigr)\nu(\xi,\eta)
a_{12,\tau}(x,\xi,\eta),
\end{align*}
and
\begin{align*}
a_{12,\tau;2}(x,\xi,\eta)=\bigl(1-\chi(\xi,\eta)\bigr)\bigl(1-\nu(\xi,\eta)\bigr)a_{12,\tau}(x,\xi,\eta),
\end{align*}
we have that
\begin{align*}
a_{12,\tau;1}(x,\xi,\eta),a_{12,\tau;2}(x,\xi,\eta)\in S^m_{1/2,1/2}(\mathbb{R}^n_x \times \mathbb{R}^n_{\xi}\times \mathbb{R}^n_{\eta}),
\end{align*} see \cite[Page~6]{cardona2026weak}
and we can decompose
\begin{align*}
T_{a_{12,\tau}}^{\phi_{1,\tau}, \phi_{2,\tau}}
\bigl(\mu(D)f,\mu(D)g\bigr)
=T_{a_{12,\tau;1}}^{\phi_{1,\tau}, \phi_{2,\tau}}
\bigl(\mu(D)f,\mu(D)g\bigr)
+T_{a_{12,\tau;2}}^{\phi_{1,\tau}, \phi_{2,\tau}}
\bigl(\mu(D)f,\mu(D)g\bigr).
\end{align*}
One introduces an even real-valued smooth function $\psi$ whose Fourier transform is
supported on the annulus
$\left\{\xi:\frac12\le |\xi|\le 2\right\}$
such that
\begin{align}\label{decomposition}
\int_{0}^{\infty}
\left|\widehat{\psi}(t\xi)\right|^{2}
\frac{dt}{t}=1,
\end{align} 
for $\xi\neq0$.
Let $\theta$ be another real-valued smooth function whose Fourier transform is equal
to one on the ball
$\{\xi:|\xi|\le\tfrac18\}$
and supported in
$\{\xi:|\xi|\le\tfrac14\}.$
Here and in the sequel we denote
$\mathcal{K}=\operatorname{supp}_{x}a_{12,\tau}.$
\begin{proposition}
Let $T_{a_{12,\tau;1}}^{\phi_{1,\tau},\phi_{2,\tau}}$ be a bilinear Fourier integral operator with symbol $a_{12,\tau;1}(x,\xi,\eta) \in S^m_{1/2,1/2}(\mathbb{R}^n_x \times \mathbb{R}^n_{\xi}\times \mathbb{R}^n_{\eta}).$ Assume that $\phi_{1,\tau}(x,\xi)$ and $\phi_{2,\tau}(x,\eta)$ are positively homogeneous of degree one in the frequency variables.
For $t>0$, we define
\begin{align}\label{atau}
a_{12,\tau;1,t}(x,\xi,\eta)
& =t^{m}a_{12,\tau;1}
\left(x,\frac{\xi}{t},\frac{\eta}{t}\right)
\left|\nabla_x\phi_{1,\tau}(x,\xi)\right|^{-m}
\nonumber\\&\quad\times\widehat{\psi}
\bigl(\nabla_x\phi_{1,\tau}(x,\xi)\bigr)
\widehat{\theta}
\bigl(\nabla_x\phi_{2,\tau}(x,\eta)\bigr).
\end{align} Then
\begin{align}\label{a}
&T_{a_{12,\tau;1}}^{\phi_{1,\tau}, \phi_{2,\tau}}
\bigl(\mu(D)f,\mu(D)g\bigr)(x)
\nonumber\\&=\int_{0}^{\infty}
\iint a_{12,\tau;1,t}(x,t\xi,t\eta)
\Bigl(\widehat{\psi}\bigl(t\nabla_x\phi_{1,\tau}(x,\xi)\bigr)\nonumber\\&\quad\times\left|\nabla_x\phi_{1,\tau}(x,\xi)\right|^{m}\mu(\xi)\widehat{f}(\xi)e^{i\phi_{1,\tau}(x,\xi)}\Bigr)
\nonumber\\&\quad\times
\Bigl(\widehat{\theta}\bigl(t\nabla_x\phi_{2,\tau}(x,\eta)\bigr)\mu(\eta)\widehat{g}(\eta)e^{i\phi_{2,\tau}(x,\eta)}\Bigr)\frac{dt}{t}\,d\xi\,d\eta.
\end{align}
\end{proposition}
\begin{proof}
By the definition of bilinear Fourier integral operator, we have
\begin{align*}
&T_{a_{12,\tau;1}}^{\phi_{1,\tau},\phi_{2,\tau}}
\bigl(\mu(D)f,\mu(D)g\bigr)(x)
\nonumber\\
&\quad =\iint
a_{12,\tau;1}(x,\xi,\eta)
\mu(\xi)\widehat{f}(\xi)
\mu(\eta)\widehat{g}(\eta)
e^{i\phi_{1,\tau}(x,\xi)}
e^{i\phi_{2,\tau}(x,\eta)}
\,d\xi\,d\eta,
\end{align*}
using the identity \eqref{decomposition} for $\widehat{\psi}\bigl(t\nabla_x\phi_{1,\tau}(x,\xi)$ with $\nabla_x\phi_{1,\tau}(x,\xi)\neq 0,$ by substituting $\displaystyle\int_0^\infty \left|\widehat{\psi}\bigl(t \nabla_x \phi_{1,\tau}(x,\xi)\bigr)\right|^2 \frac{dt}{t}=1,$ into the above expression gives us
\begin{align*}
    &= \iint
\bigg(\int_0^\infty\left|\widehat{\psi}\bigl(t \nabla_x \phi_{1,\tau}(x,\xi)\bigr)\right|^2 \frac{dt}{t}\bigg)a_{12,\tau;1}(x,\xi,\eta)
\mu(\xi)\widehat{f}(\xi)\\&\quad \times
\mu(\eta)\widehat{g}(\eta)
e^{i\phi_{1,\tau}(x,\xi)}
e^{i\phi_{2,\tau}(x,\eta)}
\,d\xi\,d\eta
\end{align*}
Since $\psi$ is an even real-valued function, we write $$\displaystyle \left|\widehat{\psi}\bigl(t \nabla_x \phi_{1,\tau}(x,\xi)\bigr)\right|^2 = \widehat{\psi}\bigl(t \nabla_x \phi_{1,\tau}(x,\xi)\bigr)\widehat{\psi}\bigl(t \nabla_x \phi_{1,\tau}(x,\xi)\bigr).$$ One has that,
\begin{align*}
    &=\iint
\bigg(\int_0^\infty \widehat{\psi}\bigl(t \nabla_x \phi_{1,\tau}(x,\xi)\bigr)\widehat{\psi}\bigl(t \nabla_x \phi_{1,\tau}(x,\xi)\bigr) \frac{dt}{t}\bigg)a_{12,\tau;1}(x,\xi,\eta)
\mu(\xi)\widehat{f}(\xi)\\&\quad \times
\mu(\eta)\widehat{g}(\eta)
e^{i\phi_{1,\tau}(x,\xi)}
e^{i\phi_{2,\tau}(x,\eta)}
\,d\xi\,d\eta,
\end{align*}
since $\widehat{\theta}$ is equal to one on the support of the frequency variables, implies that $\widehat{\theta}
\bigl(t\nabla_x\phi_{2,\tau}(x,\eta)\bigr)=1,$ and applying Fubini's theorem to interchange the order of integration yields,
\begin{align*}
    &= \int_0^\infty\iint
\bigg( a_{12,\tau;1}(x,\xi,\eta)\widehat{\psi}\bigl(t \nabla_x \phi_{1,\tau}(x,\xi)\bigr) \widehat{\theta}
\bigl(t\nabla_x\phi_{2,\tau}(x,\eta)\bigr)\bigg)
\\&\quad\times \widehat{\psi}\bigl(t \nabla_x \phi_{1,\tau}(x,\xi)\bigr)\mu(\xi)\widehat{f}(\xi)e^{i\phi_{1,\tau}(x,\xi)}\\&\quad \times
\widehat{\theta}
\bigl(t\nabla_x\phi_{2,\tau}(x,\eta)\bigr)\mu(\eta)\widehat{g}(\eta)
e^{i\phi_{2,\tau}(x,\eta)}
\frac{dt}{t}\,d\xi\,d\eta
\end{align*}
using the identity from Remark \ref{remar}. That is,
\begin{align*}
a_{12,\tau;1}\left(x,\xi,\eta\right)\widehat{\psi}\bigl(t\nabla_x\phi_{1,\tau}(x,\xi)\bigr)\widehat{\theta}\bigl(t\nabla_x\phi_{2,\tau}(x,\eta)\bigr) = a_{12,\tau;1,t}(x,t\xi,t\eta)\left|\nabla_x\phi_{1,\tau}(x,\xi)\right|^{m},
\end{align*}
yields
\begin{align*}
&=\int_{0}^{\infty}
\iint a_{12,\tau;1,t}(x,t\xi,t\eta)
\Bigl(\widehat{\psi}\bigl(t\nabla_x\phi_{1,\tau}(x,\xi)\bigr)\left|\nabla_x\phi_{1,\tau}(x,\xi)\right|^{m}\mu(\xi)\widehat{f}(\xi)e^{i\phi_{1,\tau}(x,\xi)}\Bigr)
\nonumber\\&\quad\times
\Bigl(\widehat{\theta}\bigl(t\nabla_x\phi_{2,\tau}(x,\eta)\bigr)\mu(\eta)\widehat{g}(\eta)e^{i\phi_{2,\tau}(x,\eta)}\Bigr)\frac{dt}{t}\,d\xi\,d\eta.
\end{align*}
This completes the proof.
\end{proof}
\begin{remark}\label{remar}
    From \eqref{atau}, we evaluate  $a_{12,\tau;1,t}(x,t\xi,t\eta)$ as follows:
    \begin{align*}
a_{12,\tau;1,t}(x,t\xi,t\eta)
& =t^{m}a_{12,\tau;1}
\left(x,\xi,\eta\right)
\left|\nabla_x\phi_{1,\tau}(x,t\xi)\right|^{-m}
\nonumber\\&\quad\times\widehat{\psi}
\bigl(\nabla_x\phi_{1,\tau}(x,t\xi)\bigr)
\widehat{\theta}
\bigl(\nabla_x\phi_{2,\tau}(x,t\eta)\bigr)
\\&= t^{m}a_{12,\tau;1}
\left(x,\xi,\eta\right)
\left|t\nabla_x\phi_{1,\tau}(x,\xi)\right|^{-m}
\nonumber\\&\quad\times\widehat{\psi}
\bigl(t\nabla_x\phi_{1,\tau}(x,\xi)\bigr)
\widehat{\theta}
\bigl(t\nabla_x\phi_{2,\tau}(x,\eta)\bigr)\\&
= t^{m}a_{12,\tau;1}
\left(x,\xi,\eta\right)
t^{-m}\left|\nabla_x\phi_{1,\tau}(x,\xi)\right|^{-m}
\nonumber\\&\quad\times\widehat{\psi}
\bigl(t\nabla_x\phi_{1,\tau}(x,\xi)\bigr)
\widehat{\theta}
\bigl(t\nabla_x\phi_{2,\tau}(x,\eta)\bigr)\\&
=a_{12,\tau;1}
\left(x,\xi,\eta\right)
\left|\nabla_x\phi_{1,\tau}(x,\xi)\right|^{-m}
\nonumber\\&\quad\times\widehat{\psi}
\bigl(t\nabla_x\phi_{1,\tau}(x,\xi)\bigr)
\widehat{\theta}
\bigl(t\nabla_x\phi_{2,\tau}(x,\eta)\bigr).
\end{align*}
Consequently, $$a_{12,\tau;1,t}(x,t\xi,t\eta)\left|\nabla_x\phi_{1,\tau}(x,\xi)\right|^{m}=a_{12,\tau;1}
\left(x,\xi,\eta\right)\widehat{\psi}
\bigl(t\nabla_x\phi_{1,\tau}(x,\xi)\bigr)
\widehat{\theta}
\bigl(t\nabla_x\phi_{2,\tau}(x,\eta)\bigr).$$
\end{remark}
We define:
\begin{align}\label{Finish}
\tilde{a}_{12,\tau;1}(t,x,\xi,\eta)
:=a_{12,\tau;1,t}
\bigl(x,\Psi_{1}(x,\xi),\Psi_{2}(x,\eta)
\bigr),
\end{align}
where
\begin{align*}
\Psi_j(x,\xi)
=\bigl(\nabla_x\phi_{j,\tau}(x,\cdot)\bigr)^{-1}(\xi),\qquad j=1,2,
\end{align*}
whose existence follows from the assumption of compact support in $x$, homogeneity of degree one in $\xi$, and the non-degeneracy.
Note that this argument is the same as in \cite[Page 12]{rodriguez2011boundedness}. 
By the Fourier inversion formula and integration by parts, we have that,
\begin{align}\label{F12}
  &\widetilde{a}_{12,\tau;1}(t,x,\xi,\eta)\nonumber\\&=\int_{\mathbb{R}^n}\int_{\mathbb{R}^n} e^{2\pi i\xi\cdot u+ 2 \pi i\eta \cdot v}\mathscr{F}(\widetilde{a}_{12,\tau;1}(t,x,\cdot,\cdot))(u,v) du \,dv\nonumber\\&
  =\int_{\mathbb{R}^n}\int_{\mathbb{R}^n} \dfrac{1}{(1+4\pi^2|u|^2+4\pi^2|v|^2)^N}(1-\Delta_{\xi}-\Delta_{\eta})^N[e^{2\pi i\xi\cdot u+ 2\pi i\eta \cdot v}]\nonumber\\& \quad \times \mathscr{F}(\widetilde{a}_{12,\tau;1}(t,x,\cdot,\cdot))(u,v) du\,dv\nonumber\\&
  = \int_{\mathbb{R}^n}\int_{\mathbb{R}^n} \dfrac{e^{ 2\pi i\xi\cdot u+ 2\pi i\eta \cdot v}}{(1+4\pi^2|u|^2+ 4\pi^2|v|^2)^N}\nonumber\\&
  \hspace{5cm}\times[(1-\Delta_{\xi}-\Delta_{\eta})^N(\mathscr{F}(\widetilde{a}_{12,\tau;1}(t,x,\cdot,\cdot)))(u,v)] du\,dv  \nonumber\\&
  =\int_{\mathbb{R}^n}\int_{\mathbb{R}^n} e^{ 2\pi i\xi\cdot u+ 2\pi i\eta \cdot v}\dfrac{m_{\tau}(t,x,u,v)}{(1+4\pi^2|u|^2+4\pi^2|v|^2)^N}\,du\,dv,
\end{align} 
where \begin{align}
    m_{\tau}(t,x,u,v)&= (1-\Delta_{\xi}-\Delta_{\eta})^N(\mathscr{F}(\widetilde{a}_{12,\tau;1}(t,x,\cdot,\cdot)))(u,v)\nonumber\\&
    =(1-\Delta_{\xi}-\Delta_{\eta})^N \int_{\mathbb{R}^n}\int_{\mathbb{R}^n} e^{-2\pi i\xi\cdot u- 2\pi i\eta \cdot v}\widetilde{a}_{12,\tau;1}(t,x,\xi,\eta)\,d\xi\,d\eta\nonumber\\&
    = \int_{\mathbb{R}^n}\int_{\mathbb{R}^n} e^{-2\pi i\xi\cdot u- 2\pi i\eta \cdot v} [(1-\Delta_{\xi}-\Delta_{\eta})^N\widetilde{a}_{12,\tau;1}(t,x,\xi,\eta)]\,d\xi\,d\eta,\label{m:tau}
\end{align} for a fixed large fixed $N\in \mathbb{N}$. Since the $(\xi,\eta)$-support of $\tilde{a}_{12,\tau;1}$ is contained in a compact set independent of $t$ and $x$. 

\begin{proposition}
Let us consider the symbol
$$\tilde{a}_{12,\tau;1}(t,x,\xi,\eta)
:=a_{12,\tau;1,t}
\bigl(x,\Psi_{1}(x,\xi),\Psi_{2}(x,\eta)
\bigr),$$
where
$
\Psi_j(x,\xi)
=\bigl(\nabla_x\phi_{j,\tau}(x,\cdot)\bigr)^{-1}(\xi),\quad j=1,2.$ Then there exists a smooth function $m_\tau(t,x,u,v)$
such that, we represent $a_{12,\tau;1,t}(x,\xi,\eta)$ as 
\begin{align*}
a_{12,\tau;1,t}(x,\xi,\eta)
=\iint_{\mathbb{R}^{2n}}
e^{2\pi i\nabla_x\phi_{1,\tau}(x,\xi)\cdot u
+2\pi i\nabla_x\phi_{2,\tau}(x,\eta)\cdot v}
\frac{m_{\tau}(t,x,u,v)}
{(1+4\pi^2|u|^2+4\pi^2|v|^2)^N}
\,du\,dv.
\end{align*}
\end{proposition}
\begin{proof}
   Since, for $j=1,2$, $\Psi_j(x,\cdot)$ is the inverse of the mapping
$\nabla_x\phi_{j,\tau}(x,\cdot)$, we have
\begin{align}
\Psi_1\!\left(x,\nabla_x\phi_{1,\tau}(x,\xi)\right)=\xi,
 \text{ and }
\Psi_2\!\left(x,\nabla_x\phi_{2,\tau}(x,\eta)\right)=\eta.
\end{align}
Since the identity in \eqref{F12} holds for every $(\xi,\eta)\in\mathbb{R}^{2n}$, we evaluate it at
\[
\left(\nabla_x\phi_{1,\tau}(x,\xi),
\nabla_x\phi_{2,\tau}(x,\eta)\right).
\]
That is, replacing
\[
\xi\mapsto\nabla_x\phi_{1,\tau}(x,\xi),
\text{ and }
\eta\mapsto\nabla_x\phi_{2,\tau}(x,\eta),
\] and putting them into $$\tilde{a}_{12,\tau;1}(t,x,\xi,\eta)
=a_{12,\tau;1,t}
\bigl(x,\Psi_{1}(x,\xi),\Psi_{2}(x,\eta)
\bigr),$$ yields 
\begin{align*}
&\widetilde{a}_{12,\tau;1}
\left(t,x,
\nabla_x\phi_{1,\tau}(x,\xi),
\nabla_x\phi_{2,\tau}(x,\eta)
\right) \\
&=a_{12,\tau;1,t}
\Bigl(x,\Psi_1\!\left(x,\nabla_x\phi_{1,\tau}(x,\xi)\right),
\Psi_2\!\left(x,\nabla_x\phi_{2,\tau}(x,\eta)\right)
\Bigr) \\
&=a_{12,\tau;1,t}(x,\xi,\eta).
\end{align*}
From \eqref{F12}, we obtain
\begin{align*}
&\widetilde{a}_{12,\tau;1}
\left(t,x,
\nabla_x\phi_{1,\tau}(x,\xi),
\nabla_x\phi_{2,\tau}(x,\eta)
\right)  \\
&=\iint_{\mathbb{R}^{2n}}
e^{2\pi i\nabla_x\phi_{1,\tau}(x,\xi)\cdot u
+2\pi i\nabla_x\phi_{2,\tau}(x,\eta)\cdot v}
\frac{m_{\tau}(t,x,u,v)}
{(1+4\pi^2|u|^2+4\pi^2|v|^2)^N} \,du\,dv.
\end{align*}
Combining the previous identity with
\eqref{F12}, we conclude that,
\begin{align*}
a_{12,\tau;1,t}(x,\xi,\eta)
&=\iint_{\mathbb{R}^{2n}}
e^{2\pi i\nabla_x\phi_{1,\tau}(x,\xi)\cdot u
+2\pi i\nabla_x\phi_{2,\tau}(x,\eta)\cdot v}
\frac{m_{\tau}(t,x,u,v)}
{(1+4\pi^2|u|^2+4\pi^2|v|^2)^N}
\,du\,dv.
\end{align*}
This completes the proof.
\end{proof}

\begin{lemma} \label{Rel}
 $\forall x \in \mathcal{K} \text{ and }\xi,\eta \in \mathbb{R}^n\setminus \{0\}$, $|\xi|+|\eta| \approx |\Theta(x,\xi,\eta)|.$
\end{lemma}
\begin{proof}
  Consider the mapping $(x,\xi,\eta) \mapsto (x, \Theta(x,\xi,\eta)),$ where $$\Theta:=( \Psi_1(x,\xi) ,\Psi_2(x,\eta)).$$  Since $\Psi_1(x,\xi)$ and $\Psi_2(x,\eta)$ are positively homogeneous of degree $1,$ we have that $$|\Psi_1(x,\xi)|\approx |\xi|,\, |\Psi_2(x,\eta)|\approx |\eta|,$$ uniformly on $x \in \mathcal{K}$. Now, $$|\Theta(x,\xi,\eta)|= (|\Psi_1(x,\xi)|^2+|\Psi_2(x,\eta)|^2)^{1/2}.$$ 
  
  Using the fact that  $\forall a,b \in \mathbb{R},\,\dfrac{1}{\sqrt{2}}(a+b)\le (a^2+b^2)^{1/2}\le a+b, $ we have that $$\dfrac{1}{\sqrt{2}}(|\xi|+|\eta|)\le (|\xi|^2+|\eta|^2)^{1/2}\le |\xi|+|\eta|.$$ Hence, $|\xi|+|\eta| \approx |\Theta(x,\xi,\eta)|.$ The proof is complete.
\end{proof}
\begin{remark}
    To prove Lemma \ref{Lemma: Multi}, we state the above Lemma \ref{Rel}, which guarantees that $|\xi|+|\eta| \approx |\Theta(x,\xi,\eta)|$ is uniformly smooth for the spatial variable $x$ under the change of variables and behaves in the same way as the original frequency variables.
\end{remark}

\begin{lemma}
\label{Lemma: Multi} Assume that simultaneously $\phi_1$ and $\phi_2$ satisfy the Melin-Sj\"ostrand condition for some $\tau\in \mathbb{R},$ and are phases of {\it strictly positive type}.
    For any multi-indices $\alpha, \beta,$ and $\gamma$, we have
\begin{align*}
&\sup_{0<t<1}\sup_{x\in \mathcal{K}}
\sup_{\xi,\eta\in\mathbb{R}^n\setminus\{0\}}\left|
\partial_\xi^\alpha \partial_\eta^\beta \partial_x^\gamma
\tilde{a}_{12,\tau;1}(t, x,\xi,\eta)
\right| \lesssim_{\alpha,\beta,\gamma}
1.
\end{align*} In other words,
$$
\widetilde a_{12,\tau;1}
\in
L^\infty\bigl((0,1);S^0_{1,0}\bigr),
$$
uniformly with respect to the parameter $t$.
\end{lemma}
\begin{proof}
Consider
$$\tilde{a}_{12,\tau;1}(t,x,\xi,\eta)
=a_{12,\tau;1,t}
\bigl(x,\xi_x,\eta_x), \, (\xi_x,\eta_x
):=(\Psi_{1}(x,\xi),\Psi_{2}(x,\eta)
),$$
where \begin{align}
a_{12,\tau;1,t}(x,\xi,\eta)
& =t^{m}a_{12,\tau;1}
\left(x,\frac{\xi_x}{t},\frac{\eta_x}{t}\right)
\left|\nabla_x\phi_{1,\tau}(x,\xi_x)\right|^{-m}
\nonumber\\&\quad\times\widehat{\psi}
\bigl(\nabla_x\phi_{1,\tau}(x,\xi_x)\bigr)
\widehat{\theta}
\bigl(\nabla_x\phi_{2,\tau}(x,\eta_x)\bigr).
\end{align}
To make the notation simple, and in view of Lemma \eqref{Rel}, we will drop the sub-index $x$ in $(\xi_x,\eta_x),$ since  $(\xi,\eta)\mapsto(\xi_x,\eta_x)$ is a smooth local diffeomorphism, and we make the identification $(\xi,\eta)\equiv (\xi_x,\eta_x).$ From now, let us consider this identification and then let
$$
\widetilde a_{12,\tau;1}(t,x,\xi,\eta)
=
t^m
a_{12,\tau;1}
\left(x,\frac{\xi}{t},\frac{\eta}{t}\right)
B_1(x,\xi)B_2(x,\eta),
$$
where
$$
B_1(x,\xi)
=
\left|\nabla_x\phi_{1,\tau}(x,\xi)\right|^{-m}
\widehat{\psi}
\left(\nabla_x\phi_{1,\tau}(x,\xi)\right),
$$
and
$$
B_2(x,\eta)
=
\widehat{\theta}
\left(\nabla_x\phi_{2,\tau}(x,\eta)\right).
$$
We assume that the phases are positively homogeneous of degree one in
the frequency variables and satisfy
$$
\left|
\partial_x^\gamma\partial_\xi^\delta
\phi_{j,\tau}(x,\xi)
\right|
\leq
C_{\gamma,\delta}|\xi|^{1-|\delta|},
\qquad j=1,2,
$$
together with
$$
\left|\nabla_x\phi_{j,\tau}(x,\xi)\right|
\asymp |\xi|.
$$
Since
$$
\operatorname{supp}\widehat{\psi}
\subset
\left\{
\zeta\in\mathbb{R}^n:
\frac12\leq|\zeta|\leq2
\right\},
$$
the factor $B_1(x,\xi)$ is supported where
$$
\frac12
\leq
\left|\nabla_x\phi_{1,\tau}(x,\xi)\right|
\leq2.
$$
Consequently,
$$
c_1\leq|\xi|\leq c_2
$$
on the support of $B_1(x,\xi)$, for suitable positive constants $c_1,c_2$.
Similarly, since
$$
\operatorname{supp}\widehat{\theta}
\subset
\left\{
\zeta\in\mathbb{R}^n:
|\zeta|\leq\frac14
\right\},
$$
the factor $B_2(x,\eta)$ is supported where
$$
\left|\nabla_x\phi_{2,\tau}(x,\eta)\right|
\leq\frac14.
$$
Hence
$$
|\eta|\leq C
$$
on the support of $B_2(x,\eta)$. Therefore, on the support of
$\widetilde a_{12,\tau;1}$,
$$
|\xi|\asymp1,
\qquad
|\eta|\lesssim1.
$$
Since the factor $1-\chi$ eliminates a neighborhood of
$(\xi,\eta)=(0,0)$, we also have
$$
\left|
\left(\frac{\xi}{t},\frac{\eta}{t}\right)
\right|
\asymp t^{-1},
\qquad 0<t<1.
$$
Because $a(x,\xi,\eta)\in S^m_{1,0}$ and the cutoff functions are smooth, we have
$$
\left|
\partial_x^\gamma
\partial_\xi^{\alpha'}
\partial_\eta^{\beta'}
a_{12,\tau;1}(x,\xi,\eta)
\right|
\leq
C_{\alpha',\beta',\gamma}
\left\langle(\xi,\eta)\right\rangle
^{m-|\alpha'|-|\beta'|}.
$$
Using the chain rule, we obtain
$$
\begin{aligned}
&
\left|
\partial_x^\gamma
\partial_\xi^{\alpha'}
\partial_\eta^{\beta'}
a_{12,\tau;1}
\left(x,\frac{\xi}{t},\frac{\eta}{t}\right)
\right|
\\
&\qquad\leq
C_{\alpha',\beta',\gamma}
t^{-|\alpha'|-|\beta'|}
\left\langle
\left(\frac{\xi}{t},\frac{\eta}{t}\right)
\right\rangle
^{m-|\alpha'|-|\beta'|}.
\end{aligned}
$$

On the support under consideration,
$$
\left\langle
\left(\frac{\xi}{t},\frac{\eta}{t}\right)
\right\rangle
\asymp t^{-1}.
$$
Therefore,
$$
\left|
\partial_x^\gamma
\partial_\xi^{\alpha'}
\partial_\eta^{\beta'}
a_{12,\tau;1}
\left(x,\frac{\xi}{t},\frac{\eta}{t}\right)
\right|
\leq
C_{\alpha',\beta',\gamma}t^{-m}.
$$
Consequently,
$$
\left|
t^m
\partial_x^\gamma
\partial_\xi^{\alpha'}
\partial_\eta^{\beta'}
a_{12,\tau;1}
\left(x,\frac{\xi}{t},\frac{\eta}{t}\right)
\right|
\leq
C_{\alpha',\beta',\gamma}.
$$
We next consider the exponential factor. Set
$$
E_t(x,\xi,\eta)
=
\exp\left(
-2\pi i(\tau-i)
\left[
\operatorname{Im}
\phi_{1,\tau}\left(x,\frac{\xi}{t}\right)
+
\operatorname{Im}
\phi_{2,\tau}\left(x,\frac{\eta}{t}\right)
\right]
\right).
$$
Since
$$
-i(\tau-i)=-1-i\tau,
$$
we have
$$
\begin{aligned}
|E_t(x,\xi,\eta)|
&=
\exp\left(
-2\pi
\left[
\operatorname{Im}
\phi_{1,\tau}\left(x,\frac{\xi}{t}\right)
+
\operatorname{Im}
\phi_{2,\tau}\left(x,\frac{\eta}{t}\right)
\right]
\right).
\end{aligned}
$$
By hypothesis, there exists
$c>0$ such that
$$
\operatorname{Im}\phi_{j,\tau}(x,\zeta)
\geq c|\zeta|,
\qquad
x\in\mathcal K,
\quad
\zeta\neq0.
$$
Thus,
$$
\operatorname{Im}
\phi_{1,\tau}\left(x,\frac{\xi}{t}\right)
\geq
\frac{c}{t}|\xi|.
$$
Since $|\xi|\asymp1$ on the support of $\widehat{\psi}$, it follows that
$$
|E_t(x,\xi,\eta)|
\leq
e^{-c_0/t}
$$
for some $c_0>0$.
When derivatives fall on $E_t$, the chain rule produces powers of
$t^{-1}$. Hence, for arbitrary multi-indices
$\alpha',\beta',\gamma'$, there exist constants
$C_{\alpha',\beta',\gamma'}>0$ and
$N_{\alpha',\beta',\gamma'}\geq0$ such that
$$
\left|
\partial_\xi^{\alpha'}
\partial_\eta^{\beta'}
\partial_x^{\gamma'}
E_t(x,\xi,\eta)
\right|
\leq
C_{\alpha',\beta',\gamma'}
t^{-N_{\alpha',\beta',\gamma'}}
e^{-c_0/t}.
$$
Since exponential decay dominates every negative power of $t$, for
every $L>0$,
$$
t^{-N}e^{-c_0/t}
\leq
C_{N,L}t^L,
\qquad 0<t<1.
$$
Therefore,
$$
\left|
\partial_\xi^{\alpha'}
\partial_\eta^{\beta'}
\partial_x^{\gamma'}
E_t(x,\xi,\eta)
\right|
\leq
C_{\alpha',\beta',\gamma',L}t^L.
$$
The factors $B_1(x,\xi)$ and $B_2(x,\eta)$ satisfy uniform estimates. Indeed, by
the chain rule and the symbol estimates for the phases,
$$
\left|
\partial_x^{\gamma'}
\partial_\xi^{\alpha'}B_1(x,\xi)
\right|
\leq
C_{\alpha',\gamma'},
$$
and
$$
\left|
\partial_x^{\gamma'}
\partial_\eta^{\beta'}B_2(x,\eta)
\right|
\leq
C_{\beta',\gamma'}.
$$
We now apply the Leibniz rule to
$$
\widetilde a_{12,\tau;1}
=
t^m
E_t
\,
a_{12,\tau;1}
\left(x,\frac{\xi}{t},\frac{\eta}{t}\right)
B_1(x,\xi)B_2(x,\eta).
$$
Every term in
$$
\partial_\xi^\alpha
\partial_\eta^\beta
\partial_x^\gamma
\widetilde a_{12,\tau;1}
$$
is a product of a derivative of $E_t$, a derivative of the rescaled
symbol, a derivative of $B_1$, and a derivative of $B_2$.

The factor $t^m$ cancels the factor $t^{-m}$ arising from the
rescaled symbol. The remaining powers of $t^{-1}$, generated by
differentiating the phase, are absorbed by the exponential factor
$e^{-c_0/t}$. Hence, for every $L>0$,
$$
\left|
\partial_\xi^\alpha
\partial_\eta^\beta
\partial_x^\gamma
\widetilde a_{12,\tau;1}(t,x,\xi,\eta)
\right|
\leq
C_{\alpha,\beta,\gamma,L}t^L.
$$
In particular, taking $L=0$, we conclude that
$$
\sup_{0<t<1}
\sup_{x\in\mathcal K}
\sup_{\xi,\eta\neq0}
\left|
\partial_\xi^\alpha
\partial_\eta^\beta
\partial_x^\gamma
\widetilde a_{12,\tau;1}(t,x,\xi,\eta)
\right|
\leq
C_{\alpha,\beta,\gamma}.
$$
More precisely, the stronger estimate
$$
\left|
\partial_\xi^\alpha
\partial_\eta^\beta
\partial_x^\gamma
\widetilde a_{12,\tau;1}(t,x,\xi,\eta)
\right|
\leq
C_{\alpha,\beta,\gamma,L}t^L,
\qquad
\forall L\geq 0,
$$
holds uniformly for $0<t<1$ and $x\in\mathcal K$.
Thus,
$$
\widetilde a_{12,\tau;1}
\in
L^\infty\bigl((0,1);S^0_{1,0}\bigr),
$$
uniformly with respect to the parameter $t$.
\end{proof}

\begin{remark}
 On the support of $\hat{\psi},$ we have that $ \left|\nabla_{x}\phi_{1,\tau}(x,\Psi_1(x,\xi))\right|\asymp |\xi|,$  and on the support of $\hat{\theta},$ we have that  $\left|\nabla_{x}\phi_{2,\tau}(x,\Psi_2(x,\eta))\right|\asymp |\eta|.$ This implies that, all the negative powers of $|\xi|$ and $|\eta|$ are bounded on the support of $\hat{\psi}$ and $\hat{\theta}$ respectively.
\end{remark}
\begin{remark}
Lemma~\ref{Lemma: Multi} yields that, for all multi-indices $\gamma$,
\begin{align}\label{reg:2}
\sup_{0<t<1}\;\sup_{x\in K}\;\sup_{u,v\in\mathbb{R}^{n}}
\left|\partial_{x}^{\gamma} m_{\tau}(t,x,u,v)\right|<+\infty.
\end{align} The case $\gamma=0$ will be fundamental in our further analysis.  Indeed, observe that in the proof of Lemma \ref{Lemma: Multi}, we have derived the representation
\begin{align*}
    m_{\tau}(t,x,u,v)
    = \int_{\mathbb{R}^n}\int_{\mathbb{R}^n} e^{-2\pi i\xi\cdot u- 2\pi i\eta \cdot v} [(1-\Delta_{\xi}-\Delta_{\eta})^N\widetilde{a}_{12,\tau;1}(t,x,\xi,\eta)]\,d\xi\,d\eta.
\end{align*}Since derivatives in $\Xi:=(\xi,\eta)$ of the symbol $\widetilde{a}_{12,\tau;1}(t,\cdot, \cdot, \cdot):=\widetilde{a}_{12,\tau;1}(t,x,\xi,\eta)$ are bounded and it is compactly supported in $\Xi,$ we have that
\begin{align*}
\left| m_{\tau}(t,x,u,v)\right| &\leq \int_{\mathbb{R}^n}\int_{\mathbb{R}^n} | [(1-\Delta_{\xi}-\Delta_{\eta})^N\widetilde{a}_{12,\tau;1}(t,x,\xi,\eta)]|\,d\xi\,d\eta\\
&\lesssim \sup_{|\beta|\leq 2N}\|\partial_\Xi^\beta \widetilde{a}_{12,\tau;1}(t,\cdot, \cdot, \cdot) \|_{L^\infty(\Omega)}\textnormal{measure}(\Omega),
\end{align*} where the compact set $\Omega$ is the support of $\widetilde{a}_{12,\tau;1}$ in the $(\xi,\eta)$-variable. The proof of \eqref{reg:2} now is a consequence of the previous inequalities. 
\end{remark}
\begin{proposition}
Let $N\in \mathbb{N}.$ For this $N$ fixed, there exists $m_{\tau}=m_{\tau}^{(N)}$ be the function defined in \eqref{m:tau} depending on $N;$ such that the bilinear Fourier integral operator $T^{{\phi_{1,\tau}},{\phi_{2,\tau}}}_{a_{12,\tau;1}},$ in (\ref{a}) admits the following representation 
\begin{align}
&T^{{\phi_{1,\tau}},{\phi_{2,\tau}}}_{a_{12,\tau;1}}(\mu(D)f,\mu(D)g)(x)\label{D}
\nonumber\\&=\displaystyle{ \int\limits_{0}^1\int_{\mathbb{R}^n}\int_{\mathbb{R}^n}}T^{\phi_{1,\tau}}_{\nu_{1}^{t,u}} f(x)
T^{\phi_{2,\tau}}_{\mu_{1}^{t,v}} g(x)
\frac{m_{\tau}(t,x,u,v)}{(1+4\pi^2|u|^2+4\pi^2|v|^2)^N}
\,du\,dv\,\frac{dt}{t},
\end{align}
where the symbols $\nu_1^{t,u}$ and $\mu_1^{t,v}$ are given by 
\[
\nu_1^{t,u}(x,\xi)
=\mu(\xi)
\left|\nabla_x\phi_{1,\tau}(x,\xi)\right|^{m}
e^{ 2\pi it\nabla_x\phi_{1,\tau}(x,\xi)\cdot u}
\hat{\psi}\left(t\nabla_x\phi_{1,\tau}(x,\xi)\right)
\chi(x),
\]
and
\[
\mu_1^{t,v}(x,\eta)
=\mu(\eta)
e^{2\pi it\nabla_x\phi_{2,\tau}(x,\eta)\cdot v}
\hat{\theta}\left(t\nabla_x\phi_{2,\tau}(x,\eta)\right)
\chi(x),
\]
where $\chi\in C_0^\infty(\mathbb{R}^n)$ and is equal to $1,$ for $x\in\mathcal{K}.$
\end{proposition}
\begin{proof}
    From (\ref{a}), we have 
\begin{align}\
&T_{a_{12,\tau;1}}^{\phi_{1,\tau}, \phi_{2,\tau}}
\bigl(\mu(D)f,\mu(D)g\bigr)(x)
\nonumber\\&=\int_{0}^{\infty}
\iint a_{12,\tau;1,t}(x,t\xi,t\eta)
\Bigl(\widehat{\psi}\bigl(t\nabla_x\phi_{1,\tau}(x,\xi)\bigr)\nonumber\\&\quad\times\left|\nabla_x\phi_{1,\tau}(x,\xi)\right|^{m}\mu(\xi)\widehat{f}(\xi)e^{i\phi_{1,\tau}(x,\xi)}\Bigr)
\nonumber\\&\quad\times
\Bigl(\widehat{\theta}\bigl(t\nabla_x\phi_{2,\tau}(x,\eta)\bigr)\mu(\eta)\widehat{g}(\eta)e^{i\phi_{2,\tau}(x,\eta)}\Bigr)\frac{dt}{t}\,d\xi\,d\eta.
\end{align}
 Let define 
\[
\nu_1^{t,u}(x,\xi)
=\mu(\xi)
\left|\nabla_x\phi_{1,\tau}(x,\xi)\right|^{m}
e^{ 2\pi it\nabla_x\phi_{1,\tau}(x,\xi)\cdot u}
\hat{\psi}\left(t\nabla_x\phi_{1,\tau}(x,\xi)\right)
\chi(x),
\]
and
\[
\mu_1^{t,v}(x,\eta)
=\mu(\eta)
e^{ 2\pi it\nabla_x\phi_{2,\tau}(x,\eta)\cdot v}
\hat{\theta}\left(t\nabla_x\phi_{2,\tau}(x,\eta)\right)
\chi(x),
\]
as symbols for the Fourier integral operators $T_{\nu^{t,u}_{1}}^{\phi_{1,\tau}}$ and $T_{\mu^{t,v}_{1}}^{\phi_{2,\tau}}$ respectively. The introduction of $e^{ 2\pi it\nabla_x\phi_{1,\tau}(x,\xi)\cdot u}$ and $e^{ 2\pi it\nabla_x\phi_{2,\tau}(x,\eta)\cdot v}$ will be explained later. For $f,g \in C_0^{\infty}(\mathbb{R}^n),$ we define the following Fourier integral operators as follows;
\begin{align}\label{2}
T_{\nu^{t,u}_{1}}^{\phi_{1,\tau}}f(x)
=
\int_{\mathbb R^n}
\nu^{t,u}_{1}(x,\xi)
\widehat f(\xi)
e^{i\phi_{1,\tau}(x,\xi)}
\,d\xi,
\end{align}
and
\begin{align}
T_{\mu^{t,v}_{1}}^{\phi_{2,\tau}}g(x)\label{3}
=
\int_{\mathbb R^n}
\mu^{t,v}_{1}(x,\eta)
\widehat g(\eta)
e^{i\phi_{2,\tau}(x,\eta)}
\,d\eta.
\end{align}
We know that,
\begin{align}
a_{12,\tau;1,t}(x,\xi,\eta)
=\iint_{\mathbb{R}^{2n}}
e^{2\pi i\nabla_x\phi_{1,\tau}(x,\xi)\cdot u
+2\pi i\nabla_x\phi_{2,\tau}(x,\eta)\cdot v}
\frac{m_{\tau}(t,x,u,v)}
{(1+4\pi^2|u|^2+4\pi^2|v|^2)^N}
\,du\,dv \label{AL1}
\end{align}
From (\ref{AL1}), we deduce that 
\begin{align}
&a_{12,\tau;1,t}(x,t\xi,t\eta)\nonumber\\&
=\iint_{\mathbb{R}^{2n}}
e^{2\pi i\nabla_x\phi_{1,\tau}(x,t\xi)\cdot u
+2\pi i\nabla_x\phi_{2,\tau}(x,t\eta)\cdot v}
\frac{m_{\tau}(t,x,u,v)}
{(1+4\pi^2|u|^2+4\pi^2|v|^2)^N}
\,du\,dv \nonumber\\&= 
\iint_{\mathbb{R}^{2n}}
e^{2\pi it\nabla_x\phi_{1,\tau}(x,\xi)\cdot u
+2\pi it\nabla_x\phi_{2,\tau}(x,\eta)\cdot v}
\frac{m_{\tau}(t,x,u,v)}
{(1+4\pi^2|u|^2+4\pi^2|v|^2)^N}
\,du\,dv
\label{AL2}
\end{align}
It is important to mention that, $\nabla_x\phi_{1,\tau}(x,t\xi)=t\nabla_x\phi_{1,\tau}(x,\xi)$ and $\nabla_x\phi_{2,\tau}(x,t\eta)=t\nabla_x\phi_{2,\tau}(x,\eta),$ since the phases are homogeneous of degree one in the frequency variables.
 Therefore, substituting  (\ref{AL2}) into (\ref{a}) yields
\begin{align}
&T_{a_{12,\tau;1}}^{\phi_{1,\tau}, \phi_{2,\tau}}\bigl(\mu(D)f,\mu(D)g\bigr)(x)\nonumber\\
&=\int_{0}^{1}\iint_{\mathbb{R}^{2n}} \Bigg[
\iint_{\mathbb{R}^{2n}}
e^{2\pi it\nabla_x\phi_{1,\tau}(x,\xi)\cdot u
+2\pi it\nabla_x\phi_{2,\tau}(x,\eta)\cdot v}
\frac{m_{\tau}(t,x,u,v)}
{(1+4\pi^2|u|^2+4\pi^2|v|^2)^N}
\,du\,dv
\Bigg]\nonumber\\
&\quad\times \Bigl(\widehat{\psi}\bigl(t\nabla_x\phi_{1,\tau}(x,\xi)\bigr) |\nabla_x\phi_{1,\tau}(x,\xi)|^m \mu(\xi)\widehat{f}(\xi)e^{i\phi_{1,\tau}(x,\xi)}\Bigr)\nonumber\\
&\quad\times \Bigl(\widehat{\theta}\bigl(t\nabla_x\phi_{2,\tau}(x,\eta)\bigr) \mu(\eta)\widehat{g}(\eta)e^{i\phi_{2,\tau}(x,\eta)}\Bigr) \, d\xi \, d\eta \, \frac{dt}{t}.
\end{align}
Interchanging the order of integration by Fubini's theorem and splitting $$e^{2\pi it\nabla_x\phi_{1,\tau}(x,\xi)\cdot u
+2\pi it\nabla_x\phi_{2,\tau}(x,\eta)\cdot v}$$ into $e^{2\pi it\nabla_x\phi_{1,\tau}(x,\xi)\cdot u}$ and $e^{2\pi it\nabla_x\phi_{2,\tau}(x,\eta)\cdot v}$, gives us
\begin{align}
&T_{a_{12,\tau;1}}^{\phi_{1,\tau}, \phi_{2,\tau}}\bigl(\mu(D)f,\mu(D)g\bigr)(x)\nonumber\\
&=\int_{0}^{1} \iint_{\mathbb{R}^{2n}} \frac{m_\tau(t,x,u,v)}{\left(1+4\pi^2|u|^2+4\pi^2|v|^2\right)^N}\nonumber\\
&\quad\times \left( \int_{\mathbb{R}^n} e^{2\pi it\nabla_x\phi_{1,\tau}(x,\xi)\cdot u}\widehat{\psi}\bigl(t\nabla_x\phi_{1,\tau}(x,\xi)\bigr) |\nabla_x\phi_{1,\tau}(x,\xi)|^m \mu(\xi)\widehat{f}(\xi)e^{i\phi_{1,\tau}(x,\xi)} \, d\xi \right)\nonumber\\
&\quad\times \left( \int_{\mathbb{R}^n} e^{2\pi it\nabla_x\phi_{2,\tau}(x,\eta)\cdot v}\widehat{\theta}\bigl(t\nabla_x\phi_{2,\tau}(x,\eta)\bigr) \mu(\eta)\widehat{g}(\eta)e^{i\phi_{2,\tau}(x,\eta)} \, d\eta \right) \, du \, dv \, \frac{dt}{t}. \label{addsteps}
\end{align}
Now, substituting (\ref{2}) and (\ref{3}) into (\ref{addsteps}), one has
\begin{align}
&T_{a_{12,\tau;1}}^{\phi_{1,\tau},\phi_{2,\tau}}
(\mu(D)f,\mu(D)g)(x)
\nonumber\\
&=
\int_0^1
\int_{\mathbb R^n}\int_{\mathbb R^n}
T_{\nu^{t,u}_{1}}^{\phi_{1,\tau}}f(x)
T_{\mu^{t,v}_{1}}^{\phi_{2,\tau}}g(x)
\frac{m_\tau(t,x,u,v)}
{\left(1+4\pi^2|u|^2+4\pi^2|v|^2\right)^N}
\,du\,dv\,\frac{dt}{t}.
\end{align}
The proof is complete.
\end{proof}

\begin{lemma}
    Let us consider the symbol
\begin{align*}
a_{12,\tau}(x,\xi,\eta)&:=e^{-2\pi i(\Im\phi_{1,\tau}(x,\xi)+\Im\phi_{2,\tau}(x,\eta))(\tau-i)}a(x,\xi,\eta).
\end{align*}If there exist $N\in \mathbb{N},$  such that
\begin{equation}
    \forall |\alpha|,|\beta|\leq N,\, \quad \sup_{(x,\xi,\eta)}|\partial_{x,\xi,\eta}^\alpha a_{12,\tau}(x,\xi,\eta)|\leq C_{\alpha}|\xi|^{m-|\alpha|}<\infty,\,|\xi|\geq 1, 
\end{equation} then the function $m_\tau$ is such that the following estimate holds
\begin{equation}
    |m_\tau(t,x,u,v)|\leq C_N,
\end{equation} uniformly in $0<t<1$ and $(x,\xi,\eta).$
    \end{lemma}
\begin{proof}This is a special case of Lemma 4.1 in \cite{rodriguez2011boundedness} since the symbol $a_{12,\tau}$  satisfies Kohn-Nirenberg estimates up to $N$-derivatives.    
\end{proof}

The operator $T^{\phi_{1,\tau}}_{\nu^{1}_{t,u}}$ has the following factorisation modulo
an error term 

\begin{equation}
T^{\phi_{1,\tau}}_{\nu^{1}_{t,u}}
=t^{-m_2}Q_t^{u}T^{\phi_{1,\tau}}_{\nu_{m_1}}
+t^{-m_2}E_t^{1},
\end{equation}
where
\begin{itemize}
\item $Q_t^{u}$ is a Fourier multiplier with symbol 
\[
\xi \mapsto |t\xi|^{m_2}\widehat{\psi}(t\xi)e^{2\pi it u\cdot\xi}
\in \mathscr{S}(\mathbb{R}^n).
\]
The function $\psi\in\mathscr{S}(\mathbb{R}^n)$ is chosen to be smooth,
whose Fourier transform is supported on the annulus
\[
\{\xi\in\mathbb{R}^n: \tfrac12\le |\xi|\le 2\}
\] and such that
\[
\forall t>0,\ \forall \xi\neq 0,~
\int_0^\infty |\widehat{\psi}(t\xi)|^2\frac{dt}{t}=1.
\]

\item The operator $T^{\phi_{1,\tau}}_{\nu_{m_1}}$ is a Fourier integral operator of order $m_1$
with phase $\phi_{1,\tau}$ and with symbol 
\begin{equation}
(x,\xi)\mapsto
\nu_{m_1}(x,\xi)
=\mu(\xi)\chi(x)
|\nabla_x\phi_{1,\tau}(x,\xi)|^{m_1}
\in S^{m_1}_{1,0}(\mathbb{R}^n_x\times\mathbb{R}^n_\xi).
\end{equation}
The operator $E_t^{1}$ has the same phase $\phi_{1,\tau}$ of $T^{\phi_{1,\tau}}_{\nu_{m_1}}$
and a symbol in the class
$
S^{m_1-(\frac12-\varepsilon)}_{1,0}(\mathbb{R}^n_x\times\mathbb{R}^n_\xi),
$
with semi-norms bounded by the product of $t^\varepsilon$,
$\varepsilon\in(0,\tfrac12)$, with a polynomial expression in $|u|$.
\end{itemize}

On the other hand, the operator $T^{\phi_{2,\tau}}_{\mu_{1}^{t,v}}$ admits the following
factorisation modulo an error term 

\[
t^{-m_2}T^{\phi_{2,\tau}}_{\mu_{1}^{t,v}}
=
P_t^{v}T^{\phi_{2,\tau}}_{\mu^{t}_{1}}
+
E_t^{2},
\]
where:
\begin{itemize}
 \item[ \textbullet~~(i)] the symbol of $P_t^{v}$ is given by
\[
\eta\mapsto \widehat{\theta}(t\eta)e^{2\pi itv\cdot\eta}
\in S^0_{1,0}(\mathbb{R}^n_\eta),
\]
\item[ \textbullet ~~(ii)] and $T^{\phi_{2,\tau}}_{\mu^{t}_{1}}$ is a Fourier integral operator of order
$m_2<0$, with phase function $\phi_{2,\tau}$, and with symbol
\[
(x,\eta)\mapsto
\mu^{t}_{1}(x,\eta)
=
\mu(\eta)\chi(x)\widehat{\theta}_1(t\nabla_x\phi_{2,\tau}(x,\eta))t^{-m_2}
\in S^{m_2}_{1,0}(\mathbb{R}^n_x\times\mathbb{R}^n_\eta),
\]

uniformly in $t\in(0,1]$.
Here $E_t^{2}$ is a Fourier integral operator with phase $\phi_{2,\tau}$
and a symbol in the class
$
S^{m_2-(\frac12-\varepsilon)}_{1,0}(\mathbb{R}^n_x\times\mathbb{R}^n_\eta)$,
with semi-norms that behave polynomially in $|v|$.
\end{itemize}

\subsection{Bounds using H\"older's inequality}
In this subsection, we establish the boundedness estimates for the high-frequency portion of the bilinear Fourier integral operator with complex-valued phases. 

After the decomposition introduced above, the operator can be represented in
terms of products of linear Fourier integral operators together with the kernel
$m_{\tau}(t,x,u,v)$. This reduces the bilinear estimate
to linear estimates by means of H\"older's inequality. Let consider
\begin{align*}
&T^{\phi_{1,\tau},\phi_{2,\tau}}_{a_{12,\tau}}\bigl(\mu(D)f,\mu(D)g\bigr)(x)\\
&= \int_{0}^{1} \int_{\mathbb{R}^n}\int_{\mathbb{R}^n}
T^{\phi_{1,\tau}}_{\nu^{1}_{t,u}} f(x)
T^{\phi_{2,\tau}}_{\mu^{1}_{t,v}} g(x)\,
\frac{m_{\tau}(t,x,u,v)}{(1+4\pi^2|u|^2+4\pi^2|v|^2)^N}
\,du\,dv\,\frac{dt}{t}.
\end{align*}
We want to prove that,
\begin{align*}
\bigl\|T^{\phi_{1,\tau},\phi_{2,\tau}}_{a_{12,\tau}}(\mu(D)f,\mu(D)g)\bigr\|_{L^p(\mathbb{R}^n_x)} \lesssim \|f\|_{L^{p_1}(\mathbb{R}^n_x)}
\|g\|_{L^{p_2}(\mathbb{R}^n_x)}, \text{ for } 1<p_1,p_2< \infty.
\end{align*}
\begin{proof}
Taking the $L^p$-norm and applying H\"older's inequality in the variables
$u$ and $v$, we obtain
\begin{align*}
&\bigl\|T^{\phi_{1,\tau},\phi_{2,\tau}}_{a_{12,\tau}}(\mu(D)f,\mu(D)g)\bigr\|_{L^p(\mathbb{R}^n_x)}\\
&\le \int_0^1 \int_{\mathbb{R}^n}\int_{\mathbb{R}^n}
\bigl\|T^{\phi_{1,\tau}}_{\nu^{1}_{t,u}} f(x) T^{\phi_{2,\tau}}_{\mu^{1}_{t,v}} g(x)\bigr\|_{L^p(\mathbb{R}^n_x)}
\frac{\sup_{x}|m_{\tau}(t,x,u,v)|}{(1+4\pi^2|u|^2+4\pi^2|v|^2)^N}
\,du\,dv\,\frac{dt}{t}
\\&\le \int_0^1 
 \bigl\|T^{\phi_{1,\tau}}_{\nu^{1}_{t,u}} f(x)\bigr\|_{L^{p_1}(\mathbb{R}^n_x)} \bigl\|T^{\phi_{2,\tau}}_{\mu^{1}_{t,v}} g(x)\bigr\|_{L^{p_2}(\mathbb{R}^n_x)}\\& \quad\times
\left( \int_{\mathbb{R}^n}\int_{\mathbb{R}^n}
\frac{\sup_{x}|m_{\tau}(t,x,u,v)|}{(1+4\pi^2|u|^2+4\pi^2|v|^2)^N}
\,du\,dv \right)
\frac{dt}{t}
\\&\lesssim  
\bigl\|T^{\phi_{1,\tau}}_{\nu^{1}_{t,u}} f\bigr\|_{L^{p_1}(\mathbb{R}^n_x)}
\bigl\|T^{\phi_{2,\tau}}_{\mu^{1}_{t,v}} g\bigr\|_{L^{p_2}(\mathbb{R}^n_x)}\\& \quad\times  \left( \int_0^1\int_{\mathbb{R}^n}\int_{\mathbb{R}^n}
\frac{\sup_{x}|m_{\tau}(t,x,u,v)|}{(1+4\pi^2|u|^2+4\pi^2|v|^2)^N}
\,du\,dv \right)
\,\frac{dt}{t}
\\&\lesssim \|f\|_{L^{p_1}(\mathbb{R}^n_x)}
\|g\|_{L^{p_2}(\mathbb{R}^n_x)}
\int_0^1 
\int_{\mathbb{R}^n}\int_{\mathbb{R}^n}
\frac{\sup_{x}|m_{\tau}(t,x,u,v)|}{(1+4\pi^2|u|^2+4\pi^2|v|^2)^N}
\,du\,dv\,\frac{dt}{t}.
\end{align*}
From (\ref{inf}), we deduce that $\sup_{x}|m_{\tau}(t,x,u,v)|< +\infty,$ so for $N$ large enough, we have that
\begin{align*}
    \int^1_0\int_{\mathbb{R}^n}\int_{\mathbb{R}^n}
\frac{\sup_{x}|m_{\tau}(t,x,u,v)|}{(1+4\pi^2|u|^2+4\pi^2|v|^2)^N}
\,du\,dv\, \frac{dt}{t}< \infty.
\end{align*}
This implies that, 
\begin{align*}
\bigl\|T^{\phi_{1,\tau},\phi_{2,\tau}}_{a_{12,\tau}}(\mu(D)f,\mu(D)g)\bigr\|_{L^p(\mathbb{R}^n_x)} \lesssim \|f\|_{L^{p_1}(\mathbb{R}^n_x)}
\|g\|_{L^{p_2}(\mathbb{R}^n_x)}.
\end{align*}
The proof is complete.
\end{proof}
We have used the fact that for $L^{p}$ boundedness of Fourier integral operators, we have that
\begin{align*}
\|T^{\phi_{1,\tau}}_{\nu^{1}_{t,u}} f(\cdot)\|_{L^{p_1}}
&\lesssim \|f\|_{L^{p_1}}, 
\end{align*}
and
\begin{align*}
\|T^{\phi_{2,\tau}}_{\mu^{1}_{t,v}} g(\cdot)\|_{L^{p_2}}
&\lesssim \|g\|_{L^{p_2}},
\end{align*}
and H\"older inequality holds when
\begin{align*}
    \frac{1}{p}=\frac{1}{p_1}+\frac{1}{p_2}, \text{ for } 1< p_1,p_2 < \infty \text{ and } 1< p< \infty.
\end{align*}
In the setting of linear Fourier integral operators, the main condition to ensure boundedness is when $1/p=1/p_1+1/p_2.$ This condition satisfies $L^{p_1}(\mathbb{R}^n)\times L^{p_2}(\mathbb{R}^n)\longrightarrow L^{p}(\mathbb{R}^n)$ estimates to provide estimates for the high frequency portion of the Fourier integral operators. This indicates that whenever $T^{\phi_{1,\tau}}_{\nu^{1}_{t,u}} f(x)\in L^{p_1}(\mathbb{R}^n)$ and $T^{\phi_{2,\tau}}_{\mu^{1}_{t,v}} g(x)\in L^{p_2}(\mathbb{R}^n),$ then their para-product $T^{\phi_{1,\tau}}_{\nu^{1}_{t,u}} f(x)T^{\phi_{2,\tau}}_{\mu^{1}_{t,v}} g(x)\in L^{p}(\mathbb{R}^n)$. Boundedness fails when $p_1=p_2=1$ and when $p_1=p_2=\infty,$ resulting in boundedness properties only in the range $1<p_1,p_2 < \infty.$ In particular, when the order of the bilinear Fourier integral operator, $m=0,$ we have that
\begin{align*}
\bigl\|T^{\phi_{1,\tau},\phi_{2,\tau}}_{a_{12,\tau}}(\mu(D)f,\mu(D)g)\bigr\|_{L^1(\mathbb{R}^n_x)} \lesssim \|f\|_{L^{2}(\mathbb{R}^n_x)}
\|g\|_{L^{2}(\mathbb{R}^n_x)}.
\end{align*}
\begin{remark}
    This indicates that, the norm of a bilinear Fourier integral operator is controlled by the product of the norm of its individual linear Fourier integral operator.
\end{remark}
In order to prove our main Theorem \ref{Main:1:th},  we split our analysis into low frequencies and high frequencies. The proof takes into account Lemma  \ref{Rel} and Lemma \ref{Lemma: Multi}.

\subsection{ Analysis of the low  frequencies }
Here, we introduce the low-frequency portion of the bilinear Fourier integral operator, which is reduced to linear Fourier integral operators by freezing one frequency variable and by applying H\"older's inequality. 
For the analysis of low-frequencies, it suffices to prove that a bilinear Fourier integral operator
\begin{align*}
T_{a_{12,\tau}}^{\phi_1,\tau, \phi_2,\tau}(f,g)(x)= \displaystyle{ \int\limits_{\mathbb{R}^n}\int\limits_{\mathbb{R}^n}} e^{2\pi i \phi_{1,\tau}(x,\xi)}\hat{g}(\eta)a_{12,\tau}(x,\xi,\eta)e^{2\pi i \phi_{2,\tau}(x,\eta)}\hat{f}(\xi)d\xi d\eta,
\end{align*}
with a symbol $a_{12,\tau}(x,\xi,\eta)$ compactly supported in $\xi$ or in $\eta$ satisfies the estimate
\begin{align*}
\|T^{\phi_{1,\tau}, \phi_{2,\tau}}_{a_{12,\tau}}(f,g)\|_{L^p(\mathbb{R}^n_x))} \lesssim
\|f\|_{L^{p_1}(\mathbb{R}^n_x)}
\|g\|_{L^{p_2}(\mathbb{R}^n_x)}, 
\text{ for } 1/p=1/p_1+1/p_2 \text{ and } \\ 1<p_1,p_2<\infty.
\end{align*}
The bilinear operator with complex valued phase is defined as
\begin{align*}
   &T_{a_{12,\tau}}^{\phi_1,\tau, \phi_2,\tau}(f,g)(x) \\&=\displaystyle{ \int\limits_{\mathbb{R}^n}\int\limits_{\mathbb{R}^n}} e^{2\pi i (\phi_{1,\tau}(x,\xi)+\phi_{2,\tau}(x,\eta))}e^{-2\pi i(\Im\phi_1(x,\xi)+\Im\phi_2(x,\eta))(\tau-i)}a(x,\xi,\eta)\hat{f}(\xi)\hat{g}(\eta)d\xi d\eta\\&
   = \displaystyle{ \int\limits_{\mathbb{R}^n}\int\limits_{\mathbb{R}^n}} e^{2\pi i \phi_{1,\tau}(x,\xi)}\hat{g}(\eta)a_{12,\tau}(x,\xi,\eta)e^{2\pi i \phi_{2,\tau}(x,\eta)}\hat{f}(\xi)d\xi d\eta,
\end{align*}
with the symbol $a_{12,\tau}(x,\xi,\eta)=e^{-2\pi i(\Im\phi_1(x,\xi)+\Im\phi_2(x,\eta))(\tau-i)}a(x,\xi,\eta)$ and complex-valued phase function $\phi_{12, \tau}= \phi_{1,\tau}(x,\xi)+\phi_{2,\tau}(x,\eta).$ Suppose  $\psi$ is a smooth cut-off function that is equal to $1$ on the $\xi$-support of $a_{12,\tau}$. We can rewrite  $T_{a_{12,\tau}}^{\phi_1,\tau, \phi_2,\tau}(f,g)$ as
\begin{align}
     T_{a_{12,\tau}}^{\phi_{1,\tau}, \phi_{2,\tau}}(f,g)(x)&= \int_{\mathbb{R}^n} \bigg(\int_{\mathbb{R}^n} a_{12,\tau}(x,\xi,\eta) \hat{g}(\eta)e^{2\pi i \phi_{2,\tau}(x,\eta)}\,d\eta\bigg)\nonumber\\& \quad \times \psi(\xi)\hat{f}(\xi)e^{2\pi i \phi_{1,\tau}(x,\xi)}\, d\xi\nonumber\\&
     =\int_{\mathbb{R}^n} \mathfrak{a}_{12,g}(x,\xi)\hat{f}(\xi)e^{2\pi i \phi_{1,\tau}(x,\xi)}\, d\xi=T^{\phi_{1,\tau}}_{\mathfrak{a}_g}(f)(x),\label{a123}
\end{align}
where 
\begin{align}
\mathfrak{a}_{12,g}(x,\xi)= \psi(\xi)\int_{\mathbb{R}^n} a_{12,\tau}(x,\xi,\eta) \hat{g}(\eta)e^{2\pi i \Phi_{2,\tau}(x,\eta)}\,d\eta=T^{\phi_2,\tau}_{\psi(\xi)a(\cdot,\xi, \cdot)}\, g(x) \label{a12g} . \end{align}  From \eqref{a123}, we have
\begin{align*}
    T^{\phi_{1,\tau}}_{\mathfrak{a}_g}(f)(x)&= \int_{\mathbb{R}^n} e^{2\pi i \phi_{1,\tau}(x,\xi)}\mathfrak{a}_{12,g}(x,\xi)\zeta(\xi)\hat{f}(\xi)\, d\xi\\&
    =C_d \int_{\mathbb{R}^n} e^{2\pi i \phi_{1,\tau}(x,\xi)}\int_{\mathbb{R}^n} \widehat{\mathfrak{a}}_g(x,\lambda)\, e^{2 \pi i\xi \cdot \lambda}\, d\lambda \zeta(\xi)\hat{f}(\xi)\, d\xi\\&
    =C_d \int_{\mathbb{R}^n} \widehat{\mathfrak{a}}_g(x,\lambda)\,\int_{\mathbb{R}^n} e^{2\pi i \phi_{1,\tau}(x,\xi)}\widehat{f(\cdot -\lambda)}(\xi)\zeta(\xi)\, d\xi\, d\lambda\\&
=C_d\int_{\mathbb{R}^n}\widehat{\mathfrak{a}}_g(x,\lambda)\,T^{\phi_1,\tau}_{\zeta}\bigl(f(\cdot - \lambda)\bigr)\, d\lambda
\end{align*}
Let us consider a function $\chi \in C_0^\infty$  such that  $\chi = 1$ on  $\operatorname{supp}_x(\widehat{a}_g)$. Then,  
\begin{align*}
     T^{\phi_{1,\tau}}_{\mathfrak{a}_g}(f)(x)&= C_d\int_{\mathbb{R}^n}\widehat{\mathfrak{a}}_g(x,\lambda)\,\chi(x)T^{\phi_1,\tau}_{\zeta}\bigl(f(\cdot - \lambda)\bigr)\, d\lambda\\&
=C_d\int_{\mathbb{R}^n}\widehat{\mathfrak{a}}_g(x,\lambda)\,T^{\phi_1,\tau}_{\chi \zeta}\bigl(f(\cdot - \lambda)\bigr)\, d\lambda
\end{align*}
H\"older inequality implies that
\begin{align}
\|T^{\phi_{1,\tau}}_{\mathfrak{a}_g}(f)(x)\|_{L^p(\mathbb{R}^n_x)}
&\lesssim \int_{\mathbb{R}^n}\|\widehat{\mathfrak{a}}_g(x,\lambda)\|_{L^{p_2}(\mathbb{R}^n_x)} \|T^{\varphi_1}_{\chi \zeta}(f(\cdot - \lambda))(x)\|_{L^{p_1}(\mathbb{R}^n_x)}\, d\lambda \nonumber\\
    &\lesssim \int_{\mathbb{R}^n} \|\widehat{\mathfrak{a}}_g(x,\lambda)\|_{L^{p_2}(\mathbb{R}^n_x)} \|f(\cdot - \lambda)\|_{L^{p_1}(\mathbb{R}^n_x)}\, d\lambda\nonumber\\&
\lesssim\int_{\mathbb{R}^n}\|\widehat{\mathfrak{a}}_g(x,\lambda)\|_{L^{p_2}(\mathbb{R}^n_x)}
\|f\|_{L^{p_1}(\mathbb{R}^n_x)}\, d\lambda
\end{align}

To estimate the norm 
$\|\widehat{\mathfrak{a}}_g(x,\lambda)\|_{L^{p_2}(\mathbb{R}^n_x)},$ observe that
\begin{align*}
\widehat{\mathfrak{a}}_g(x,\lambda)
&= \int_{\mathbb{R}^n} e^{ 2 \pi i\lambda \cdot \xi} \mathfrak{a}_g(x,\xi)\, d\xi \nonumber\\
&\lesssim  \langle \lambda \rangle^{-2N}
\int_{\mathbb{R}^n} e^{2\pi i\lambda \cdot \xi} (I - \Delta_\xi)^N \mathfrak{a}_g(x,\xi)\, d\xi.
\end{align*}
Then
\begin{align*}
\|\widehat{\mathfrak{a}}_g(x,\lambda)\|_{L^{p_2}(\mathbb{R}^n_x)}
\lesssim
\langle \lambda \rangle^{-2N}
\int_{\mathfrak{S}}
\|(I - \Delta_\xi)^N \mathfrak{a}_g(x,\xi)\|_{L^{p_2}(\mathbb{R}^n_x)}\, d\xi.
\end{align*}
where we have denoted the support of the symbol  $a_{12,\tau}$  in $\xi$  by $\mathfrak{S}.$  We can compute
\begin{align*}
(I - \Delta_\xi)^N \mathfrak{a}_g(x,\xi)
&= (I - \Delta_\xi)^N \bigl(\psi(\xi)\, T^{\phi_2,\tau}_{a_{12,\tau}(\cdot,\xi,\cdot)} g(x)\bigr) \nonumber\\
&= T^{\phi_{2,\tau}}_{\varkappa(\xi)} g(x).
\end{align*}
Note that the symbol of  $T^{\phi_{2,\tau}}_{\varkappa(\xi)}$  is given by
\begin{align*}
\varkappa(\xi) := (1 - \Delta_\xi)^N\bigl[\psi(\xi)\, a(\cdot,\xi,\cdot)\bigr],
\end{align*}
and is a linear combination of symbols of order  $\le m = m_1 + m_2 \le m_2,$  where we have used that $m_1 \le 0.$  Thus,  $\varkappa$  is a symbol of order  $m_2$  and we have the smoothing estimate
\begin{align*}
\|T^{\phi_{2,\tau}}_{\varkappa(\xi)} g(x)\|_{L^{p_2}}
\lesssim \|g\|_{L^{p_2}}.
\end{align*}
By taking  $N > \frac{n}{2},$  we have that
\begin{align*}
\|T^{\phi_{1,\tau}}_{{\mathfrak{a}}_g}(f)(x)\|_{L^p(\mathbb{R}^n_x)}
\lesssim
\int_{\mathbb{R}^n} \|\widehat{\mathfrak{a}}_g(x,\lambda)\|_{L^{p_2}(\mathbb{R}^n_x)}
\|f\|_{L^{p_1}(\mathbb{R}^n_x)}\, d\lambda.
\end{align*}
Summarising, we have proved that
\begin{align}
\|T^{\phi_{1,\tau}, \phi_{2,\tau}}_{a_{12,\tau}}(f,g)\|_{L^p(\mathbb{R}^n_x))}
&\lesssim \int_{\mathbb{R}^n} \langle \lambda \rangle^{-2N}
\int_{\mathfrak{S}} \|(I - \Delta_\xi)^N a_g(x,t,\xi)\|_{L^{p_2}(\mathbb{R}^n_x)}\, d\xi \,
\nonumber\\& \quad \times \|f\|_{L^{p_1}(\mathbb{R}^n_x)}\, d\lambda\nonumber \\
&= \int_{\mathbb{R}^n} \langle \lambda \rangle^{-2N}
\int_{\mathfrak{S}} \|T^{\phi_2,\tau}_{\varkappa(\xi)} g\|_{L^{p_2}(\mathbb{R}^n_x)}\, d\xi \,
\|f\|_{L^{p_1}(\mathbb{R}^n_x)}\, d\lambda \nonumber\\
&\lesssim \int_{\mathbb{R}^n} \langle \lambda \rangle^{-2N} d\lambda \,
\sup_{\xi \in \mathfrak{S}} \|T_{\varkappa(\xi)}\|_{\mathscr{B}(L^{p_2}, L^{p_2})}
|\mathfrak{S}| \,
\|g\|_{L^{p_2}(\mathbb{R}^n_x)}\,
\|f\|_{L^{p_1}(\mathbb{R}^n_x)} \nonumber\\
&\lesssim_{\mathfrak{S}}
\|f\|_{L^{p_1}(\mathbb{R}^n_x)}
\|g\|_{L^{p_2}(\mathbb{R}^n_x)}.
\end{align}
To deal with the term 
$T^{\phi_1,\tau,\phi_2,\tau}_{a_{12,g}}\big(\mu(D)f, \mu(D)g\big),$  let us introduce two smooth cut-off functions 
\[
\chi, \nu : \mathbb{R}^{2n} \to \mathbb{R},
\]
such that
\begin{align*}
\chi(\xi,\eta) =
\begin{cases}
1, & |(\xi,\eta)| \le 1, \\
0, & |(\xi,\eta)| \ge 2,
\end{cases}
\end{align*}
and
\begin{align*}
\nu(\xi,\eta) =
\begin{cases}
0, & \lambda^2 |\xi| \le 16 |\eta|, \\
1, & 64 |\eta| \le \lambda^2 |\xi|.
\end{cases}
\end{align*}
Defining
\[
a_{1,\tau}(x,\xi,\eta)
= \big(1 - \chi(\xi,\eta)\big)\,\nu(\xi,\eta)\,a_{12,\tau}(x,\xi,\eta)
\]
and
\[
a_{2,\tau}(x,\xi,\eta)
= \big(1 - \chi(\xi,\eta)\big)\,\big(1 - \nu(\xi,\eta)\big)\,a_{12,\tau}(x,\xi,\eta),
\]
we have that
\[
a_{1,\tau},\, a_{2,\tau} \in S^m_{1/2,1/2}(\mathbb{R}^n_x \times \mathbb{R}^n_{\xi}\times \mathbb{R}^n_{\eta}),
\] where $S^m_{1/2,1/2}(\mathbb{R}^n_x \times \mathbb{R}^n_{\xi}\times \mathbb{R}^n_{\eta})$ denotes the symbol class of order $m$ in the spatial variable $x \in \mathbb{R}^n_x$ with two frequency variables $(\xi,\eta) \in \mathbb{R}^n_{\xi} \times \mathbb{R}^n_{\eta},$
and we can decompose
\[T_{a}^{\phi_1, \phi_2}\big(\mu(D)f, \mu(D)g\big)
=T_{a_1}^{\phi_1, \phi_2}\big(\mu(D)f,\mu(D)g\big)+T_{a_2}^{\phi_1, \phi_2}\big(\mu(D)f, \mu(D)g\big).
\]
One introduces an even, real-valued smooth function $\psi$ whose Fourier transform is supported on the annulus
$\left\{ \xi \;\middle|\; \tfrac{1}{2} \le |\xi| \le 2 \right\}$
such that
\begin{equation}
\int_0^{\infty} \left| \widehat{\psi}(t\xi) \right|^2 \, \frac{dt}{t} = 1
\quad \text{for } \xi \neq 0. 
\end{equation}
Let $\theta$ be another real-valued smooth function whose Fourier transform is equal to one on the ball
$\left\{ \xi \;\middle|\; |\xi| \le \tfrac{1}{8} \right\}$
and supported in
$\left\{ \xi \;\middle|\; |\xi| \le \tfrac{1}{4} \right\}.$
Here and in the sequel, we denote
$\mathcal{K} = \operatorname{supp}_x a_{12,\tau}.$
\begin{align*}
&T_{a}^{\phi_1, \phi_2}\big(\mu(D)f, \mu(D)g\big)(x)\\
&= \int\limits_{0}^{\infty}\int_{\mathbb{R}^n}\int_{\mathbb{R}^n}
a_{1,\tau}(x, t\xi, t\eta) \\& \quad \times
\Big( \psi\big(t\nabla_x \phi_{1,\tau}(x,\xi)\big)
\left|\nabla_x \phi_{1,\tau}(x,\xi)\right|^{m}
\mu(\xi)\,\widehat{f}(\xi)\, e^{i\phi_{1,\tau}(x,\xi)}
\Big) \\
&\quad \times
\Big(
\theta\big(t\nabla_x \phi_{2,\tau}(x,\eta)\big)
\mu(\eta)\,\widehat{g}(\eta)\, e^{i\phi_{2,\tau}(x,\eta)}
\Big)
\,\frac{dt}{t}\, d\xi\, d\eta.
\end{align*}
For \[
a_{12,\tau,t}(x,\xi,\eta):=
t^{m}a_{12,\tau}\!\left(x,\frac{\xi}{t},\frac{\eta}{t}\right)\left|\nabla_x \phi_{1,\tau}(x,\xi)\right|^{-m}
\,\psi\!\left(\nabla_x \phi_{1,\tau}(x,\xi)\right)\hat{\theta}\left(\nabla_x \phi_{2,\tau}(x,\eta)\right).
\]
We define
\begin{equation}
\widetilde{a}_{12,\tau}(t,x,\xi,\eta)
:=a_{12,\tau,t}\bigl(x,\Psi_1(x,\xi),\Psi_2(x,\eta)\bigr),
\end{equation}
where
\[
\Psi_j(x,\xi)
=\bigl(\nabla_x \phi_{j,\tau}(x,\cdot)\bigr)^{-1}(\xi),
\qquad j=1,2,
\]
whose existence follows from the assumption of compact support in $x$, homogeneity of degree one in $\xi$, and the non-degeneracy.
\subsection{ Analysis on the  high frequencies}
In this subsection, we analyse the high-frequency portion of the bilinear Fourier integral operator with complex-valued phase functions by applying a Coifman--Meyer type decomposition.

\begin{align*}
&T^{\phi_{1,\tau},\phi_{2,\tau}}_{a_1}(\mu(D)f,\mu(D)g)(x)
\\&=\displaystyle \int\limits_{0}^1\iint
T^{\phi_{1,\tau}}_{\nu^{1}_{t,u}} f(x)
T^{\phi_{2,\tau}}_{\mu^{1}_{t,v}} g(x)
\frac{m_{\tau}(t,x,u,v)}{(1+4\pi^2|u|^2+4\pi^2|v|^2)^N}
\,du\,dv\,\frac{dt}{t}. 
\end{align*}
Coifman--Meyer decomposition is used to break each of the functions into high and low frequencies. It then combines the two functions as a para-product to separate the high frequencies of $f,$ to interact with high frequencies from $g,$ high frequencies from $f,$ to interact with low frequencies from $g,$ high frequencies from $g,$ to interact with low frequencies of $f,$ and the low frequencies from both functions $f$ and $g$ to interact. 

The strategy above helps us to analyse the operator with high and low frequencies separately, which is useful to obtain the boundedness of bilinear Fourier integral operators. Now, using the decomposition \`a la Coifman-Meyer (see \cite[Page~21]{rodriguez2011boundedness}) yields,
\begin{align*}
\int_{0}^{1} 
T^{\phi_{1,\tau}}_{v^{t,u}_1} f(x)\,
T^{\phi_{2,\tau}}_{\mu^{t,v}_1} g(x)\,
m_{\tau}(x)\,\frac{dt}{t}
= I_1 + I_2 + I_3,
\end{align*}
with
\begin{align*}
I_1 &= \int_{0}^{1} 
Q_t^u T^{\phi_{1,\tau}}_{\nu_{m_1}} f(x)\,
P_t^v T^{\phi_{2,\tau}}_{\mu_t^{\,1}} g(x)\,
m_{\tau}(x)\,\frac{dt}{t}, \\
I_2 &= \int_{0}^{1} 
Q_t^u T^{\phi_{1,\tau}}_{\nu_{m_1}} f(x)\,
E_t^2 g(x)\,
m_{\tau}(x)\,\frac{dt}{t},
\end{align*}
and
\begin{align*}
I_3 &= \int_{0}^{1} 
E_t^1 f(x)\,
t^{-m_2} T^{\phi_{2,\tau}}_{\mu^{t,v}_1} g(x)\,
m_{\tau}(x)\,\frac{dt}{t}.
\end{align*}
Note that  with $m<0$ and $m_2<0$ we can have this representation which is useful for estimating $I_1.$
To estimate $I_1$, we further decompose the function $g$ using the smoothing operators $S_{r,t}^{j}$, $j=1,2,3$. 

This second decomposition is necessary because $I_1$ represents the interaction of the high-frequency components of both $f$ and $g$, which is difficult to handle. The operators $S_{r,t}^{j}$ split the high-frequency contribution of $g$ into three simpler pieces, which continue to interact with the high-frequency component of $f$. Assuming these components form dyadic frequency pieces, we then sum over all dyadic frequency scales through a continuous frequency decomposition. 

The integration with respect to the parameter $r$ accounts for all frequency scales, that is, when $|\xi|\le |\eta|$, $|\eta|\le |\xi|$ and $|\eta|\approx|\xi|$  ensuring that all possible frequency interactions are considered.
\begin{align*}
I_1&= \sum_{j=1}^{3} \int_{\kappa}^{\infty} \int_{0}^{1}
\chi_{(0,1)}(rt)
\left(Q_t^u T^{\phi_1,\tau}_{\nu_{m_1}} f\right)(x)
\left(P_t^v S_{r,t}^j g\right)(x)
\, m_{\tau}(x)\,
\frac{dt\,dr}{tr},
\end{align*}
where, $S_{r,t}^1=r^{m_2}P_{1,t}Q_{tr}T^{\phi_2,\tau}_{\gamma},$
 with $T^{\phi_{2,\tau}}_{\gamma}$ being a Fourier integral operator of order $m_2$ with symbol
 \begin{align*}
\gamma(x,\eta)
= \chi(x)\mu(\eta)\left|\nabla_x \phi_{2,\tau}(x,\eta)\right|^{m_2}
\in S^{m_2}_{1,0}(\mathbb{R}^n_x  \times \mathbb{R}^n_\eta).
\end{align*}
The operator  $P_{1,t}$  is the Fourier multiplier with symbol  $\hat{\theta}_1(t\eta)$ and $Q_{tr}$ represents the Fourier multiplier  with symbol $\hat{\psi}(tr\eta)$.
The operators  $S_{r,t}^2$ and  $S_{r,t}^3$ have the form:

\begin{align*}
S_{r,t}^2 g(x)
&= r^{m_2} P_{1,t}
\left(
\int e^{i\phi_{2,\tau}(\cdot,\eta)} R^{2}_{tr}(\cdot,\eta)\,\widehat{g}(\eta)\, d\eta
\right)(x), 
\end{align*}
and
\begin{align*}
S_{r,t}^3 g(x)
&= r^{m_2}
\left(
\int e^{i\phi_{2,\tau}(x,\eta)} R^{1}_{t}(x,\eta)\,\widehat{g}(\eta)\, d\eta
\right),
\end{align*}
where
$ R^{1}_{t}(x,\eta)\text{ and }  \; R^{2}_{tr}(x,\eta)$  are the remainder symbols which belong to\\ 
 $S^{m_2 - \left(\frac{1}{2} - \varepsilon\right)}_{1,0}
\left(\mathbb{R}^n_x  \times \mathbb{R}^n_\eta\right)$.
Moreover, these symbols have their seminorms bounded by $t^{\varepsilon}$
and  $(tr)^{\varepsilon}$, respectively, with 
$\varepsilon \in (0, \tfrac{1}{2}).$
Finding the bound of (\ref{D}) is analogous to
\begin{align*}
\bigg|\int_{0}^{1} 
T^{\phi_{1,\tau}}_{v^{t,u}_1} f(x)\,
T^{\phi_{2,\tau}}_{\mu^{t,v}_1} g(x)\,
m_{\tau}(x)\,\frac{dt}{t}\bigg|
\le  |I_1| + |I_2| + |I_3|
\end{align*}
We aimed to find the bounds of each $I_i$ for $i=1,2,3$ separately. Let us consider $I_1$ decompose as follows:
\begin{align*}
I_1 &= \int_{0}^{1} 
Q_t^u T^{\phi_{1,\tau}}_{\nu_{m_1}} f(x)\,
P_t^v T^{\phi_{2,\tau}}_{\mu_t^{\,1}} g(x)\,
m_{\tau}(x)\,\frac{dt}{t}\\&
=\sum_{j=1}^{3}
\int_{\kappa}^{\infty}
\int_{0}^{1}
\chi_{(0,1)}(rt)\,
\bigl(Q^{u}_{t} T^{\phi_{1,\tau}}_{\nu_{m_1}} f\bigr)(x)\,
\bigl(P^{v}_{t} S^{j}_{r,t} g\bigr)(x)\,
m_{\tau}(x)\,
\frac{dt\,dr}{t r}\\&
= \int_{\kappa}^{\infty}
\int_{0}^{1}
\chi_{(0,1)}(rt)\,
\bigl(Q^{u}_{t} T^{\phi_{1,\tau}}_{\nu_{m_1}} f\bigr)(x)\,
\bigl(P^{v}_{t} S^{1}_{r,t} g\bigr)(x)\,
m_{\tau}(x)\,
\frac{dt\,dr}{t r}\\&+\int_{\kappa}^{\infty}
\int_{0}^{1}
\chi_{(0,1)}(rt)\,
\bigl(Q^{u}_{t} T^{\phi_{1,\tau}}_{\nu_{m_1}} f\bigr)(x)\,
\bigl(P^{v}_{t} S^{2}_{r,t} g\bigr)(x)\,
m_{\tau}(x)\,
\frac{dt\,dr}{t r} \\&+\int_{\kappa}^{\infty}
\int_{0}^{1}
\chi_{(0,1)}(rt)\,
\bigl(Q^{u}_{t} T^{\phi_{1,\tau}}_{\nu_{m_1}} f\bigr)(x)\,
\bigl(P^{v}_{t} S^{3}_{r,t} g\bigr)(x)\,
m_{\tau}(x)\,
\frac{dt\,dr}{t r}.
\end{align*}
Note that $I_1$ is the composition of a bilinear paraproduct with two Fourier integral operators. Of which, we split the above expression of $I_1$ into three terms with the first term containing the main contribution factor from $g,$ while the second and the third terms contain the minor contributions from $g$. 

The methodology described above helps to estimate each sub-term of $I_1$ separately, such that the first term can be estimated by using the duality argument and the remaining terms are estimated by using the H\"older's inequality. Here, we are using the fact that the following $L^{p_2}(\mathbb{R}^n_x)$-estimates are satisfied:
\begin{itemize}
    \item $\|P^{v}_{t} S^{3}_{r,t} g(x)\|_{L^{p_2}(\mathbb{R}^n_x)} \lesssim \|S^{3}_{r,t} g(x)\|_{L^{p_2}(\mathbb{R}^n_x)}$,
    
    \item $\|P^{v}_{t} S^{2}_{r,t} g(x)\|_{L^{p_2}(\mathbb{R}^n_x)} \lesssim \|S^{2}_{r,t} g(x)\|_{L^{p_2}(\mathbb{R}^n_x)}$,
    
    \item $\|Q^{u}_{t} T^{\phi_1}_{\nu_{m_1}} f(x)\|_{L^{p_1}(\mathbb{R}^n_x)} \lesssim\|T^{\phi_1}_{\nu_{m_1}} f(x)\|_{L^{p_1}(\mathbb{R}^n_x)}$.
\end{itemize}
The above properties are necessary to indicate that $P^v_t$ and $Q^u_t$ are bounded on $L^{p_2}(\mathbb{R}^n_x)$ and $L^{p_1}(\mathbb{R}^n_x)$ respectively (see \cite[Page~16]{Cardona:2026}): 
\begin{itemize}
    \item $\|S^{3}_{r,t} g(x)\|_{L^{p_2}(\mathbb{R}^n_x)} \lesssim r^{m_2} t^{\varepsilon} \|g\|_{L^{p_2}(\mathbb{R}^n_x)}$,
    
    \item $\|S^{2}_{r,t} g(x)\|_{L^{p_2}(\mathbb{R}^n_x )} \lesssim r^{m_2+\varepsilon} t^{\varepsilon} \|g\|_{L^{p_2}(\mathbb{R}^n_x)}$,
    
    \item $\|T^{\phi_1}_{\nu_{m_1}} f(x)\|_{L^{p_1}(\mathbb{R}^n_x)} \lesssim \|f\|_{L^{p_1}(\mathbb{R}^n_x)}$.
\end{itemize}
The above inequalities indicate that $S^{2}_{r,t}$ and $S^{3}_{r,t}$ are controlled by $r^{m_2} t^{\varepsilon}$ and $r^{m_2+\varepsilon} t^{\varepsilon}$ respectively and $T^{\phi_1}_{\nu_{m_1}}$ is bounded by $L^{p_1}(\mathbb{R}^n_x)$ (see \cite[Page~16]{Cardona:2026}).
To estimate the bound of $I_1$, we need to know the estimate of the following expressions. 

Notice also that, below 
\begin{align}\label{inf}
\|m_{\tau}\|_{L^\infty} := \|m_{\tau}\|_{L^\infty_U},
\end{align}
denotes the $L^\infty$-norm of $m_{\tau}$ as a function of $x$, and it is bounded polynomially as a function of $U = 1 + 4\pi^2(|u| + |v|)$.

\[
\begin{aligned}
&\left\|
\int_{\kappa}^{\infty}
\int_{0}^{1}
\chi_{(0,1)}(rt)\,
\bigl(Q^{u}_{t} T^{\phi_{1,\tau}}_{\nu_{m_1}} f\bigr)(x)\,
\bigl(P^{v}_{t} S^{3}_{r,t} g\bigr)(x)\,
m_{\tau}(x)\,
\frac{dt\,dr}{t r}
\right\|_{L^{p}(\mathbb{R}_x^n)} \\
&\le
\int_{\kappa}^{\infty}
\int\limits_{0}^{\frac{1}{r}}
\left\|
\bigl(Q^{u}_{t} T^{\phi_{1,\tau}}_{\nu_{m_1}} f\bigr)(x)\,
\bigl(P^{v}_{t} S^{3}_{r,t} g\bigr)(x)\,
m(x)
\right\|_{L^{p}(\mathbb{R}_x^n)}
\frac{dt\,dr}{t} \\
&\lesssim
\int_{\kappa}^{\infty}
\int\limits_{0}^{\frac{1}{r}}
\|
\bigl(Q^{u}_{t} T^{\phi_{1,\tau}}_{\nu_{m_1}} f\bigr)(x)
\|_{L^{p_1}(\mathbb{R}_x^n)}
\|
\bigl(P^{v}_{t} S^{3}_{r,t} g\bigr)(x)
\|_{L^{p_2}(\mathbb{R}_x^n)}
\|m\|_{L^\infty}
\frac{dt\,dr}{t} \\
&\lesssim
\int_{\kappa}^{\infty}
\int\limits_{0}^{\frac{1}{r}}
r^{m_2} t^{\varepsilon}
\frac{dt\,dr}{t r}\,
\|f\|_{L^{p_1}(\mathbb{R}^n_x)}\,
\|g\|_{L^{p_2}(\mathbb{R}^n_x)}.
\end{aligned}
\]
We have used the fact that,
\[
\int_{\kappa}^{\infty} \int_{0}^{1/r} r^{m_2} t^{\varepsilon} \,\frac{dt\,dr}{t r} < \infty,
\] since $m_2 < 0$. 

\begin{remark}
We have
\[
\int_{\kappa}^{\infty} \int_{0}^{1/r} r^{m_2} t^{\varepsilon} \,\frac{dt\,dr}{t r} < \infty,
\]
in particular for $m_2 = 0$, provided that $m_2 < \varepsilon$, where $\varepsilon > 0$, is small enough.
\end{remark}
Again, we have the estimate below as
\begin{align*}
&\left\|
\int_{\kappa}^{\infty}
\int_{0}^{1}
\chi_{(0,1)}(rt)\,
\bigl(Q^{u}_{t} T^{\phi_{1,\tau}}_{\nu_{m_1}} f\bigr)(x)\,
\bigl(P^{v}_{t} S^{2}_{r,t} g\bigr)(x)\,
m_{\tau}(x)\,
\frac{dt\,dr}{t r}
\right\|_{L^{p}(\mathbb{R}_x^n)} \\
&\le
\int_{\kappa}^{\infty}
\int\limits_{0}^{\frac{1}{r}}
\left\|
\bigl(Q^{u}_{t} T^{\phi_{1,\tau}}_{\nu_{m_1}} f\bigr)(x)\,
\bigl(P^{v}_{t} S^{2}_{r,t} g\bigr)(x)\,
m_{\tau}(x)
\right\|_{L^{p}(\mathbb{R}_x^n)}
\frac{dt\,dr}{t} \\
&\lesssim
\int_{\kappa}^{\infty}
\int\limits_{0}^{\frac{1}{r}}
\|
\bigl(Q^{u}_{t} T^{\phi_{1,\tau}}_{\nu_{m_1}} f\bigr)(x)
\|_{L^{p_1}(\mathbb{R}_x^n)}
\|
\bigl(P^{v}_{t} S^{2}_{r,t} g\bigr)(x)
\|_{L^{p_2}(\mathbb{R}_x^n)}
\|m\|_{L^\infty}
\frac{dt\,dr}{t} \\
&\lesssim
\int_{\kappa}^{\infty}
\int\limits_{0}^{\frac{1}{r}}
r^{m_2+\varepsilon} t^{\varepsilon}
\frac{dt\,dr}{t r}\,
\|f\|_{L^{p_1}(\mathbb{R}^n_x)}\,
\|g\|_{L^{p_2}(\mathbb{R}^n_x)}.
\end{align*}

Again, we use the fact that, $$\int_{\kappa}^{\infty}
\int\limits_{0}^{\frac{1}{r}}
r^{m_2+\varepsilon} t^{\varepsilon}
\frac{dt\,dr}{t r}< \infty$$ since $m_2<0.$ Indeed, note that:
\begin{align*}
\int_{\kappa}^{\infty}\int_{0}^{1/r}r^{m_{2}+\varepsilon}t^{\varepsilon}\frac{dt\,dr}{tr}&= \int_{\kappa}^{\infty}\int_{0}^{1/r}r^{m_{2}+\varepsilon-1}t^{\varepsilon-1} dt\,dr\\
&=\int_{\kappa}^{\infty}
r^{m_{2}+\varepsilon-1}
\left(\int_{0}^{1/r}
t^{\varepsilon-1}\,dt\right)dr\\&
= \frac1{\varepsilon}\int_{\kappa}^{\infty}
r^{m_{2}+\varepsilon-1}\bigg[t^{\varepsilon}\bigg]_{0}^{1/r}\,dr
\\&=\frac1{\varepsilon}
\int_{\kappa}^{\infty}
r^{m_{2}+\varepsilon-1}
\left(\frac1r\right)^{\varepsilon}
dr \\&=\frac1{\varepsilon}
\int_{\kappa}^{\infty} r^{m_{2}-1}\,dr.
\end{align*}
Since $m_2<0,$ implies that $m_2-1<-1.$ Moreover, $\varepsilon>0.$ One has that
\begin{align*}
\frac{1}{\varepsilon}\int_{\kappa}^{\infty}
r^{m_{2}-1}\,dr=\frac{1}{\varepsilon}\left[
\frac{r^{m_{2}}}{m_{2}}
\right]_{\kappa}^{\infty}
=-\frac{\kappa^{m_{2}}}{ \varepsilon m_{2}}<\infty.
\end{align*}
This completes the above estimate.

Let us consider this last estimation in order to complete the bound for $I_1.$ In this case, 
\[
S^{1}_{r,t} = r^{m_2} P_{1,t} Q_{tr} T^{\phi_{2,\tau}}_{\gamma},
\]

\[
\left\|
\int_{\kappa}^{\infty}
r^{m_2}
\int_{0}^{1}
\chi_{(0,1)}(rt)\,
\bigl(Q^{u}_{t} T^{\phi_{1,\tau}}_{\nu_{m_1}} f\bigr)(x)\,
\bigl(P^{v}_{t} [P_{1,t} Q_{tr} T^{\phi_{2,\tau}}_{\gamma}] g\bigr)(x)\,
m_{\tau}(x)\,
\frac{dt\,dr}{t r}
\right\|_{L^{p}(\mathbb{R}_x^n)}.
\]
We use the duality argument in this estimation, so we take an arbitrary $h \in L^{p'}(\mathbb{R}_x^n )$, assumed to be supported in $x$.
Now, by duality, we have
\[
\|F\|_{L^p(\mathbb{R}_x^n)}
= \sup_{\|h\|_{L^{p'}(\mathbb{R}_x^n)} \le 1}
\left|
\int_{\mathbb{R}^n} F(x)\, h(x)\, dx
\right| \text{ with } \frac{1}{p}+\frac{1}{p^{'}}=1.
\]
Hence, it suffices to estimate
\[
\sup_{\|h\|_{L^{p'}(\mathbb{R}_x^n)} \le 1}
\left|
\int_{\mathbb{R}^n}
\int_{0}^{1}
h(x)\,
\bigl(Q^{u}_{t} T^{\phi_{1,\tau}}_{\nu_{m_1}} f\bigr)(x)\,
\bigl(\tilde{P}^{v}_{t} [R_t T^{\phi_{2,\tau}}_{\gamma}] g\bigr)(x)\,
m_{\tau}(x)\,
\frac{dt}{t}\,
dx
\right|.
\]
This quantity is bounded by
\[
\lesssim
\|h\|_{L^{p'}(\mathbb{R}_x^n)}\,
\|f\|_{L^{p_1}(\mathbb{R}^n_x)}\,
\|g\|_{L^{p_2}(\mathbb{R}^n_x)}.
\]
Here, we use the fact that, $\|h\|_{L^{p'}(\mathbb{R}_x^n)}\le 1$. Let us define: \\ $$R_t(G)(x) =  \int\limits_{\kappa}^{\frac{1}{t}} r^{m_2} Q_{t r}(G)(x)\,\frac{dr}{r}, 
\quad \text{and ,} \quad 
\tilde{P^v_t} := P_{t}^v P_{1,t}.$$
A further decomposition was applied to $\tilde{P^v_t}$ which results $\tilde{P^v_t}= P^v_t+ Q^v_t$. Now, let us  define $\mathscr{S}$ as follows:
\begin{align*}
    \mathscr{S}:= \int_{\mathbb{R}^n} \int_0^1 
h(x)\,(Q^u_t T^{\phi_{1,\tau}}_{\nu_{m_1}} f)(x)\,
\big( \widetilde{P}^v_t [R_t T^{\phi_{2,\tau}}_\gamma] g \big)(x)\,
m_{\tau}( x)\,\frac{dt}{t}\,dx.
\end{align*}
Finding the bound of $\mathscr{S}$ gives us,
\begin{align*}
|\mathscr{S}| 
&= \left|
\int_{\mathbb{R}^n} \int_0^1 
h(x)\,(Q^u_t T^{\phi_{1,\tau}}_{\nu_{m_1}} f)(x)\,
(\widetilde{P}^v_t [R_t T^{\phi_{2,\tau}}_\gamma] g)(x)\,
m_{\tau}( x)\,\frac{dt}{t}\,dx
\right| \\
&\leq \left|
\int_{\mathbb{R}^n} \int_0^1 
h(x)\,(Q^u_t T^{\phi_{1,\tau}}_{\nu_{m_1}} f)(x)\,
(P^{v,2}_t [R_t T^{\phi_{2,\tau}}_\gamma] g)(x)\,
m_{\tau}( x)\,\frac{dt}{t}\,dx
\right| \\
&\quad + \left|
\int_{\mathbb{R}^n} \int_0^1 
h(x)\,(Q^u_t T^{\phi_{1,\tau}}_{\nu_{m_1}} f)(x)\,
(Q^{v,2}_t [R_t T^{\phi_{2,\tau}}_\gamma] g)(x)\,
m_{\tau}( x)\,\frac{dt}{t}\,dx
\right| \\
&= \left|
\int_{\mathbb{R}^n} \int_0^1 
\widetilde{Q}_t\big((Q^u_t T^{\phi_{1,\tau}}_{\nu_{m_1}} f)(x)\big)\,
(P^{v,2}_t [R_t T^{\phi_2,\tau}_\gamma] g)(x)\,
h(x)\,m_{\tau}( x)\,\frac{dt}{t}\,dx
\right| \\
&\quad + \left|
\int_{\mathbb{R}^n} \int_0^1 
\widetilde{P}_t\big((Q^u_t T^{\phi_{1,\tau}}_{\nu_{m_1}} f)(x)\big)\,
(Q^{v,2}_t [R_t T^{\phi_{2,\tau}}_\gamma] g)(x)\,
h(x)\,m_{\tau}( x)\,\frac{dt}{t}\,dx
\right| \\
&= \left|
\int_{\mathbb{R}^n} \int_0^1 
(Q^u_t T^{\phi_{1,\tau}}_{\nu_{m_1}} f)(x)\,
(P^{v,2}_t [R_t T^{\phi_{2,\tau}}_\gamma] g)(x)\,
\widetilde{Q}_t\big(h(\cdot)m_{\tau}(\cdot)\big)(x)\,
\frac{dt}{t}\,dx
\right| \\
&\quad + \left|
\int_{\mathbb{R}^n} \int_0^1 
(Q^u_t T^{\phi_{1,\tau}}_{\nu_{m_1}} f)(x)\,
(Q^{v,2}_t [R_t T^{\phi_{2,\tau}}_\gamma] g)(x)\,
\widetilde{P}_t\big(h(\cdot)m_{\tau}(\cdot)\big)(x)\,
\frac{dt}{t}\,dx
\right|.
\end{align*}
Let us define
\begin{equation*}
I = \int_{\mathbb{R}^n} \int_0^1 
(Q^u_t T^{\phi_{1,\tau}}_{\nu_{m_1}} f)(x)\,
(P^{v,2}_t [R_t T^{\phi_{2,\tau}}_\gamma] g)(x)\,
\widetilde{Q}_t\big(h(\cdot)m_{\tau}( \cdot)\big)(x)\,
\frac{dt}{t}\,dx,
\end{equation*}
and
\begin{equation*}
II = \int_{\mathbb{R}^n} \int_0^1 
(Q^u_t T^{\phi_{1,\tau}}_{\nu_{m_1}} f)(x)\,
(Q^{v,2}_t [R_t T^{\phi_{2,\tau}}_\gamma] g)(x)\,
\widetilde{P}_t\big(h(\cdot)m_{\tau}(\cdot)\big)(x)\,
\frac{dt}{t}\,dx.
\end{equation*}
Using the relation 
\begin{align}\label{Q}
\tilde{Q_t}((h(\cdot)m_{\tau}(\cdot))(x)= m_{\tau}(\cdot)\tilde{Q_t}(h(\cdot))+ [ \tilde{Q_t}, M_{m_{\tau}}]h(\cdot),
\end{align}
and 
\begin{align}\label{P}
\tilde{P_t}((h(\cdot)m_{\tau}(\cdot))(x)= m_{\tau}(\cdot)\tilde{P_t}(h(\cdot))+ [ \tilde{P_t}, M_{m_{\tau}}]h(\cdot). \end{align}
Using \eqref{Q} and \eqref{P}, we express $I$ and $II$ as follows:
\begin{align*}
I 
&= \int_{\mathbb{R}^n} \int_0^1 
(Q^u_t T^{\phi_{1,\tau}}_{\nu_{m_1}} f)(x)\,
(P^{v,2}_t [R_t T^{\phi_{2,\tau}}_\gamma] g)(x)\,
m_{\tau}( x)\,(\widetilde{Q}_t h(\cdot))(x)\,
\frac{dt}{t}\,dx \nonumber \\
&\quad + \int_{\mathbb{R}^n} \int_0^1 
(Q^u_t T^{\phi_{1,\tau}}_{\nu_{m_1}} f)(x)\,
(P^{v,2}_t [R_t T^{\phi_{2,\tau}}_\gamma] g)(x)\,
\big([\widetilde{Q}_t, M_{m_{\tau}}]h(\cdot)\big)(x)\,
\frac{dt}{t}\,dx  \nonumber\\
&= I_4 + I_5. 
\end{align*}
Similarly,
\begin{align*}
II 
&= \int_{\mathbb{R}^n} \int_0^1 
(Q^u_t T^{\phi_{1,\tau}}_{\nu_{m_1}} f)(x)\,
(Q^{v,2}_t [R_t T^{\phi_{2,\tau}}_\gamma] g)(x)\,
m_{\tau}( x)\,(\widetilde{P}_t h(\cdot))(x)\,
\frac{dt}{t}\,dx \nonumber\\
&\quad + \int_{\mathbb{R}^n} \int_0^1 
(Q^u_t T^{\phi_{1,\tau}}_{\nu_{m_1}} f)(x)\,
(Q^{v,2}_t [R_t T^{\phi_{2,\tau}}_\gamma] g)(x)\,
\big([\widetilde{P}_t, M_{m_{\tau}}]h(\cdot)\big)(x)\,
\frac{dt}{t}\,dx \nonumber\\
&= I_6 + I_7.
\end{align*}
Now, $I_5$ and $I_7$ can be estimated separately since they have commutator terms, using $L^{p'}\text{-boundedness}$, and $I_4$ and $I_6$ can also be estimated using the Cauchy-Schwarz inequality. By considering the boundedness properties below,
\begin{itemize}
    \item 
    \[
    \|[\widetilde{P}_t, M_{m_{\tau}}]h(\cdot)\|_{L^{p'}(\mathbb{R}^n_x)} 
    \lesssim t \|h(\cdot)\|_{L^{p'}(\mathbb{R}^n_x)}.
    \]
\item 
    \[
    \|[\widetilde{Q}_t, M_{m_{\tau}}]h(\cdot)\|_{L^{p'}(\mathbb{R}^n_x)} 
    \lesssim t \|h(\cdot)\|_{L^{p'}(\mathbb{R}^n_x)}.
    \]
\item 
    \[
    \|P^{v,2}_t [R_t T^{\phi_{2,\tau}}_\gamma] g\|_{L^{p_2}(\mathbb{R}_x^n)} 
    \lesssim C(v)\|T^{\phi_{2,\tau}}_\gamma g\|_{L^{p_2}(\mathbb{R}_x^n)} 
    \lesssim C(v)\|g\|_{L^{p_2}(\mathbb{R}^n_x)}.
    \]
 \item 
    \[
    \|Q^{v,2}_t [R_t T^{\phi_{2,\tau}}_\gamma] g\|_{L^{p_2}(\mathbb{R}_x^n)} 
    \lesssim C(v)\|T^{\phi_{2,\tau}}_\gamma g\|_{L^{p_2}(\mathbb{R}_x^n)} 
    \lesssim C(v)\|g\|_{L^{p_2}(\mathbb{R}^n_x)}.
    \]
\end{itemize}
we obtain
\[
|I_5| + |I_7| \lesssim C(u, v, \varepsilon)\,
\|f\|_{L^{p_1}(\mathbb{R}^n_x)}\,
\|g\|_{L^{p_2}(\mathbb{R}^n_x)}\,
\|h\|_{L^{p'}(\mathbb{R}^n_x)}.
\]
In order to finish our analysis of  $\mathscr{S}$,  let us estimate the terms  $I_4$ and  $I_6$.  By using the Cauchy--Schwarz inequality, we have that
\begin{align*}
|I_4| 
&= \left|
\int_{\mathbb{R}^n} 
\int_0^1 
(Q^u_{t} T^{\phi_{1,\tau}}_{ \nu^{m_1}} f)(x) 
\big(P^{v,2}_t [R_t T^{\phi_{2,\tau}}_{\gamma}] g\big)(x) 
\, m_{\tau}( x) \, \widetilde{Q}_t h(\cdot))(x) \, \frac{dt}{t} \, dx \right| \\
&= \left|\int_{\mathbb{R}^n} \int_0^1 
(Q^u_{t} T^{\phi_{1,\tau}}_{\nu^{m_1}} f)(x) 
\big(P^{v,2}_t [R_t T^{\phi_{2,\tau}}_{\gamma}] g\big)(x) 
\, m_{\tau}( x) \, (\widetilde{Q}_t h(\cdot))(x) \, \frac{dt}{t} \, dx \right|\\
&\le \left(\int_{\mathbb{R}^n} \int_0^1 
\big| (Q^u_{t} T^{\phi_{1,\tau}}_{ \nu^{m_1}} f)(x) \big|^2 
\, \frac{dt}{t} \, dx 
\right)^{\frac12} \\& \quad
\times
\left(\int_{\mathbb{R}^n} \int_0^1 
\big| (P^{v,2}_t [R_t T^{\phi_{2,\tau}}_{\gamma}] g)(x) \, m_{\tau}(x) \, (\widetilde{Q}_t h(\cdot))(x) \big|^2 
\, \frac{dt}{t} \, dx 
\right)^{\frac12}.
\end{align*}
For the first term, we obtain
\begin{align*}
&\left( \int_{\mathbb{R}^n} \int_0^1 \big| (Q^u_{t} T^{\phi_{1,\tau}}_{ \nu_{m_1}} f)(x) \big|^2 \, \frac{dt}{t} \, dx  \right)^{\frac12} \\
&\quad \le 
\int_{\mathbb{R}^n} \left( \int_0^1 \big| (Q^u_{t} T^{\phi_{1,\tau}}_{ \nu_{m_1}} f)(x) \big|^2 \, \frac{dt}{t} \right)^{\frac12} dx   \\
&\quad =    \left\| \left( \int_0^1 \big| (Q^u_{t} T^{\phi_{1,\tau}}_{ \nu_{m_1}} f)(\cdot) \big|^2 \, \frac{dt}{t} \right)^{\frac12} \right\|_{L^{2}(\mathbb{R}_x^n)} \\
&\quad \lesssim \| T^{\phi_{1,\tau}}_{\nu_{m_1}} f(\cdot) \|_{L^{p_1}(\mathbb{R}_x^n)}. 
\end{align*}
For the second factor, observe that
\begin{align*}
&\left( \int_{\mathbb{R}^n} \int_0^1 
\big| (P^{v,2}_{t} [R_t T^{\phi_{2,\tau}}_{\gamma}] g)(x)\, m_{\tau}(x)\, (\widetilde{Q}_t  h(\cdot))(x) \big|^2 
\, \frac{dt}{t} \, dx  \right)^{\frac12} \\
&\quad \le \|m_{\tau}\|_{L^\infty}
\left( \int_{\mathbb{R}^n}\int_0^1 
\big| (P^{v,2}_{t} [R_t T^{\phi_{2,\tau}}_{\gamma}] g)(x)\, (\widetilde{Q}_t  h(\cdot))(x) \big|^2 
\, \frac{dt}{t} \, dx  \right)^{\frac12}.
\end{align*}
Let
\[
A_t := (P^{v,2}_{t} [R_t T^{\phi_{2,\tau}}_{\gamma}] g)(x), 
\qquad 
B_t := (\widetilde{Q}_t  h(\cdot))(x).
\]
Then, by  H\"older’s inequality,
\begin{align*}
&\left( \int_{\mathbb{R}^n} \int_0^1 |A_t B_t|^2 \, \frac{dt}{t} \, dx \right)^{\frac12} \\
&\quad \le \int_{\mathbb{R}^n}  
\left( \int_0^1 |A_t B_t|^2 \, \frac{dt}{t} \right)^{\frac12} dx  \\
&\quad \le \left\|\sup_{0<t<1} |A_t|\right\|_{{L}^{p_2}(\mathbb{R}^n_x)} 
\left\|\left( \int_0^1 |B_t|^2 \, \frac{dt}{t} \right)^{\frac12} dx\right\|_{{L}^{p^{\prime}}(\mathbb{R}^n_x)}.  
\end{align*}
Using a similar argument as in \cite[Remark~4.4]{Cardona:2026}, by applying Bourgain's maximal function theorem in \cite{Bourgain1985}, allows us to conclude that 
\begin{align*}
\left\|\sup_{0<t<1}
\left|P_t^{v,2}\bigl[R_tT_{\gamma}^{\phi_{2,\tau}}g\bigr]\right|\right\|_{L^{p_2}(\mathbb{R}^n_x)}
\lesssim\|g\|_{L^{p_2}(\mathbb{R}^n_x)},
\end{align*}
and that one has also the inequality
\begin{align*}
\left\|\left(
\int_0^1|\widetilde{Q}_th(x)|^2
\frac{dt}{t}\right)^{1/2}
\right\|_{L^{p'}(\mathbb{R}^n_x)}
\lesssim\|h\|_{L^{p'}(\mathbb{R}^n_x)}.
\end{align*}
Combining all estimates above, we conclude that
\begin{align*}
|I_4|
\lesssim 
\| T^{\phi_{1,\tau}}_{\nu_{m_1}} f \|_{L^{p_1}(\mathbb{R}^n_x)}
\left\| \sup_{0<t<1} |P^{v,2}_{t} [R_t T^{\phi_{2,\tau}}_{\gamma}] g| \right\|_{L^{p_2}(\mathbb{R}^n_x)}
\left\|
\left( \int_0^1 | \widetilde{Q}_t h |^2 \, \frac{dt}{t} \right)^{\frac12}
\right\|_{L^{p'}(\mathbb{R}^n_x)}.
\end{align*}
The above inequality holds  due to the assumption that $1/p_1+1/p_2+1/p'=1$ (see \cite[Page~20]{Cardona:2026}).
In consequence,
$|I_4| \lesssim \|h\|_{L^{p'}(\mathbb{R}^n_x)} \, \|f\|_{L^{p_1}(\mathbb{R}^n_x)} \, \|g\|_{L^{p_2}(\mathbb{R}^n_x)}.$
Since \(h\) is assumed to be supported in \(\mathbb{R}^n_x\), and since the same analysis carried out for \(I_4\) also applies to \(I_6\), we obtain
\[
|I_4| + |I_6| \lesssim \|h\|_{L^{p'}(\mathbb{R}^n_x}) \, \|f\|_{L^{p_1}(\mathbb{R}^n_x)} \, \|g\|_{L^{p_2}(\mathbb{R}^n_x)}.
\]
Hence, we have completed the proof of the estimate
\[
|I| \lesssim \|h\|_{L^{p'}(\mathbb{R}^n_x)} \, \|f\|_{L^{p_1}(\mathbb{R}^n_x)} \, \|g\|_{L^{p_2}(\mathbb{R}^n_x)}.
\]
To estimate $I_2$, we obtain
\[
\| Q^u_{t} T^{\phi_1}_{\nu_{m_1}} f(x) \|_{L^{p_1}(\mathbb{R}^n_x)} 
\lesssim_u \|f\|_{L^{p_1}(\mathbb{R}^n_x)}.
\]
Then Hölder's inequality gives
\begin{align*}
\|I_2\|_{L^p(\mathbb{R}^n_x)} 
&\le \int_0^1 
\| Q^u_{t} T^{\phi_{1,\tau}}_{\nu_{m_1}} f(x)\, E_t^2 g(x)\, m_{\tau}(x) \|_{L^p(\mathbb{R}^n_x)} 
\, \frac{dt}{t} \\
&\le \int_0^1 
\| Q^u_{t} T^{\phi_{1,\tau}}_{\nu_{m_1}} f(x) \|_{L^{p_1}(\mathbb{R}^n_x)} 
\| E_t^2 g(x)\, m_{\tau}(x) \|_{L^{p_2}(\mathbb{R}^n_x)} 
\, \frac{dt}{t} \\
&\le \int_0^1 
\| Q^u_{t} T^{\phi_{1,\tau}}_{\nu_{m_1}} f(x) \|_{L^{p_1}(\mathbb{R}^n_x)} 
\| E_t^2 g(x) \|_{L^{p_2}(\mathbb{R}^n_x)} 
\|m_{\tau}\|_{L^\infty} 
\, \frac{dt}{t} \\
&\le \int_0^1 
\| E_t^2 g(x) \|_{L^{p_2}(\mathbb{R}^n_x)} 
\|m_{\tau}\|_{L^\infty} 
\, \frac{dt}{t} \, \|f\|_{L^{p_1}(\mathbb{R}^n_x)} \\
&\lesssim_u 
\int_0^1 t^{\varepsilon - 1} \|m_{\tau}\|_{L^\infty} \, dt \,\|f\|_{L^{p_1}(\mathbb{R}^n_x)} 
\|g\|_{L^{p_2}(\mathbb{R}^n_x)}.
\end{align*}
The same analysis above applies to \(I_3\). Recall that
$ \|m_{\tau}\|_{L^\infty} := \|m_{\tau}\|_{L^\infty_U}$
denotes the \(L^\infty\)-norm of \(m_{\tau}\) as a function of \((x)\), and it is
bounded polynomially as a function of \(U = 1 + 4\pi^2|u| + 4\pi^2|v|\). 
In consequence, taking \(N\) large enough, we have 
\begin{align*}
&\big\| T^{\phi_{1,\tau},\phi_{2,\tau}}_{a_{12,\tau}}(\mu(D)f, \mu(D)g)(x) \big\|_{L^p(\mathbb{R}^n_x)}
\\&\le \iint_{\mathbb{R}^{2n}} 
\frac{\|m_{\tau}\|_{L^\infty_U}}{(1 + 4\pi^2|u|^2 + 4\pi^2|v|^2)^N} \, du \, dv \times \|f\|_{L^{p_1}(\mathbb{R}^n_x)} 
\|g\|_{L^{p_2}(\mathbb{R}^n_x)} \\
&\lesssim_N 
\|f\|_{L^{p_1}(\mathbb{R}^n_x)} 
\|g\|_{L^{p_2}(\mathbb{R}^n_x)}.
\end{align*}

\subsection{ Main theorem}
We now formulate the main boundedness theorem for bilinear Fourier integral operators with complex-valued phase functions. The proof of this statement is the combination of the previous subsections. The following theorem formulates in detail the assumptions in Theorem \ref{Main:1:th}.

\begin{theorem}[Bilinear Seeger--Sogge--Stein theorem for complex phases]\label{Th:FIO:C:2}
Let \begin{align*}
T^{\phi_{1}, \phi_{2}}_{a}(f,g)(x)
&= \int_{\mathbb{R}^n}\int_{\mathbb{R}^n}
e^{2\pi i (\phi_{1}(x,\xi)+\phi_{2}(x,\eta))}
\, a(x,\xi,\eta)\, \widehat{f}(\xi)\, \widehat{g}(\eta)
\, d\xi \, d\eta,
\end{align*}
be a bilinear Fourier integral operator, where $a(x,\xi,\eta)\in S^m_{1,0}(\mathbb{R}^n_x \times \mathbb{R}^n_{\xi}\times \mathbb{R}^n_{\eta})$ is compactly supported in the spatial variable $x.$ Assume that there exists $\tau\in \mathbb{R},$ such that $\phi_{1,\tau}(x,\xi) = \Re \phi_1(x,\xi) + \tau \Im \phi_1(x,\xi) $ and $\phi_{2,\tau}(x,\eta) = \Re \phi_2(x,\eta) + \tau \Im \phi_2(x,\eta)$ satisfy:
\begin{itemize}
        \item[1.] $\phi_{1,\tau}(x,\xi)$ is positively homogeneous of degree one in $\xi,$ and $\phi_{2,\tau}(x,\eta)$ is also positively homogeneous of degree one in $\eta,$
        \item[2.]  $\phi_{1,\tau}(x,\xi) $ and $\phi_{2,\tau}(x,\eta)$ satisfy the non-degeneracy condition.
    \end{itemize} 
    Then whenever $$1< p_1,p_2< \infty,\quad 1<p < \infty, \quad\frac{1}{p}=\frac{1}{p_1}+\frac{1}{p_2}$$ and
\begin{align*}
m\le-(n-1)\left(
\left|\frac{1}{p_1}-\frac{1}{2}\right|
+\left|\frac{1}{p_2}-\frac{1}{2}\right|
\right), \quad n\ge 2,
\end{align*} the bilinear Fourier integral operator $T^{\phi_{1}, \phi_{2}}_{a}
:L^{p_1}(\mathbb{R}^n)\times L^{p_2}(\mathbb{R}^n)
\longrightarrow L^p(\mathbb{R}^n)$ extends to a bounded bilinear operator. Moreover, $\|T^{\phi_{1},\phi_{2}}_{a}(f,g)\|_{L^p}
\lesssim \|f\|_{L^{p_1}}\|g\|_{L^{p_2}}.$
\end{theorem}
\section{Endpoint estimates for $m=-(n-1),$ $n\geq 2$}\label{endpoint:estimates}
In this section, we discuss the endpoint estimates of the bilinear Fourier integral operators with complex-valued phases, and provide a counterexample in spaces where boundedness fails. 
Now, we seek to find the boundedness of bilinear Fourier integral operator $T_a^{\phi_1,\phi_2},$ for high frequency support of $T_a^{\phi_1,\phi_2}$, in this case, we have $|\eta|\lesssim |\xi|$ and since $T_{a_{1}}^{\phi_1,\phi_2}$ and $T_{a_{2}}^{\phi_1,\phi_2}$ are symmetric, we analyse just $T_{a_{1}}^{\phi_1,\phi_2}$ and the other follows the symmetry considerations. 
\begin{remark}
    Our analysis of this section takes into account that $n\ge 2$ and that $m=-(n-1).$ Note that this is the bilinear counterpart of the Seeger-Sogge-Stein  order $-\frac{n-1}{2}$ in the linear setting. 
\end{remark}

\subsection{The estimate of $L^2\times L^\infty \to L^2$}
Let $f, h\in L^2$ and $g\in L^\infty$. Using the representation in $I_4$ and applying Cauchy-Schwarz inequality yields
$$
|I_4|\lesssim\left(
\iint_{0}^{1}\left|Q_t^{u} T_{\nu{m_1}}^{\phi_1} f(x)\right|^2\,\frac{dt}{t}\,dx\right)^{\frac12}\left(
\iint_{0}^{1}\left|P_t^{2,v} R_t T_{\gamma}^{\phi_2} g(x)\right|^2\left|\widetilde{Q}_t h(x)\right|^2\,\frac{dt}{t}\,dx\right)^{\frac12}.
$$
Using the fact that, $\|T_{a_1}^{\phi_1}f\|_{L^2}\le C\|f\|_{L^2}$ when $1<p<\infty$ and 
$$\left\| P_t^{2,v} R_t T_{\gamma}^{\phi_2} g \right\|_{L^\infty}\lesssim
\left\| R_t T_{\gamma}^{\phi_2} g\right\|_{L^\infty}
\lesssim\left\| T_{\gamma}^{\phi_2} g\right\|_{\mathrm{BMO}}\lesssim
\|g\|_{L^\infty}.$$
Clearly, 
\begin{align*}
\left(
\iint_{0}^{1}\left|P_t^{2,v} R_t T_{\gamma}^{\phi_2} g(x)\right|^2\left|\widetilde{Q}_t h(x)\right|^2\,\frac{dt}{t}\,dx\right)^{\frac12}\lesssim \|g\|_{L^\infty}\|h\|_{L^2}.
\end{align*}
Hence, $|I_4|\lesssim \|f\|_{L^2}\|g\|_{L^\infty}\|h\|_{L^2}.$
We again estimate $I_6$ using the standard quadratic estimate for Carleson measures. This gives us, 
\begin{align*}
|I_6|
\lesssim\left(\int_{\mathbb{R}^n}\int_0^1\left|
Q_t^{u} T_{\nu_{m_1}}^{\phi_1} f(x)\right|^2\,\frac{dt}{t}\,dx\right)^{\frac{1}{2}}\left(\int_{\mathbb{R}^n}\int_0^1\left|\widetilde{P}_t h(x)\right|^2\left|
Q_t^{2,v} R_t T_{\gamma}^{\phi_2} g(x)\right|^2\,\frac{dt}{t}\,dx\right)^{\frac{1}{2}}.
\end{align*}
The first factor is bounded using the same idea as $I_4$, our attention now is on the second factor, which involves the use of Carleson measure. Indeed, we know that,
$$ d\mu(x,t) =\left|
Q_t^{2,v} R_t T_{\gamma}^{\phi_2} g(x)\right|^2\,\frac{dt}{t}\,dx$$ is a Carleson measure with norm bounded by a constant times
$\|g\|_{L^\infty}^2 = 1.$ The previous analysis implies that 
$$\left(\int_{\mathbb{R}^n}\int_0^1\left|\widetilde{P}_t h(x)\right|^2\left|
Q_t^{2,v} R_t T_{\gamma}^{\phi_2} g(x)\right|^2\,\frac{dt}{t}\,dx\right)^{\frac{1}{2}} \lesssim \|g\|_{L^\infty}\|h\|_{L^2}.$$ Hence, 
$$|I_6|\lesssim\|f\|_{L^2}\|g\|_{L^\infty}\|h\|_{L^2}.$$
\subsection{The estimate of $L^\infty\times L^\infty \to \mathrm{BMO}$} In this subsection, we prove that the bilinear Fourier integral operator $T^{\phi_1,\phi_2}_a:L^\infty\times L^\infty \to \mathrm{BMO}$ is bounded whenever $m=-(n-1).$ To do this, we estimate both $I_4$ and $I_6$ such that $m_1=m_2=-\dfrac{(n-1)}{2}$ for $f,g \in L^{\infty}$ and $h\in H^1.$ We may assume without loss of generality that $\|f\|_{L^{\infty}}=\|g\|_{L^{\infty}}=1.$ We have
\begin{align*}
    I_4=\int_{\mathbb{R}^n} \int_0^1 
(Q^u_t T^{\phi_1,\tau}_{\nu_{m_1}} f)(x)\,
(P^{2,v}_t [R_t T^{\phi_2,\tau}_\gamma] g)(x)\,
m_{\tau}( x)\,(\widetilde{Q}_t h(\cdot))(x)\,
\frac{dt}{t}\,dx,
\end{align*}
taking the norm on both sides gives us
\begin{align*}
    |I_4|
\leq\left\|P_t^{v,2}R_tT_{\gamma}^{\phi_2}g\right\|_{L^\infty}\left|
\iint_0^1\left(Q_t^uT_{\nu_{m_1}}^{\phi_1}f(x)\right)\widetilde Q_t h(x)\,m_{\tau}(x)
\,\frac{dt}{t}\,dx\right|.
\end{align*}
Using the fact that $\left\|P_t^{2,v}R_tT_{\gamma}^{\phi_2}g\right\|_{L^\infty}
\lesssim
\left\|R_tT_{\gamma}^{\phi_2}g\right\|_{L^\infty}
\lesssim
\left\|T_{\gamma}^{\phi_2}g\right\|_{\mathrm{BMO}}
\lesssim
\|g\|_{L^\infty}.$ 
The above estimate is valid since $P^{2,v}_t$ is a Fourier multiplier which is bounded by $L^\infty$ norm. $ R_t$ is also bounded from  $L^{\infty}$ to $\mathrm{BMO},$ (see (75) of \cite[Page ~28]{rodriguez2011boundedness}). Again, from T.2 in Theorem~2.5 in \cite{rodriguez2011boundedness}, follows that $\left\|T_{\gamma}^{\phi_2}g\right\|_{\mathrm{BMO}}
\lesssim \|g\|_{L^\infty},$ since $T_{\gamma}^{\phi_2}$ is a Fourier integral operator.
One has that
\begin{align*}
    |I_4|
\leq \|g\|_{\infty}\left|
\iint_0^1\left(Q_t^uT_{\nu_{m_1}}^{\phi_1}f(x)\right)\widetilde Q_t h(x)\,m_{\tau}(x)
\,\frac{dt}{t}\,dx\right|.
\end{align*}
 Consequently, by applying Corollary~4.12 of \cite[Page ~35]{rodriguez2011boundedness} to the above estimates, we have that
\begin{align*}
    |I_4|
\lesssim
\left\|T_{\nu_{m_1}}^{\phi_1}f\right\|_{\mathrm{BMO}}
\|h\|_{H^1}\|g\|_{L^\infty}
\lesssim\|f\|_{L^\infty}\|h\|_{H^1}\|g\|_{L^\infty}.
\end{align*}
Next, we estimate $I_6,$ which is given as
\begin{align*}
   I_6= \int_{\mathbb{R}^n} \int_0^1 
(Q^u_t T^{\phi_{1,\tau}}_{\nu_{m_1}} f)(x)\,
(Q^{v,2}_t [R_t T^{\phi_{2,\tau}}_\gamma] g)(x)\,
m_{\tau}( x)\,(\widetilde{P}_t h(\cdot))(x)\,
\frac{dt}{t}\,dx, 
\end{align*}
taking norm on both sides gives us
\begin{align*}
|I_6|
&\leq\|m_\tau\|_{L^\infty}
\int_{\mathbb{R}^n}\int_0^1\bigl|
(Q_t^uT_{\nu_{m_1}}^{\phi_{1,\tau}}f)(x)\,
(Q_t^{v,2}[R_tT_\gamma^{\phi_{2,\tau}}]g)(x)\bigr|
\,|\widetilde P_th(x)|\frac{dt}{t}\,dx.
\end{align*}
Since,
\begin{align*}
\bigl|
(Q_t^uT_{\nu_{m_1}}^{\phi_{1,\tau}}f)(x)\,
(Q_t^{v,2}[R_tT_\gamma^{\phi_{2,\tau}}]g)(x)\bigr|
\leq\frac12\left(
|Q_t^uT_{\nu_{m_1}}^{\phi_{1,\tau}}f(x)|^2
+|Q_t^{v,2}[R_tT_\gamma^{\phi_{2,\tau}}]g(x)|^2
\right).
\end{align*}
We use the fact that,
$$
|Q_t^uT_{\nu_{m_1}}^{\phi_{1,\tau}}f(x)|^2
\,dx\,\frac{dt}{t}
\quad\text{and}\quad
|Q_t^{v,2}[R_tT_\gamma^{\phi_{2,\tau}}]g(x)|^2
\,dx\,\frac{dt}{t}
$$ are both Carleson measures, which implies that,
$$
d\mu(x,t)
=(Q_t^uT_{\nu_{m_1}}^{\phi_{1,\tau}}f)(x)\,
(Q_t^{v,2}[R_tT_\gamma^{\phi_{2,\tau}}]g)(x)\,
dx\,\frac{dt}{t}, $$
is also a Carleson measure satisfying the property that,
$$\|d\mu\|_{\mathcal C}
\lesssim\|f\|_{L^\infty}^2
+ \|g\|_{L^\infty}^2=2, \text{ since we assume previously that, }$$
$$ \|f\|_{L^{\infty}}=\|g\|_{L^{\infty}}=1.$$ 
Now, from \cite[Theorem~4]{fefferman1972h} and \eqref{reg:2}, we have
$$
|I_6|\lesssim\iint_0^1\left|\widetilde{P}_t h(x)\right|
\,d|\mu|(x,t)\lesssim\int_{\mathbb{R}^n}\sup_{0<t<1}
\sup_{|x-y|<t}\left|\widetilde{P}_t h(y)\right|
\,dx\lesssim\|h\|_{H^1}.
$$
In this case, the restriction $n\ge 2$ is necessary, here we state a counterexample when $n=1.$
\begin{proposition}\label{pro}
There exists $f,g \in L_c^\infty$ and
\(a(x,\xi,\eta) \in S^{0}_{1,0}(\mathbb{R}\times (\mathbb{R}\setminus\{0\})\times (\mathbb{R}\setminus\{0\}))\) such that
$T^{\phi_1, \phi_2}_a(f,g)
\notin \operatorname{BMO}$ where $\phi_{1}(x,\xi)=x\xi+|\xi| \text{ and }\phi_{2}(x,\eta)=x\eta.$
\end{proposition}
\begin{remark}
    The counterexample for the real case also serves as a counterexample in the complex case (see Proposition~5.1 of \cite[Page ~38]{rodriguez2011boundedness}).
\end{remark}
\subsection{The estimate of $ H^1 \times L^\infty \to L^1$} To estimate both $I_4$ and $I_6,$ we consider the case when $m_1=m_2=-\dfrac{(n-1)}{2}$ and for $f\in H^1, g,h\in L^{\infty}.$ Given that
\begin{align*}
I_4=\int_{\mathbb{R}^n}\int_0^1
\left(Q_t^u\left((1-\widehat{\theta})(D)T_{\nu_{m_1}}^{\phi_1}f(x)\right)\right)
\left(P_t^{2,v}R_tT_{\gamma}^{\phi_2}g(x)\right)
\widetilde Q_th(x)\,m_\tau(x)\,\frac{dt}{t}\,dx.
\end{align*}
Notably, it has been established in \cite{Goldberg1979} that a function $f$ belongs to the local Hardy space $h^1$ if and only if
$f \in L^1(\mathbb{R}^n)$  and $\mathfrak{R_j}\bigl((1-\widehat{\theta})(D)f\bigr)\in L^1(\mathbb{R}^n)$
where $\widehat{\theta}$ is a smooth bump function supported in a neighbourhood of the origin, and $\mathfrak{R_j}$ denotes the $j$-th Riesz transform, that is, \begin{align*}
\|f\|_{h^1}\approx\|f\|_{L^1}+\sum_{j=1}^{n}
\left\|\mathfrak{R_j}\bigl((1-\widehat{\theta})(D)f\bigr)\right\|_{L^1}
\approx\|\widehat{\theta}(D)f\|_{L^1}
+\|(1-\widehat{\theta})(D)f\|_{H^1},
\end{align*}
where the last estimate follows from the classical characterization of the $H^1$ norm in terms of the Riesz transforms (see, for example, \cite[Page~221]{Stein1970}). It is worth mentioning that different choices of $\theta$ yield equivalent norms. Since the frequency support of the symbol associated with $Q_t^u$ is supported in an annulus, there exists a smooth compactly supported bump function $\widehat{\theta}$ such that, for any $0<t<1,$ we have
$Q_t^u=Q_t^u(1-\widehat{\theta})(D).$ Due to this, we rewrite $I_4$ as
\begin{align*}
I_4=\int_{\mathbb{R}^n}\int_0^1
\left(Q_t^u(1-\widehat{\theta})(D)T_{\nu_{m_1}}^{\phi_1}f(x)\right)
\left(P_t^{2,v}R_tT_{\gamma}^{\phi_2}g(x)\right)
\widetilde Q_th(x)\,m_\tau(x)
\,\frac{dt}{t}\,dx.
\end{align*}
We apply Corollary~4.12 in \cite{rodriguez2011boundedness}, with $F=h$, $G=T_{\nu_{m_1}}^{\phi_1}f(x)$ and $$v(t,x)=P_t^{2,v}R_tT_{\gamma}^{\phi_2}g(x)\,m_\tau(x),$$ gives us
\begin{align*}
|I_4|\lesssim\|T_{\nu_{m_1}}^{\phi_1}f\|_{H^1}\|g\|_{L^\infty}
\|h\|_{\mathrm{BMO}}.
\end{align*}
For the natural embedding of
$L^\infty$ into $\mathrm{BMO},$ and using the characterization of the local Hardy space,
$\|T_{\nu_{m_1}}^{\phi_1}f\|_{H^1}
\lesssim \|T_{\nu_{m_1}}^{\phi_1}f\|_{h^1},$ one has that
\begin{align*}
|I_4|\lesssim\left\|T_{\nu_{m_1}}^{\phi_1}f\right\|_{H^1}\|g\|_{L^\infty}\|h\|_{\mathrm{BMO}}\lesssim\left\|T_{\nu_{m_1}}^{\phi_1}f\right\|_{h^1}\|g\|_{L^\infty}\|h\|_{L^\infty}.
\end{align*}
Again, using the embedding of $H^1$ into $h^1$, one obtains
\begin{align*}
|I_4| \lesssim \|f\|_{H^1}\,\|g\|_{L^\infty}\,\|h\|_{L^\infty}.
\end{align*}
To control $I_6,$ we observe that a similar analysis to that of $I_4$ allows us to write
\begin{align*}
I_6 =\int_{\mathbb{R}^n}\int_0^1 
\big(Q^u_{t}(1-\widehat{\theta})(D)T^{\phi_1}_{\nu m_1} f(x)\big)\big(Q^{2,v}_t R_t T^{\phi_2}_{\gamma} g(x)\big)
\tilde{P}_t h(x)\, m_\tau(x)\, dt\, dx.
\end{align*}
Applying Proposition~4.11 in \cite{rodriguez2011boundedness}, assume $$v(t,x)= \tilde{P}_t h(x)\, m_\tau(x), \, F= (1-\widehat{\theta})(D)T^{\phi_1}_{\nu m_1} f(x)$$ and $$G(t,x)= Q^{2,v}_t R_t T^{\phi_2}_{\gamma} g(x)$$ and using the fact that $
|G(t,x)|^{2}\,dx\,\frac{dt}{t}$
is a Carleson measure. Following a similar argument as in the case of $I_4,$ yields
\begin{align*}
|I_6|\lesssim
\big\|(1-\widehat{\theta})(D)T^{\phi_1}_{\nu m_1} f\big\|_{H^1}\,\|g\|_{L^\infty}\,\|h\|_{BMO}
\lesssim\|f\|_{H^1}\,\|g\|_{L^\infty}\,\|h\|_{L^\infty}.\end{align*}

\subsection{The estimate of $ L^{\infty} \times H^{1}\to L^1$} In this case, $g \in H^1$, $f, h \in L^\infty$, and $m_1=m_2=-\dfrac{(n-1)}{2}$. Without loss of generality, we assume that $\|f\|_{L^\infty} = \|h\|_{L^\infty} = 1$. Arguing just as in $  H^1 \times L^{\infty}\to L^1.$ Here, the smooth compactly supported bump function $\widehat{\theta}$ is acting on $R_t$ such that
\begin{align*}
I_6 =\int_{\mathbb{R}^n}\int_0^1 
\big(Q^u_{t} T^{\phi_1}_{\nu_{m_1}} f(x)\big)
\big(Q_t^{2,v} R_t (1-\widehat{\theta})(D) T^{\phi_2}_{\gamma} g(x)\big)
\widehat{P}_t h(x)\, m_\tau(x)\, \frac{dt}{t}\, dx,
\end{align*}
taking norm on both sides gives us
\begin{align*}
|I_6|\lesssim\left\| T^{\phi_1}_{\nu_{m_1}} f \right\|_{BMO}\left\| (1-\widehat{\theta})(D) T^{\phi_2}_{\nu_{m_1}} g \right\|_{H^1}\|h\|_{L^\infty}.
\end{align*}
Using the fact that, 
\begin{align*}
\left\| T^{\phi_1}_{\nu_{m_1}} f \right\|_{BMO}
\lesssim \|f\|_{L^\infty} \text{ and }
\left\| (1-\widehat{\theta})(D) T^{\phi_2}_{\nu_{m_1}} g \right\|_{H^1}
\lesssim \|g\|_{H^1}.
\end{align*}
So, we obtain
\begin{align*}
|I_6|\lesssim
\left\| T^{\phi_1}_{\nu_{m_1}} f \right\|_{BMO}
\left\| (1-\widehat{\theta})(D) T^{\phi_1}_{\nu_{m_1}} g \right\|_{H^1}
\|h\|_{L^\infty}
\lesssim\|f\|_{L^\infty}\,\|g\|_{H^1}\,\|h\|_{L^\infty}.\end{align*}
To control $I_4$, we follow the same argument as before, such that we rewrite $I_4$ as
\begin{align*}
I_4=\int_{\mathbb{R}^n}\int_0^1
\left(Q_t^uT_{\nu_{m_1}}^{\phi_1}f(x)\right)
\left(P_t^{2,v}R_t(1-\widehat{\theta})(D)T_{\gamma}^{\phi_2}g(x)\right)
\widetilde Q_th(x)\,m_\tau(x)
\,\frac{dt}{t}\,dx,
\end{align*}
Following a similar argument as in the case of $I_6$ and \cite[Theorem~4]{fefferman1972h}, and using the tangential maximal function characterization of the local Hardy space (see~(9) of~\cite[Page.~8]{rodriguez2011boundedness}) yields
\begin{align*}
|I_4|\lesssim\int_{\mathbb{R}^n}
\sup_{0<t<1}\ \sup_{|x-y|<t}\big|P_t^{2,v} T^{\phi_2}_{\gamma} g(x)\big|\, dx
\lesssim\|T^{\phi_2}_{\gamma} g\|_{h^1}
\lesssim\|g\|_{H^1}.
\end{align*}
\begin{remark}
    It is known that $H^p$ is a subspace of $h^p$ and, for any $f\in H^p, \,\|f\|_{h^p}\lesssim \|f\|_{H^p}$. In particular, when $p=1,$ we have that $\|f\|_{h^1}\lesssim \|f\|_{H^1}.$
\end{remark}

\section{Bilinear Fourier integral operators on the torus} \label{4}
In this section, we consider the boundedness of FIO's on the torus $\mathbb{T}^n$. For this, we consider an extension result similar to the one proved by Cardona \cite{Cardona:2026}. First, we state an auxiliary lemma.

\label{section:main}\label{ssection:extension}
\subsection{Transference theorems} We start this section by establishing some transference theorems about the boundedness of Fourier integral operators from $\mathbb{R}^n$ to the torus $\mathbb{T}^n.$ The following lemma will be required for our further analysis. 

\begin{lemma}
    Suppose $f$ is a continuous periodic function on $\mathbb{R}^n$ and let $\{g_m\}_m$ be a sequence of uniformly bounded periodic continuous functions on $\mathbb{R}^n$. If $\{g_m\}_m$ converges point-wise to a function $g$ defined on $\mathbb{R}^n$ and $\{\varepsilon_m\}_m$ is a sequence of positive real numbers, then 
    \begin{equation*}
        \lim_{m\to\infty} \varepsilon_m^{n/2} \int_{\mathbb{R}^n} e^{-\varepsilon|x|^2} f(x)g_m(x)\diff x = \int_{\mathbb{T}^n} f(x)g(x)\diff x,
    \end{equation*}
    provided that $\varepsilon_m \to 0$.
    \label{lem:convergence}
\end{lemma}
Next, we establish a criterion for the transference of boundedness of Fourier integral operators with complex phases. Since we are interested in obtaining boundedness properties on Lebesgue spaces, the condition $\mathrm{Im}(\phi) \geq M,$ on a phase function for some real number $M,$ will be reasonably required. 
\begin{theorem}\label{torus}
    Let $1<p<\infty$, and let $\phi$ a complex valued function defined in $\mathbb{T}^n\times \mathbb{R}^{nr}$ such that $\mathrm{Im}(\phi) \geq M$ where $M\in \mathbb{R}.$ Let us assume that $a:\mathbb{T}^n\times \mathbb{R}^{nr}\to \mathbb{C}$ is bounded and continuous, and the corresponding multi-linear Fourier integral operator
    \begin{equation*}
        Tf(x) = \int_{\mathbb{R}^{nr}} e^{i\phi(x, \xi)}a(x, \xi)\widehat{f_1}(\xi_1)\cdots \widehat{f_r}(\xi_r) \diff \xi
    \end{equation*}
    extends to a bounded linear operator from $L^{p_1}(\mathbb{R}^n)\times\cdots\times L^{p_r}(\mathbb{R}^n)$ into $L^{p}(\mathbb{R}^n)$. Then the periodic multi-linear Fourier integral operator 
    \begin{equation*}
        Af(x) = \sum_{\xi \in\mathbb{Z}^{nr}} e^{i\phi(x, \xi)}a(x, \xi)(\mathcal{F}_{\mathbb{T}^n}f_1)(\xi_1)\cdots (\mathcal{F}_{\mathbb{T}^n}f_r)(\xi_r) 
    \end{equation*}
    extends to a bounded linear operator from $L^{p_1}(\mathbb{T}^n)\times\cdots\times L^{p_r}(\mathbb{T}^n)$ into $L^{p}(\mathbb{T}^n)$, given that
    \begin{equation*}
        \frac{1}{p_1} + \cdots +\frac{1}{p_r} = \frac{1}{p}, \quad 1\leq p_i <\infty.
    \end{equation*}
    Moreover, there exists a positive constant $C_p$ such that 
    \begin{equation*}
        \|A\|_{\mathcal{B}(L^{p_1}(\mathbb{T}^n)\times\cdots\times L^{p_r}(\mathbb{T}^n), L^{p}(\mathbb{T}^n))} \leq C_p \|T\|_{\mathcal{B}(L^{p_1}(\mathbb{R}^n)\times\cdots\times L^{p_r}(\mathbb{R}^n), L^{p}(\mathbb{R}^n))}
    \end{equation*}
\end{theorem}

\begin{proof}
    First, we prove the following identity for trigonometric polynomials $P_j$ and $Q$ on $\mathbb{T}^n$: 
    \begin{align}
        &\lim_{\varepsilon\to0} \varepsilon^{n/2}\int_{\mathbb{R}^n} [T(P_1w_{\alpha_1\varepsilon}, \ldots, P_rw_{\alpha_r\varepsilon})](x)\overline{Q(x)}w_{\beta\varepsilon}\diff x  \nonumber\\&= c_{n,r,p}\int_{\mathbb{T}^n} [A(P_1,\ldots,P_r)](x)\overline{Q(x)}\diff x,
        \label{eq:polynomial-identity}
    \end{align}
    for some positive constant $c_{n,r,p}$, where $w_\delta(x) = e^{-\delta|x|^2}$ and $\sum a_j + \beta = 1.$
    By linearity, we may assume that $P_j = e^{i2\pi m_jx_j}$ and $Q=e^{i2\pi kx}$ for $k,m_j\in\mathbb{Z}^n$. Notice that 
    \begin{align*}
        \int_{\mathbb{T}^n} [A(P_1,\ldots,P_r)](x)\overline{Q(x)}\diff x & = \int_{\mathbb{T}^n} \sum_{\xi\in\mathbb{Z}^{nr}} e^{i\phi(x, \xi)} a(x, \xi) \delta_{m,\xi}\overline{Q(x)}\diff x \\
        & = \int_{\mathbb{T}^n} e^{i\phi(x, m)} a(x, m)\overline{Q(x)} \diff x\\
        & = \int_{\mathbb{T}^n} e^{i\phi(x, m)-i2\pi kx} a(x, m) \diff x.
    \end{align*}
    Now, recall that the Euclidean Fourier transform of $P_j w_{\alpha_j\varepsilon}$ is given by 
    \begin{equation*}
        (\mathcal{F}_{\mathbb{R}^n}P_j w_{\alpha_j\varepsilon})(\xi_j) = (\alpha_j\varepsilon)^{-n/2}e^{-|\xi_j-m_j|^2/\alpha_j\varepsilon}.
    \end{equation*}
    Hence, we obtain that 
    \begin{align*}
        &\varepsilon^{n/2}\int_{\mathbb{R}^n} [T(P_1w_{\alpha_1\varepsilon}, \ldots, P_rw_{\alpha_r\varepsilon})](x)\overline{Q(x)}w_{\beta\varepsilon}\diff x \\ 
        & = \varepsilon^{n/2}\int_{\mathbb{R}^n} \int_{\mathbb{R}^{nr}} e^{i\phi(x, \xi)}a(x, \xi)\prod(\alpha_j\varepsilon)^{-n/2}e^{-|\xi_j-m_j|^2/\alpha_j\varepsilon}\overline{Q(x)}w_{\beta\varepsilon}\diff \xi\diff x \\ 
        & = \varepsilon^{n/2}\int_{\mathbb{R}^n} \int_{\mathbb{R}^{nr}} e^{i\phi(x, \xi)-i2\pi kx} e^{-\beta\varepsilon|x|^2}a(x, \xi)\prod(\alpha_j\varepsilon)^{-n/2}e^{-|\xi_j-m_j|^2/\alpha_j\varepsilon}\diff \xi\diff x\\
        & = \varepsilon^{n/2}\int_{\mathbb{R}^{nr}} \int_{\mathbb{R}^{n}} e^{i\phi(x, \xi)-i2\pi kx} e^{-\beta\varepsilon|x|^2}a(x, \xi)\diff x \prod(\alpha_j\varepsilon)^{-n/2}e^{-|\xi_j-m_j|^2/\alpha_j\varepsilon}\diff \xi \\
        & = \beta^{-n/2}(\beta\varepsilon)^{n/2}\int_{\mathbb{R}^{nr}} \int_{\mathbb{R}^{n}} e^{i\phi(x, (\alpha\varepsilon)^{1/2} \eta + m)-i2\pi kx} e^{-\beta\varepsilon|x|^2}a(x, (\alpha\varepsilon)^{1/2} \eta + m)\diff x e^{-|\eta|^2}\diff \eta
    \end{align*}
    where $(\alpha\varepsilon)^{1/2} \eta + m = (\alpha_1\varepsilon)^{1/2} \eta_1 + m_1 + \cdots + (\alpha_r\varepsilon)^{1/2} \eta_r + m_r$. Now, notice that  
    \begin{equation*}
        |e^{i\phi(x, (\alpha\varepsilon)^{1/2} \eta + m)-i2\pi kx} e^{-\beta\varepsilon|x|^2}a(x, (\alpha\varepsilon)^{1/2} \eta + m)|  \leq e^{-M}\|a\|_{L^\infty} 
    \end{equation*}
    for every $x\in\mathbb{R}^n$ and $\varepsilon>0$. Hence, by Lemma \ref{lem:convergence}, we have that 
    \begin{align*}
        &\lim_{\varepsilon\to0}\varepsilon^{n/2}\int_{\mathbb{R}^n} [T(P_1w_{\alpha_1\varepsilon}, \ldots, P_rw_{\alpha_r\varepsilon})](x)\overline{Q(x)}w_{\beta\varepsilon}\diff x \\ 
        & =  \beta^{-n/2}\int_{\mathbb{R}^{nr}} \lim_{\varepsilon\to0}(\beta\varepsilon)^{n/2}\int_{\mathbb{R}^{n}} e^{i\phi(x, (\alpha\varepsilon)^{1/2} \eta + m)-i2\pi kx} e^{-\beta\varepsilon|x|^2}\\& \quad\times a(x, (\alpha\varepsilon)^{1/2} \eta + m)\diff x e^{-|\eta|^2}\diff \eta \\
        & =  \beta^{-n/2}\int_{\mathbb{R}^{nr}} \int_{\mathbb{T}^{n}} e^{i\phi(x, m)-i2\pi kx} a(x, m)\diff x e^{-|\eta|^2}\diff \eta\\
        & =  \beta^{-n/2}\pi^{nr/2} \int_{\mathbb{T}^{n}} e^{i\phi(x, m)-i2\pi kx} a(x, m)\diff x. 
    \end{align*}
    We can choose $\alpha_j = 1/p_j$ and $\beta=1/p'$ so that with $C_p =(\pi^r p')^{-n/2}$ we have that 
    \begin{align*}
        &\|A(P_1, \ldots, P_r)\|_{L^{p}(\mathbb{T}^n)} \\
        & = \sup_{\|Q\|_{L^{p'}(\mathbb{T}^n)}=1} \left|\int_{\mathbb{T}^n} [A(P_1,\ldots,P_r)](x)\overline{Q(x)}\diff x\right| \\
        & = C_p  \sup_{\|Q\|_{L^{p'}(\mathbb{T}^n)}=1} \lim_{\varepsilon\to0} \varepsilon^{n/2} \left|\int_{\mathbb{R}^n} [T(P_1w_{\alpha_1\varepsilon}, \ldots, P_rw_{\alpha_r\varepsilon})](x)\overline{Q(x)}w_{\beta\varepsilon}\diff x\right| \\
        & \leq C_p\sup_{\|Q\|_{L^{p'}(\mathbb{T}^n)}=1} \lim_{\varepsilon\to0} \varepsilon^{n/2} \|T\|_{\mathcal{B}(L^{p_1}\times\cdots\times L^{p_r}, L^{p})} \prod \|P_jw_{\varepsilon/p_j}\|_{L^{p_j}(\mathbb{R}^n)} \|Qw_{\varepsilon/p'}\|_{L^{p'}(\mathbb{R}^n)}      \\
        & =C_p \sup_{\|Q\|_{L^{p'}(\mathbb{T}^n)}=1} \|T\|_{\mathcal{B}}\lim_{\varepsilon\to0}   \prod \left(\varepsilon^{n/2}\int_{\mathbb{R}^n}|P_jw_{\varepsilon/p_j}|^{p_j}\diff x\right)^{1/p_j} \\& \quad \times \left(\varepsilon^{n/2}\int_{\mathbb{R}^n}|Qw_{\varepsilon/p'}|^{p'}\diff x\right)^{1/p'} \\
        & =C_p \sup_{\|Q\|_{L^{p'}(\mathbb{T}^n)}=1} \|T\|_{\mathcal{B}}\lim_{\varepsilon\to0}   \prod \left(\varepsilon^{n/2}\int_{\mathbb{R}^n}|P_j|^{p_j}e^{-\varepsilon|x|^2}\diff x\right)^{1/p_j} \\& \quad\times\left(\varepsilon^{n/2}\int_{\mathbb{R}^n}|Q|^{p'}e^{-\varepsilon|x|^2}\diff x\right)^{1/p'}      \\
        & =C_p \sup_{\|Q\|_{L^{p'}(\mathbb{T}^n)}=1} \|T\|_{\mathcal{B}}  \prod \left(\int_{\mathbb{T}^n}|P_j|^{p_j}\diff x\right)^{1/p_j} \left(\int_{\mathbb{T}^n}|Q|^{p'}\diff x\right)^{1/p'} \\
        & = C_p \|T\|_{\mathcal{B}(L^{p_1}\times\cdots\times L^{p_r}, L^{p})}  \prod \|P_j\|_{L^{p_j}(\mathbb{T}^n)}.
    \end{align*}
    Thus, there exists a unique bounded multi-linear extension from $L^{p_1}(\mathbb{T}^n)\times\cdots\times L^{p_r}(\mathbb{T}^n)$ into $L^{p}(\mathbb{T}^n)$.
\end{proof}

\subsection{Boundedness properties on the torus}
Combining the result in the previous subsection with Theorem {\color{blue} \ref{torus}}, we obtain the following theorem for bilinear FIO's with complex phases on the torus.

\begin{corollary}
    Let $T^{\phi_1,\phi_2}_a$ be a bilinear Fourier integral operator defined on the torus as 
    \begin{equation*}
        T^{\phi_1,\phi_2}_a (f,g)(x) = \sum_{\xi,\eta\in\mathbb{Z}^n} e^{2\pi i[\phi_1(x, \xi)+\phi_2(x,\eta)]} a(x, \xi,\eta) (\mathcal{F}_{\mathbb{T}^n}f)(\xi)(\mathcal{F}_{\mathbb{T}^n}g)(\eta),
    \end{equation*}
    where $a\in S^m_{1,0}(\mathbb{T}^n_x \times \mathbb{R}^n_{\xi}\times \mathbb{R}^n_{\eta})$ is compactly supported in $x,$
    and there exists $\tau\in \mathbb{R},$ such that $\phi_{1,\tau}(x,\xi) = \Re \phi_1(x,\xi) + \tau \Im \phi_1(x,\xi) $ and $\phi_{2,\tau}(x,\eta) = \Re \phi_2(x,\eta) + \tau \Im \phi_2(x,\eta)$ satisfy:
\begin{itemize}
        \item[1.] $\phi_{1,\tau}(x,\xi)$ is positively homogeneous of degree one in $\xi,$ and $\phi_{2,\tau}(x,\eta)$ is also positively homogeneous of degree one in $\eta,$
        \item[2.]  $\phi_{1,\tau}(x,\xi) $ and $\phi_{2,\tau}(x,\eta)$ satisfy the non-degeneracy condition.
        \item[3.] $\phi_1, \phi_2\in C^\infty( \mathbb{T}^n\times \mathbb{R}^n\setminus\{0\},\mathbb{C}).$
        \item[4.]We assume that the phases $\phi_j,$ with $\textnormal{Im}(\phi_1), \textnormal{Im}(\phi_2)\geq 0,$ are of {\it strictly positive type}.
    \end{itemize} 
    Such that, when
    \begin{equation*}
        1< p_1,p_2 < \infty, \quad 1<p< \infty, \quad \frac{1}{p}=\frac{1}{p_1} + \frac{1}{p_2},
    \end{equation*}
    and that 
    \begin{equation*}
        m \le - (n-1)\left( \left| \frac{1}{p_1}-\frac{1}{2} \right| + \left| \frac{1}{p_2}-\frac{1}{2} \right|
        \right), \quad n\ge 2.
    \end{equation*}
    Then $T^{\phi_1,\phi_2}_a$ extends to a bounded bilinear operator from $L^{p_1}(\mathbb{T}^n) \times L^{p_2}(\mathbb{T}^n)$ to $L^{p}(\mathbb{T}^n)$. Moreover, we have that 
    \begin{equation*}
        \| T^{\phi_1,\phi_2}_a(f,g)\|_{L^p} \lesssim \|f\|_{L^{p_1}}\|g\|_{L^{p_2}}.
    \end{equation*}
\end{corollary}

\subsection{Conclusions and further references}
The study of boundedness properties for Fourier integral operators and their bilinear counterparts has been extensively developed in recent years, we refer to \cite{Cardona:2026} for the smoothing conjecture for bilinear Fourier integral operators as well as for a discussion about different results. We also refer to Grafakos and Peloso \cite{grafakos2010bilinear} for the motivation of considering bilinear Fourier integral operators as operators for dealing with restriction estimates for wave-type equations. 

In the periodic setting, fundamental weak type (1,1) and H\"older--Besov bounds, as well as general \(L^p\)-boundedness criteria for periodic pseudo-differential operators, were established in \cite{Cardona:2014, Cardona:2016, Cardona:2018}. These results were substantially extended to the broader framework of graded Lie groups, where \(L^p\)-bounds for pseudo-differential operators were proved in \cite{Cardona:Delgado:Ruzhansky:2021}, and later \(L^p\)--\(L^q\) estimates were obtained in \cite{Cardona:Kumar:Ruzhansky:2025}. 

On compact Lie groups, \(L^p\)--\(L^q\) boundedness results were derived in \cite{Cardona:Delgado:Kumar:Ruzhansky:2025}, while refined kernel estimates and weak (1,1) bounds were provided in \cite{Cardona:Yeghoyan:Ruzhansky:2026}. The multilinear theory has also seen significant progress, with \(L^p\)-boundedness and nuclearity of multilinear pseudo-differential operators on \(\mathbb{Z}^n\) and the torus being investigated in \cite{Cardona:Kumar:2019}. 

A systematic review of the boundedness of toroidal/periodic pseudo-differential operators has been undertaken in a series of works \cite{Cardona:Martinez:2025a, Cardona:Martinez:2025b, Cardona:Martinez:2025c}, while periodic Fourier integral operators in \(L^p\)-spaces were analyzed in \cite{Cardona:Messiouene:Senoussaoui:2021}. Oscillatory singular integrals on compact Lie groups were studied in \cite{Cardona:Ruzhansky:2022}, and the subelliptic theory for pseudo-differential and Fourier integral operators on compact Lie groups has been comprehensively developed in the monograph \cite{Cardona:Ruzhansky:2026}. 

Despite this extensive body of work, the specific problem of Seeger-Sogge-Stein type orders for bilinear Fourier integral operators with complex phases remains largely unexplored, particularly, in the setting of the $L^{\frac{1}{2}}\times L^{\frac{1}{2}}\rightarrow L^{1,\infty} $-estimate. 

It is important to mention that the classical boundedness results for Fourier integral operators was established by Seeger--Sogge--Stein \cite{seeger1991regularity}. In particular, when $p=2$ is due to L.~H\"ormander \cite{hormander1971fourier}. M.~M.~Peloso and S.~Secco in \cite{peloso2008boundedness} proved the boundedness of Fourier integral operators on Hardy spaces, for $0<p\le 1.$ 

L.~Grafakos and M.~M,~Peloso, in \cite{grafakos2010bilinear} studied bilinear Fourier integral operators by extending known results concerning bilinear pseudodifferential operators. In \cite{rodriguez2015boundedness}, they extend the results of R.~Coifman and Y.~Meyer for bilinear pseudodifferential operators to the case of oscillatory integral operators. Some endpoint estimates for bilinear paraproducts on the local BMO space were established in \cite{rodriguez2015some}, and one of the goals here is to provide further extensions in that direction.

Recently, the work on Fourier integral operators with complex phases has been studied by the first author and Ruzhansky, in \cite{cardona2026weak}, establishing the weak $(1,1)$ boundedness of Fourier integral operators with complex phases. This result completes the endpoint estimates previously established by Ruzhansky \cite{Ruzhansky2001}. Further development was also made by the first and the third author in \cite{cardona2025boundedness}, which seeks to extend Fourier integral operators with real-valued phases on Fourier Lebesgue spaces to the complex-valued setting. 

Despite these advances and others, applying the problem of Seeger--Sogge--Stein type orders for bilinear Fourier integral operators on Besov spaces and Triebel--Lizorkin spaces with complex phases is missing, and using similar techniques developed in this paper to estimate the endpoints could help investigate the mapping properties for bilinear Fourier integral operators on these spaces. 

The present paper aims to fill this gap by establishing sharp boundedness criteria in the spirit of the classical Seeger-Sogge-Stein theorem for the bilinear setting with complex phases. For local smoothing properties for bilinear Fourier integral operators we refer the reader to \cite{Cardona:2026}.\\

\noindent{\bf Acknowledgement.}
 Duv\'an Cardona has been  supported  by the IR Center of Industrial Nuclear Energy and by the Department of Mathematics at the King Fahd University of Petroleum and Minerals. He also was partially supported by the  Odysseus  1  grant  G.0H94.18N:  Analysis  and  Partial Differential Equations of the Flemish Foundation FWO,  by the Methusalem programme of the Ghent University Special Research Fund (BOF)
(Grant number 01M01021), by the FWO Fellowship
Grant No 1204824N and by the FWO Grant K183725N of the Belgian Research Foundation FWO. He also has been supported by the Oberwolfach Leibniz Fellow of the Mathematical Institute of Research of Oberwolfach, MFO-Germany, Project F2511 (2026).  Manuel Mart\'inez-Flores has been supported by Liderazgo en Ciencias scholarship from Universidad del Valle de Guatemala. This work was partially supported by the Friends of IMU (FIMU), under the IMU Breakout Fellowship Program. Frederick Opoku has been supported by the Breakout Graduate Fellowship of the International Mathematical Union (IMU).
He has also been supported by the Kwame Nkrumah University of Science and Technology (KNUST), Ghana.

\subsubsection{Conflict of interests statement - Data Availability Statements}  The authors state that there is no conflict of interest.  Data sharing does not apply to this article as no datasets were generated or
analysed during the current study.

\bibliographystyle{amsplain}

\end{document}